%% file: cmqv.tex
\documentclass [12pt]{amsart}

\usepackage{amssymb,amsxtra,amsfonts}
\usepackage{epsf}              %
\usepackage{epsfig}             %
\usepackage{graphics}

\usepackage[colorlinks]{hyperref}

\input{macros}

\input{macros-thm1.1}

\begin{document}
\renewenvironment{ppb}[1]{}{}

\title[Global Weyl groups and multiplicative quiver varieties]{Global Weyl groups and a new theory of multiplicative quiver varieties}
\author{Philip Boalch}
\thanks{\tiny Supported by ANR grants 
13-IS01-0001-01, 
13-BS01-0001-01,  
09-JCJC-0102-01, 
08-BLAN-0317-01/02}

\maketitle

\begin{abstract}
In previous work a relation between a large class of
Kac--Moody algebras and meromorphic connections on global curves was established---notably the Weyl group gives isomorphisms between different  moduli spaces of connections, and the root system is also seen to play a role.
This involved a modular interpretation of
many Nakajima quiver varieties, as moduli spaces of connections, whenever the underlying graph was a complete k-partite graph (or more generally a supernova graph).
However in the isomonodromy story, or wild nonabelian Hodge theory, slightly larger moduli spaces of connections are considered.
This raises the question  of whether the full moduli spaces admit Weyl group isomorphisms, rather than just the open parts isomorphic to quiver varieties.
This question will be solved here, by developing a 
``multiplicative version'' 
of the previous approach.
This amounts to constructing many algebraic symplectic isomorphisms between wild character varieties. 
This approach also enables us to state a conjecture for certain irregular Deligne--Simpson problems
and introduce some noncommutative algebras, generalising the ``generalised double affine Hecke algebras''.
\end{abstract}

\section{Introduction}

The Nakajima quiver varieties \cite{nakaj-duke94, nakaj-duke98} are a large class of \hk manifolds attached to graphs and some data on the graph.
Some early examples arose \cite{Kron.ale} as ALE \hk four manifolds and 
\cite{kron-nakaj90} as a way to describe moduli spaces of instantons, i.e. solutions of the anti-self-dual Yang--Mills equations on the ALE spaces, generalising the Atiyah--Drinfeld--Hitchin--Manin construction of instantons \cite{adhm}.
The underlying graph  may be interpreted as a Dynkin diagram for a (symmetric) Kac--Moody algebra, and the quiver varieties attached to the graph have proven to be important in the study of Kac--Moody algebras (see e.g. \cite{nakaj-sugaku}).
In previous work \cite{rsode, slims} it was shown that some special Nakajima quiver varieties
arise as moduli spaces $\cM^*$ of meromorphic connections on the trivial bundle on the Riemann sphere, for example 
whenever the quiver is a complete $k$-partite graph for any $k$.
(This modular interpretation is different to that of \cite{Kron.ale, kron-nakaj90}, which only involves certain quivers closely related to affine Dynkin quivers.)
In brief, work of Kraft--Procesi \cite{Kraft-Procesi-InvMath79}, 
Nakajima \cite{nakaj-duke94} and Crawley-Boevey \cite{CB-additiveDS} 
showed that 
moduli spaces of Fuchsian systems
$$\cM^* \ \cong\  \cO_1\times \cdots \times \cO_m \spq G$$
were quiver varieties for star-shaped graphs (where $\cO_i\subset \Lie(G)^*, G=\GL_n(\IC)$)
and this was extended in \cite{rsode, slims} to some of the more general moduli spaces $\cM^*$ of \cite{smid} involving irregular connections.
This quiver approach to meromorphic connections is useful since it enables the construction of many isomorphisms between different moduli spaces 
$\cM^*$, often with different ranks and pole divisors, simply by reordering the symplectic quotient.
It also hints at where one might find non-affine Kac--Moody algebras in the geometry of connections on global curves (complementing the standard interpretation of affine Kac--Moody algebras in terms of loop algebras
\cite{pr-segal} p.77).

However these moduli spaces $\cM^*$ are open pieces of the full moduli spaces $\cM$ of connections on curves that appear in the isomonodromy story 
(\cite{smid} Rmk 2.1, Cor. 4.9) or wild nonabelian Hodge theory 
(\cite{wnabh, ihptalk}). These are also spaces of solutions to the anti-self-dual Yang--Mills equations, but different to those mentioned above---in brief  Hitchin \cite{Hit-sde} considered instantons that were invariant under two translations and found that the resulting `reduced' 
anti-self-dual Yang--Mills equations made sense on any compact Riemann surface, where they became equations, {\em the Hitchin equations}, for a pair consisting of a unitary connection together with a Higgs field.
As shown in \cite{Hit-sde}, moduli spaces of solutions to the Hitchin equations are also \hk  manifolds.
It has gradually been understood (\cite{Sim-hboncc, wnabh}) 
how to extend Hitchin's work to the case of noncompact Riemann surfaces, introducing boundary conditions so as to still obtain finite dimensional moduli spaces with (complete) \hk metrics (cf. the short survey \cite{ihptalk}). 
In brief this gives a way to attach a \hk manifold to a Riemann surface, and some other data on the surface (e.g. to specify the boundary conditions). We will call these \hk manifolds the {\em wild Hitchin moduli spaces}.
In one complex structure in the \hk family they are spaces of meromorphic Higgs bundles and in another they are spaces of meromorphic connections.
As one might expect, since Riemann surfaces are more complicated than 
graphs, the wild Hitchin moduli spaces are more complicated than the Nakajima quiver varieties.
More precisely there is no known finite dimensional construction of the metrics on the wild Hitchin spaces, whereas all the quiver varieties arise as \hk quotients of a finite dimensional \hk vector space.

The aim of the present article is to show that nonetheless a large class of the wild Hitchin spaces may be considered as precise ``multiplicative'' versions of quiver varieties, and then to use this to construct isomorphisms between different wild Hitchin spaces.
In particular a theory of ``multiplicative quiver varieties'' 
will be developed,
which will attach an algebraic symplectic manifold to a graph and some data on the graph. 
Then we will show
for a special class of graphs, the ``supernova graphs'' (including all the complete $k$-partite graphs for any $k$), that the multiplicative quiver variety attached to the graph is isomorphic to a wild Hitchin moduli space, thereby giving a {\em graphical} way to recover the underlying holomorphic symplectic manifold.
Conversely this shows that such multiplicative quiver varieties have (complete) \hk metrics.
A notion of multiplicative quiver varieties is also suggested by the work \cite{CB-Shaw}, but beyond the case of star-shaped graphs their approach does not give the desired spaces, as we will show explicitly by considering the case when the graph is a triangle 
(\S\ref{ssn: example mqv}).
In particular it is unclear if the spaces suggested by  \cite{CB-Shaw} are  \hk beyond the star-shaped case.

One of the main motivations for developing this approach is to better understand the automorphisms and isomorphisms between the 
wild character varieties (i.e. between the wild Hitchin spaces in their Betti algebraic structures).
As an application of the graphical approach we will show that the Weyl group of the Kac--Moody algebra attached to the graph acts to give algebraic symplectic isomorphisms between wild character varieties.
This is a multiplicative analogue of the approach of \cite{rsode,slims}; in essence the isomorphisms now arise by reordering the multiplicative 
symplectic quotients.
The fission idea (\cite{fission, gbs}) gives an inductive approach to prove that such isomorphisms are symplectic, enabling a reduction to the simplest case. 
Many (not necessarily affine) Kac--Moody Weyl groups thus arise in the context of the global geometry of a Riemann surface, in contrast to 
the familiar local understanding of affine Weyl groups related to loop groups (i.e. to gauge theory on a boundary circle).

The results of this paper hinge on certain developments in the theory of 
quasi-Hamiltonian geometry (a.k.a. multiplicative Hamiltonian  geometry or the theory of group valued moment maps)
and may be viewed as an exploration of a  class of examples. 
This theory was introduced in 1997 \cite{AMM} as an alternative finite dimensional approach to construct symplectic moduli spaces of flat connections on principal $K$-bundles over Riemann surfaces for compact groups $K$.
Subsequently the analogous holomorphic/algebraic theory for complex reductive groups $G$ such as $\GL_n(\IC)$ was developed (starting in 2002 \cite{saqh02}) and many new examples of quasi-Hamiltonian spaces were constructed in this context. 
These examples give great flexibility to construct many new algebraic symplectic manifolds (beyond spaces of representations of fundamental groups of Riemann surfaces) and the present article explores some of them. 

Specifically this extension made contact with the spaces of algebraic connections considered both in algebraic geometry and in the classical theory of differential equations, rather than just the holomorphic or $C^\infty$ connections familiar in differential geometry.
The key point is that certain connections on algebraic $G$-bundles on smooth complex algebraic curves may be classified by topological objects, {\em Stokes $G$-local systems}.
In turn the moduli spaces of Stokes local systems (the wild character varieties) have a quite explicit description, that {\em looks like} a multiplicative symplectic quotient.
Thus by extending the theory of multiplicative symplectic quotients to the complex algebraic world, and then constructing new examples, it is possible to prove that the wild character varieties {\em are} multiplicative symplectic quotients (\cite{saqh02, saqh, fission, gbs}). This gives an algebro-geometric alternative to the earlier analytic construction of such symplectic manifolds, 
{\em \`a la Atiyah--Bott}, in \cite{thesis, smid}, from which the quasi-Hamiltonian structures were derived (\cite{saqh02, saqh} \S4).
The wild nonabelian Hodge correspondence of \cite{Sab99, wnabh} implies that, in one complex structure in the \hk family, the wild Hitchin moduli spaces are spaces of algebraic connections on curves (cf. also the survey \cite{ihptalk}), and are thus isomorphic to wild character varieties. 

The further step we are taking in the present article is to notice that in some cases the wild character varieties {\em look like} multiplicative quiver varieties, and then develop the requisite theory so that they {\em are}  multiplicative quiver varieties.
As we will show this enables the construction of many (surprising) algebraic symplectic isomorphisms between wild character varieties.

As a consequence of this link between graphs and spaces of connections it is possible to 
make a conjecture (\S\ref{sn: gdsp}) about when some of the moduli spaces connections considered in 
\cite{smid, wnabh} are nonempty.
This conjecture involves the Kac--Moody root system attached to the graph, and generalises one of Crawley--Boevey \cite{CB-ihes} to some cases of the irregular Deligne--Simpson problem.
Further we will define some new noncommutative algebras (\S\ref{sn: fissionalgs})
which control the  multiplicative quiver varieties, generalising the 
multiplicative preprojective algebras of \cite{CB-Shaw} (some special cases of which 
contain the generalised double affine Hecke algebras of Etingof--Oblomkov--Rains \cite{eor-gda}).

\subsection{Statement of main results}

The first aim is to show that certain wild character varieties are multiplicative analogues of quiver varieties.
Since it is little extra trouble we will set-up the theory of multiplicative quiver varieties in full generality, even though it is only a special class
that appear as wild character varieties. (This is analogous to \cite{CB-Shaw} where only the star-shaped case is used.)
Thus the first step is
the construction of the multiplicative quiver varieties. 
In \S\ref{defn: coloured quiver} we will define the notion of a coloured quiver, which is a graph $\Ga$ with nodes $I$ and a colouring of the edges, together with some extra data involving certain orderings.
The colouring is a map $\ga:\Ga\to C$ from the set of edges of $\Ga$ to the set $C$ of colours.
In the simplest case the main property is that each monochromatic subgraph 
$\ga^{-1}(c)\subset \Ga$ consists of a disjoint union of complete $k$-partite graphs for any colour $c\in C$. 
The ``classical case'' is when  each monochromatic subgraph is a disjoint union of individual edges (and in this case the orderings amount to an orientation of each edge).

\begin{thm}
Suppose $\Ga$ is a coloured quiver with nodes $I$ and  we choose data $d\in \IZ^I_{\ge 0}$, and  $q\in (\IC^*)^I$ consisting of an integer and an invertible complex number for each node.
Then this data determines an algebraic variety, the ``multiplicative quiver variety'' $\cM(\Ga,q,d)$ together with a canonical open subset of stable points 
$$\cM^{st}(\Ga,q,d) \subset \cM(\Ga,q,d).$$
If nonempty $\cM^{st}(\Ga,q,d)$ is a smooth symplectic algebraic variety of dimension equal to $2-(d,d)$, where $(\,\,,\,\,)$ is the bilinear form of the Kac--Moody Cartan matrix of $\Ga$. 
Up to isomorphism the symplectic varieties
$\cM^{st}(\Ga,q,d),\cM(\Ga,q,d)$ only depend on $q,d$ and the underlying coloured graph and not on any of the choices of ordering.
\end{thm}

In the classical case (where each monochromatic subgraph consists of disjoint edges) this construction reduces to that suggested by \cite{CB-Shaw} and further studied in \cite{vdb-doublepoisson, yamakawa-mpa}.
We will also prove that if a graph is star-shaped then up to isomorphism it makes no difference how it is coloured (see Theorem \ref{thm: recolouring} for a more general statement).
In general the multiplicative quiver varieties {\em do} depend on the choice of colouring though:

\begin{prop}
Suppose $\Ga$ is a triangle. Then there are two inequivalent ways to colour $\Ga$ and the resulting multiplicative quiver varieties are not isomorphic in general, for example if $d=(1,1,1)$ and $q$ is generic.
\end{prop}

Further, beyond the star-shaped case, it is the ``nonclassical'' multiplicative quiver varieties that appear as wild character varieties, and thus have complete \hk metrics (via the irregular Riemann--Hilbert correspondence and the wild nonabelian Hodge theorem \cite{Sab99, wnabh}).
This depends on the consideration of a special class of graphs, 
the {\em supernova graphs} introduced in \cite{rsode, slims}: the simplest examples of such graphs 
consist of a central core which is a
complete $k$-partite graph,
together with a leg glued on to each node of the core 
(see Definition \ref{dfn: sn graph} in general).

\begin{thm}\label{thm: isom to wcv}
Suppose $\Ga$ is a simply-laced supernova graph, coloured so that its core is monochromatic, and $q,d$ are arbitrary.
Then the multiplicative quiver variety $\cM^{st}(\Ga,q,d)$
is isomorphic to a wild character variety, and is thus hyperk\"ahler. 
\end{thm}

More precisely, if the core of $\Ga$ has $k$ parts then there are $k+1$ ways to ``read'' the multiplicative quiver variety as a wild character variety, typically with different ranks/pole configurations (see the dictionary in \S\ref{ssn: dictionary}).
In particular this gives an alternative ``graphical'' way of thinking about a large class of wild character varieties, in terms of supernova graphs 
(rather than the usual perspective involving irregular curves \cite{gbs}).
This new viewpoint is useful since the graph determines a Kac--Moody algebra in a standard way and in particular a Kac--Moody root system and Weyl group (as described in \S\ref{sn: km}).
Passing between the different possible readings of a multiplicative quiver variety leads
to the fact that the Kac--Moody Weyl group 
acts to give symplectic isomorphisms, in the following sense:

\begin{thm}\label{thm: refln isoms-intro}
Suppose $\Ga$ is a simply-laced supernova graph with nodes $I$, coloured so that its core is monochromatic, $q,d$ are arbitrary, and 
$s_i\in \Aut(\IZ^I)$, $r_i\in\Aut((\IC^*)^I)$ are the corresponding simple reflections, generating the Kac--Moody Weyl group (cf. \S\ref{sn: km}).
Then if $q_i\neq 1$ the multiplicative quiver varieties
$$\cM^{st}(\Ga,q,d)\qquad\text{and}\qquad \cM^{st}(\Ga,r_i(q),s_i(d))$$
are isomorphic smooth symplectic algebraic varieties, for any node $i\in I$.
\end{thm}

In turn via Theorem \ref{thm: isom to wcv} this yields many more symplectic isomorphisms between wild character varieties, the typical orbit being infinite.
In certain special cases (when the core is complete bipartite)
similar explicit isomorphisms are known to arise from the action of the Fourier--Laplace transform on meromorphic connections \cite{BJL81, malg-book}.
Even in such special cases the proof that they are symplectic is new---the approach here is the first to establish 
precise algebraic symplectic isomorphisms between the full symplectic (Betti) moduli spaces.
Once lifted up to the quasi-Hamiltonian world a nice way to establish such isomorphisms becomes possible 
(\S\ref{sn: isos}), purely in the context of smooth affine varieties, although it will take some effort to set up the framework for this.

\subsection{Simple example}
The wild character varieties are a class of symplectic/Poisson algebraic varieties that generalise the character varieties of Riemann surfaces, i.e. the 
spaces of complex fundamental group representations of Riemann surfaces.
Deligne proved in 1970 \cite{Del70} that complex fundamental group representations of punctured smooth algebraic curves parameterise the connections on algebraic vector bundles which have {\em regular singularities} at each puncture. The wild character varieties parameterise more general, {\em irregular}, connections.
The main Poisson manifolds appearing in the theory of quantum groups are simple examples of wild character varieties (see \cite{hdr} \S4).

A basic fact, known to $\cD$-module experts, 
that we wish to emphasize
is that it is not so easy to delineate between tame and wild character varieties: some of the simplest wild character varieties are {\em isomorphic} to tame character varieties, i.e. to spaces of fundamental group representations. 
A fundamental example of this basic fact is the following statement:

\begin{thm}\label{thm: intro isom}
Suppose $\cM_{\text{tame}}^{st}$ is a (symplectic) tame character variety parameterising irreducible representations of the fundamental group of a punctured Riemann sphere (with arbitrary fixed local conjugacy classes of monodromy on arbitrary, fixed rank vector bundles). 
Then there is a wild character variety $\cM_{2+1}^{st}$
parameterising meromorphic connections on vector bundles on the Riemann sphere with just two poles of order $2$ and $1$ respectively and an isomorphism
$$\cM_\text{tame}^{st}\ \cong\ \cM^{st}_{2+1}$$
of algebraic symplectic manifolds.
Moreover all such $2+1$ (untwisted) wild character varieties are isomorphic to such tame cases. 
\end{thm}

On the other hand there are many wild cases which are not isomorphic to tame cases. 
(Due to \cite{smid,Sab99, wnabh, gbs} the wild character varieties enjoy the same key properties as the tame case, such as having complete \hk metrics with an underlying algebraic symplectic structure of a topological nature.) 
This theorem, and several generalisations, will be proved
by showing both sides are isomorphic to the same multiplicative quiver variety,  in \S\ref{sn: mqvs and wcvs}.
The (more general) dictionary relating the ranks of the bundles, number of poles, and local conjugacy classes will also be explained there.
Such isomorphisms may be viewed as ``multiplicative analogues'' of the isomorphisms of \cite{rsode, slims}.
As further concrete motivation,
the map appearing in Theorem \ref{thm: intro isom} will be described as directly as possible in the first section below.

The rough layout of this article is as follows.
The core of the article is \S\ref{sn: isos} which establishes some quasi-Hamiltonian isomorphisms. 
The three isomorphism theorems there relate incredibly complicated explicit expressions for
the multiplicative symplectic forms.
(Perhaps the main discovery of the article is that one can establish such isomorphisms inductively 
using the fission and fusion operations to reduce to the simplest case, which is still highly nontrivial.)
The earlier sections prepare the way, describing a language in 
which to make such inductive arguments, and sets 
up the theory of multiplicative quiver varieties.
Later sections give some applications of these isomorphisms, and discuss some related topics.

\ 

{\small
\noindent
{\bf Acknowledgments.}
\noindent
The author's line of thought here can be traced back to 
a question of N. Hitchin (from as long ago as 1995 \cite{Hit95long}) which was really about how to remove the Stokes data from the theory of Frobenius manifolds, using the Fourier--Laplace transform (cf. \cite{thesis} \S8). 
This tempted the author to go in the opposite direction and extend known technology to the wild world.
In any case, many thanks are due to him for this thought-provoking  question, and  equally to B. Malgrange for gently pointing out \cite{malg-book}, after the author had been enthusing about %
\cite{BJL81}\footnote{\cite{BJL81} was used in \cite{pecr},\cite{k2p} to find new finite braid orbits and to understand Okamoto's affine $D_4$ action in terms of Fourier--Laplace on the level of Betti data, giving an ``intrinsic'' counterpoint to the earlier ``extrinsic'' De\,Rham approach of Arinkin--Lysenko \cite{AL-p6ims}, (cf. \cite{quad} Remark 1). This example could be seen as a simple prototype of what we do here. 
}.
Thanks are also due to J\'er\'emy Blanc (for explaining Proposition 
\ref{prop: ncps} back in 2008 \cite{blanc-email}).
Many of these results were announced in \cite{rsode}, and lectured on
during the June 2013 
Moduli Spaces meeting at the CRM, Montreal.
The idea of constructing new noncommutative algebras using Stokes data, generalising the multiplicative 
preprojective algebras and the generalised DAHAs, was discussed in a talk in Nice in May 2009 \cite{nicetalk09}.
}

\section{Example explicit isomorphism}
We will sketch here the direct construction of the map between the spaces
in the example of Theorem \ref{thm: intro isom}, as an illustration of the more general results to be proved. 
The proof that such maps, and generalisations,  are symplectic is one of the main motivations, and will appear in \S\ref{sn: mqvs and wcvs}.

In this example the underlying map
arises by computing the action of Fourier--Laplace on Betti data 
\cite{BJL81, malg-book}. 
The space $\cM_\text{tame}$ may be described as follows.
There is an integer $m$, a complex vector space $V$ and conjugacy classes
$\cC_1,\ldots,\cC_m,\cC_\infty\subset \GL(V)$ so that
$$\cM_\text{tame} \cong 
\{ (T_1,\ldots,T_m)\in \GL(V)^m\st T_i\in \cC_i,\ 
 T_m\cdots T_2T_1\in\cC_\infty \}/\GL(V)$$
and a point $(T_1,\ldots,T_m)$ is stable if there is no proper nontrivial subspace $U\subset V$ with $T_i(U)\subset U$ for all $i$.
The corresponding Stokes data may be constructed as follows.
Given $T_i\in \cC_i$ let $d_i = \rank(T_i-1)$ and set 
$W_i = \IC^{d_i}$ for $i=1,\ldots,m$.
Then $T_i$ may be written as 
$T_i = 1+b_ia_i$
for linear maps $a_i:V\to W_i$ and $b_i:W_i\to V$ such that 
$a_i$ is surjective and $b_i$ is injective.
Moreover the injectivity/surjectivity conditions imply (cf. Appendix \ref{apx: relating orbits})
that the 
conjugacy class of $T_i$ is uniquely determined by (and uniquely determines) the conjugacy class  of 
$$h_i := 1+a_i b_i\in \GL(W_i).$$
Let $\breve\cC_i\subset \GL(W_i)$ denote the inverse conjugacy class, of 
$h_i^{-1}$, and define $W=\bigoplus W_i$ to be the (external) 
direct sum of the spaces $W_i$. 
Let $U_\pm\subset \GL(W)$ be the block triangular unipotent subgroups determined by the ordered grading of $W$.
Then define elements  $u_\pm \in U_\pm$ by the prescription
$$hu_+-u_- = [a_ib_j]\in \End(W)$$
where $a_ib_j\in \Hom(W_j,W_i)$, and 
$h\in H:=\Prod \GL(W_i)$ has 
components $h_i$.
The remarkable algebraic fact then (see  \S\ref{ssn: lin alg}) is that:
$$1+AB = u_-^{-1}hu_+\in \GL(W),\qquad
T_m\cdots T_2T_1 = 1+BA \in \GL(V)$$
where $A:V\to W$ and $B:W\to V$ are the  maps with components $(a_i T_{i-1}\cdots T_1)\in \Hom(V,W_i)$ and 
$b_i\in \Hom(W_i,V)$ respectively.
Moreover the fact that $(T_1,\ldots,T_m)$ is stable implies that
$B$ is surjective and $A$ is injective.
Thus fixing the conjugacy class $\cC_\infty$ of $T_m\cdots T_2T_1$ is equivalent to fixing the conjugacy class $\cC$ of $u_-^{-1}hu_+$
(cf. Appendix \ref{apx: relating orbits}).
Thus, by defining $S_1\in U_+,S_2\in U_-$ 
so that $hS_2S_1 = u_-^{-1}hu_+$, 
we obtain a point of the wild character variety 
$$\cM_{2+1}\cong \{ 
(S_1,S_2,h)\in U_+\times U_-\times H\st
h^{-1}_i\in \breve \cC_i,\  hS_2S_1\in \cC\}/H$$
where $H$ acts by diagonal conjugation,
and one may check the stability conditions match up, and so obtain a genuine isomorphism of moduli spaces of stable points.

Now, both sides have natural holomorphic symplectic structures. On the tame side it essentially goes back to Atiyah--Bott \cite{AB83}, 
and may be obtained as the multiplicative symplectic reduction
$$\cM_\text{tame} \cong 
\left(\cC_m\fus{}\cdots\fus{}\cC_2\fus{}\cC_1\right) 
\spqa{\cC_\infty} \GL(V)$$

\noindent
of the product of the conjugacy classes.
On the irregular side the symplectic structure goes back to the  extension of the Atiyah--Bott construction in \cite{smid}, and may be obtained as the multiplicative symplectic reduction
$$ \cM_{2+1}\  \cong\  G \sqpa{\cC} \gah\spqa{\breve \cC} H$$

\noindent
of the fission space $\gah\cong G\times U_+\times U_-\times H$ of \cite{saqh02, saqh, fission}, where $G=\GL(W)$.
The theorem asserts that these two algebraic symplectic structures match up under the given isomorphism (of stable points).
Our proof of this (and its generalisations) involves showing that both spaces arise as the {\em same}
multiplicative symplectic quotient, 
and are just two ways to ``read'' a certain multiplicative quiver variety.
\section{Nakajima quiver varieties}

We will recall the complex symplectic approach to Nakajima quiver varieties
\cite{nakaj-duke94, nakaj-duke98}, in a way that is convenient for the multiplicative analogue we have in mind.

\subsection{Representations of graphs}\label{ssn: reps of graphs}

Suppose $\Ga$  is a graph with nodes $I$ (and edges $\Ga$).
(We will always suppose both $I$ and $\Ga$ are finite sets
and, unless otherwise stated, that every edge connects two distinct nodes.)
Let $\bar\Ga$ be the set of oriented edges of $\Ga$, i.e. the set of pairs $(e,o)$ such that  $e\in \Ga$ is an edge of $\Ga$ and $o$ is a choice of one of the two possible orientations of $e$.
Thus if $a\in \bar \Ga$ is an oriented edge, the head $h(a)\in I$ and tail $t(a)\in I$ nodes of $a$  are well defined. 
For our purposes it is convenient to define a {\em representation} $\rho$ of the graph $\Ga$ to be the following data:

1) an $I$-graded vector space $V= \bigoplus_{i\in I} V_i$, and

2) for each oriented edge $a\in \bar\Ga$,
a linear map $\rho(a)=v_a : V_{t(a)}\to V_{h(a)}$ between the vector spaces at the head and the tail of $a$.

Thus the data in 2) amounts to choosing a linear map in both directions along each edge of $\Ga$.
(Equivalently $\bar \Ga$ may be viewed as a quiver, the double of $\Ga$ with any orientation, and a representation of the graph $\Ga$ is the same thing as a quiver representation of $\bar \Ga$.) 
A {\em subrepresentation} of a representation $V$ of $\Ga$ consists of an $I$-graded subspace $V'\subset V$ which is preserved by the linear maps, i.e. such that 
$v_a(V'_{t(a)})\subset V'_{h(a)}$ for each oriented edge 
$a\in \bar \Ga$.
A representation $V$ is {\em irreducible} if it has no proper nontrivial subrepresentations.
Given an $I$-graded vector space $V$ we may consider the set 
$\Rep(\Ga,V)$ of all
representations of $\Ga$ on $V$. 
This is just the vector space
$$\Rep(\Ga,V)=\bigoplus_{a\in \bar\Ga} \Hom(V_{t(a)},V_{h(a)})$$  
of all possible maps, in each direction along each edge of $\Ga$. 
(In the present article the notation $\Rep(\Ga,V)$ always denotes the representation of the graph $\Ga$, even if $\Ga$ has the additional structure of a quiver.)

Given a graph $\Ga$ and an $I$-graded vector space $V$ the group $H=\Prod\GL(V_i)$ acts on $\Rep(\Ga,V)$ via its natural action on $V$ preserving the grading.
Further a choice of orientation of the graph $\Ga$ determines a holomorphic symplectic structure on $\Rep(\Ga,V)$, and then the action of $H$ is Hamiltonian with a moment map
$$\mu : \Rep(\Ga,V) \ \to\ \lh^*=\Lie(H)^*\cong \Prod_{i\in I}\End(V_i).$$

\subsection{Additive/Nakajima quiver varieties.}\label{sn: nqv}
The Nakajima quiver varieties are defined by choosing a central value $\la\in \IC^I$ of the moment map and taking the symplectic quotient:
$$\cN(\Ga,\la,d) = \Rep(\Ga,V) \spqa{\la} H 
= \{ \rho\in \Rep(\Ga,V) \st \mu(\rho) = \la\}/H$$
where  $\la$ is identified with the central element
$\sum\la_i\Id_{V_i}$ of $\Lie(H)^*$.
Here the symbol $d\in \IZ^I$ denotes the dimension vector, the vector of dimensions of the components of $V$: $d_i=\dim(V_i)$, and the quotient is the affine quotient, taking the variety associated to the ring of $H$ invariant functions.
Further one can consider the open subset of stable points 
$\Rep(\Ga,V)^{st}\subset\Rep(\Ga,V)$ 
for the action of $H$ (defined as the points whose $H$-orbit is closed of dimension $\dim(H)-1$).
By results of King \cite{king-quivers} a graph representation $\rho$ is stable if{f} it is {\em irreducible}.
Considering stable orbits defines 
$\cN^{st}(\Ga,\la,d)\subset \cN(\Ga,\la,d)$.

\begin{rmk}\label{rmk: theta-nqv}
More generally one can consider other geometric invariant theory (GIT) quotients: given a weight $\th\in \IQ^I$ with $\sum \th_i d_i=0$, consider the nontrivial linearization $\cL=\Rep(\Ga,V)\times\IC\to \Rep(\Ga,V)$ 
by lifting the $H$ action 
(to the total space of the trivial line bundle $\cL$)
via the character $h\mapsto \prod \det(h_i)^{-m\th_i}$, for a positive integer $m$ so that $m\th_i\in \IZ$  for all $i\in I$.
Taking the proj of the graded ring of $H$-invariant sections of powers of $\cL$ restricted to $\mu^{-1}(\la)$ yields a more general (quasi-projective)
quiver variety  
$\cN(\Ga,\th,\la,d)$, and we can also consider its open subset of 
 $\th$-stable points $\cN^{st}(\Ga,\th,\la,d)$.
More explicitly 
$\th$ determines the
$\th$-stable and $\th$-semistable points
$\Rep(\Ga,V)^{\th\text{-st/}\th\text{-ss}}\subset \Rep(\Ga,V)$, and 
via \cite{king-quivers}
a representation $(\rho,V)$ of $\Ga$ 
is $\th$-stable if{f} any nontrivial proper subrepresentation $(\rho',V')$ satisfies $\dim(V')\cdot \th < 0$ (and $\rho$ is $\th$-semistable if instead $\dim(V')\cdot \th \le 0$).
Set-theoretically $\cN^{st}(\Ga,\th,\la,d)$ is the set of 
$H$-orbits in $\Rep(\Ga,V)^{\th\text{-st}}\cap \mu^{-1}(\la)$ 
and $\cN(\Ga,\th,\la,d)$ is the set of $H$-orbits in 
$\Rep(\Ga,V)^{\th\text{-ss}}\cap \mu^{-1}(\la)$ which are closed (in this subset).
We will mainly focus on the affine case ($\th=0$)  here.

\end{rmk}

\subsection{Fission and supernova graphs}\label{sn: snova}
In 2008 a dictionary was discovered \cite{rsode} relating 
the Nakajima quiver varieties for a special class of graphs, the {\em supernova graphs}, to moduli spaces of connections on curves. 
These graphs may be defined as follows. 
Let $V$ be a complex vector space and let $\lt\subset\End(V)$ be a Cartan subalgebra of the Lie algebra of $\GL(V)$, such as the diagonal matrices.
Choose a $\lt$-valued polynomial
$$Q(w) = \sum_1^r A_i w^i\in \lt[w]$$
with zero constant term.
(We will often call $Q$ an `irregular type' cf. \cite{gbs}.)
Given $Q$ we obtain a grading $V=\bigoplus_I V_i$ of $V$ into the 
eigenspaces of $Q$, so that 
$$Q = \sum q_i(w) \Id_i$$
where $\Id_i$ is the idempotent for $V_i$ and the $q_i$ are distinct elements of $\IC[w]$.
Then we obtain a graph $\Ga(Q)$ with nodes $I$ by joining exactly 
\beq\label{eq: fission mults}
\deg(q_i-q_j)-1
\eeq
edges between the nodes $i$ and $j$ for any $i\neq j\in I$.
The graph $\Ga(Q)$ is the ``fission graph'' attached to the irregular type $Q$, and they are basic examples of supernova graphs.\footnote{
Fission graphs were first introduced in an equivalent fashion in \cite{rsode} p.29---the irregular type appearing there is $dQ$ written in the coordinate $z=1/w$. The procedure described there to construct them by sequentially breaking the vector space $V$ into smaller pieces explains the name ``fission''.} 
\footnote{
The Stokes data of a connection with irregular part $dQ$ may be parameterised by 
representations of the graph with $\deg(q_i-q_j)$ edges 
between nodes $i,j$: this is the ``Stokes graph'' of $Q$, but our focus here is on the fission graphs.}

For example if $Q = A_2w^2+A_1w$ has degree at most two, then the fission graph of $Q$ will be simply-laced (have no multiple edges).
This is the main case of interest to us here.
In fact the simply-laced fission graphs are exactly the complete $k$-partite graphs:

Recall that a graph $\Ga$ with nodes $I$ is a complete $k$-partite graph 
if there is a partition $I=\bigsqcup_{j\in J} I_j$ of its nodes into $k$ nonempty parts $I_j$ labelled by a set $J$ with $\#J=k$, 
such that two nodes are connected by a single edge if and only if they are not in the same part.
If $Q = A_2w^2+A_1w$ then define $J$ to be the eigenspaces of $A_2$ and let $I_j$ be the eigenspaces of $A_1$ inside the $j$ eigenspace of $A_2$ (so that $I=\bigsqcup_{j\in J} I_j$ is the set of simultaneous 
eigenspaces of $A_1,A_2$).
This identifies the fission graph $\Ga(Q)$ with the complete $k$-partite graph determined by this partition of $I$, where $k=\#J$ is the number of eigenspaces of $A_2$.
Many  complete $k$-partite graphs are drawn in 
the figures of \cite{rsode, slims}: they are determined by a partition $P$ 
of an integer into $k$ parts (corresponding to the partition of the simultaneous eigenspaces of $A_1,A_2$ into the eigenspaces of $A_2$).
For example the graph $\Ga(1,1)$ corresponding to the partition $1+1$ is just a single edge connecting two nodes (the interval), and similarly 
$\Ga(1,1,1)$ is the triangle and 
$\Ga(2,2)$ is the square (a complete bipartite graph).
The star-shaped graph with $n$ legs of length one is the bipartite graph 
$\Ga(1,n)$. 
The graphs $\Ga(n)$ have $n$ nodes and no edges, and
the graphs $\Ga(1,1,\ldots,1)$ are the complete graphs (with every pair of nodes connected by a single edge).

\begin{defn}\label{dfn: sn graph}
A simply-laced  {\em supernova graph}
 is a graph obtained by gluing a leg (of length $\ge 0$) on to each node of a simply-laced fission graph\footnote{
In general the supernova graphs are defined as follows (this is equivalent to \cite{rsode} Appendix C): 
take any fission graph $\Ga(Q)$ with nodes $I_0$ and introduce some new nodes $I_1$. Define a new graph $\Ga(Q)'$ with nodes $I_0\sqcup I_1$ by adding a single edge from each node of $I_0$ to each node of $I_1$. 
Finally the supernova graph is obtained by  gluing a leg (of length $\ge 0$) on to each node of the core $\Ga(Q)'$.}.
\end{defn}

Here a ``leg'' of length $l$ is just a  Dynkin graph of type 
$A_{l+1}$, with $l$ edges. 
The initial fission graph will be referred to as the {\em core} of the 
supernova graph.
Note that by definition a graph is ``star-shaped'' if and only if it is a  supernova graph with core  $\Ga(1,n)$ for some $n$.

\begin{figure}[ht]
	\centering
	\input{quivers2.pstex_t}
	\caption{Example supernova graph, with core $\Ga(4,2,2)$}
\end{figure}
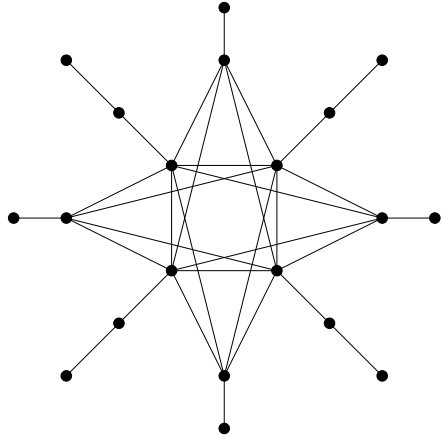

The basic result relating supernova graphs to moduli spaces of connections is:
\begin{thm}[\cite{rsode, slims}]
Suppose $\Ga$ is a simply-laced supernova graph with nodes $I$, and we choose a dimension vector $d\in \IZ^I$ and parameters $\la\in \IC^I$.  
Then the Nakajima quiver variety $\cN^{st}(\Ga,\la,d)$ is isomorphic to a moduli space of stable meromorphic connections $\cM_{st}^*$ on the holomorphically trivial bundle on the Riemann sphere.
\end{thm}

In the star-shaped case this result follows by combining
the known relation between $\GL_n(\IC)$ coadjoint orbits and quivers (\cite{Kraft-Procesi-InvMath79}\cite{nakaj-duke94} \S7), with the standard symplectic quotient construction of moduli spaces of Fuchsian systems (see e.g. \cite{Hit95long, smid}), and it was first used by Crawley-Boevey
\cite{CB-additiveDS}. 
The general result rests on the symplectic approach of \cite{smid} to more general  connections.
See \cite{rsode}, \cite{slims} \S9 for full details, and for discussion of the reflection isomorphisms on such  moduli spaces of connections, corresponding to the quiver variety reflection functors.
Taking nontrivial weights for the quiver varieties corresponds to taking nontrivial parabolic weights for the connections (cf. \cite{ihptalk}).
The dictionary relating general (non-simply-laced) supernova graphs to connections was given in \cite{rsode} Appendix C. (A proof of a result announced  there appeared recently in \cite{hi-ya-nslcase}).

This dictionary gives a different way to organise moduli spaces of connections, and suggests where to look for interesting examples (by taking interesting graphs):
For example if the graph is a triangle (with dimension vector $d=(1,1,1)$) the corresponding quiver variety is a
Gibbons-Hawking gravitational multi-instanton of type $A_2$ (which is a deformation of the minimal resolution of the Kleinian singularity $\IC^2/\IZ_3$): the above dictionary says, for example, that this very natural space is isomorphic to a moduli space of connections on rank two vector bundles on the Riemann sphere with one order three pole and a first order pole.

Our basic aim is to establish a ``multiplicative analogue'' of this result, thereby giving a graphical approach that captures all the points of some of the full moduli spaces of connections $\cM$ considered in \cite{smid, wnabh, ihptalk} (which allow connections on nontrivial holomorphic bundles), rather than just the open pieces $\cM^*\subset \cM$.
First we will discuss  the theory of multiplicative symplectic geometry that will be required.

\begin{rmk}
Note that a standard way (cf. e.g. \cite{nakaj-duke98} p.520)  to think about representations of graphs
is in terms of block matrices, i.e. via  
the inclusion $\Rep(\Ga,V) \subset \End(V)$ (if there are no multiple edges), whence the symplectic structure pairs the corresponding entries in the upper triangular part and the lower triangular part (\cite{nakaj-duke98} (3.3)).
As will be explained below (in Prop. \ref{prop: open subset of reps}), 
Stokes data give a natural 
multiplicative version of this picture, with the  upper/lower triangular block matrices above replaced by an open subset of the  upper/lower triangular unipotent block matrices.
\end{rmk}

\section{Quasi-Hamiltonian geometry}

Many of our results will be framed in the language of quasi-Hamiltonian geometry, so we will briefly review the main features of this theory.
In essence this theory is a multiplicative version of the usual Hamiltonian theory, with moment maps taking values in Lie groups rather than the dual of the Lie algebra.
The axioms for the analogue of the symplectic form and its interaction with the group action and the moment map are  more complicated.
The upshot is a direct and explicit algebraic approach to constructing certain quite exotic symplectic manifolds.
The original article \cite{AMM} worked with compact groups $K$ and then 
gave a new/alternative finite dimensional construction of the 
(real) symplectic structure on moduli spaces of flat $C^\infty$ connections on principal $K$-bundles over Riemann surfaces with fixed local monodromy conjugacy classes
(such moduli spaces admit K\"ahler structures, via the Narasimhan--Seshadri correspondence or its parabolic analogues).
Here, as in \cite{saqh02, saqh, fission, gbs}, we will work with complex reductive groups, such as  $G=\GL_n(\IC)$, as a means to construct algebraic complex symplectic manifolds---in those articles it was shown that the quasi-Hamiltonian approach is actually quite flexible and can be extended to construct new moduli spaces, not just moduli spaces of flat $C^\infty$ connections.
Specifically the quasi-Hamiltonian toolbox was enlarged: besides the 
conjugacy classes and the doubles, some new basic pieces were introduced
and it was shown that these were sufficient to give a finite dimensional  construction of complex symplectic structures on moduli spaces of meromorphic connections appearing in the wild nonabelian Hodge correspondence of \cite{Sab99, wnabh}.
We will show these  new pieces yield a new relation to quivers.

\subsection{Quasi-Hamiltonian spaces}
Fix a connected complex reductive group $G$ (in this article $G$ will always be a product of general linear groups).
Our notation and conventions  are as in \cite{gbs}, so in particular 
$\theta,\bar\theta$ denote 
the left and right-invariant Maurer-Cartan forms on $G$, and we fix a nondegenerate invariant bilinear form $(\ ,\ )$ on $\g=\Lie(G)$.
A complex manifold $M$ is a 
complex quasi-Hamiltonian $G$-space
if there is an action of $G$ on $M$, 
a $G$-equivariant map $\mu:M\to G$ (where $G$ acts on itself by
conjugation) and a $G$-invariant holomorphic two-form
$\omega\in \Omega^2(M)$ such that:

\noindent(QH1). 
The exterior derivative of $\omega$ is the pullback along the moment map of the canonical three-form on $G$:
$d\omega = \mu^*(\th^3)/6,$

\noindent(QH2).
For all $X\in \g$,
$\omega(v_X,\cdot\,) = \frac{1}{2}\mu^*(\theta+\overline\theta, X)
\in \Omega^1(M)$ where $v_X$ is the fundamental vector field of $X$,

\noindent(QH3).
At each point $m\in M$: 
$\ker \omega_m \cap \ker  d\mu = \{0\} \subset T_mM$.

As a basic example \cite{AMM}, any conjugacy class $\cC\subset G$ is a quasi-Hamiltonian $G$-space with the moment map given by the inclusion and action given by conjugation, analogously to coadjoint orbits in the usual Hamiltonian story.
There is also a multiplicative version of the symplectic quotient:

\begin{thm}[\cite{AMM}]
Let $M$ be a quasi-Hamiltonian $G\times H$-space with moment map
$(\mu,\mu_H):M\to G\times H$ and suppose
that the quotient by $G$ of the inverse image $\mu^{-1}(1)$  
of the identity under the first moment map is a manifold.
Then the restriction of the two-form $\omega$ to $\mu^{-1}(1)$ 
descends to the {\em reduced space}
\begin{equation}\label{eq: qh quot}
M\spq G := \mu^{-1}(1)/G 
\end{equation}
and makes it into a quasi-Hamiltonian $H$-space. 
\end{thm}

In particular, if $H$ is abelian (or in particular trivial) 
then the quasi-Hamiltonian axioms imply $M\spq G$ is a complex symplectic manifold; this will be our main method of constructing symplectic manifolds. 
We will mainly work in the complex algebraic category, with $M$ a smooth affine variety (and $\mu, \omega$ algebraic), so that $\mu^{-1}(1)$ is an affine subvariety. 
The quotient $\mu^{-1}(1)/G$ is then taken to be  
the  geometric invariant theory quotient, namely the affine variety associated to the ring of $G$-invariant functions on $\mu^{-1}(1)$.
The points of this quotient  
correspond bijectively to the {\em closed} $G$-orbits in $\mu^{-1}(1)$, and so in general it is different to the set-theoretic quotient. Alternatively one may view the 
points of the  geometric invariant theory quotient as parameterising the quotient of $\mu^{-1}(1)$ by a stronger equivalence relation than orbit equivalence (S-equivalence): two points are S-equivalent if their orbit closures intersect.

The fusion product, which puts a ring structure on the  
category of quasi-Hamiltonian $G$-spaces, is defined as follows. 

\begin{thm}[\cite{AMM}]\label{thm: fusion}
Let $M$ be a quasi-Hamiltonian $G\times G\times H$-space, 
with moment map $\mu=(\mu_1,\mu_2,\mu_3)$. 
Let $G\times H$ act by 
the diagonal embedding $(g,h)\to (g,g,h)$. 
Then $M$ with two-form 
\begin{equation} \label{eqn: fusion 2form}
\wt{\omega}= \omega - \frac{1}{2}(\mu_1^* \theta, \mu_2^* \overline \theta)
\end{equation}
and moment map
$$\wt{\mu} = (\mu_1\cdot \mu_2,\mu_3):M\to G\times H$$
is a quasi-Hamiltonian $G\times H$-space. 
\end{thm}

If $M_i$ is a quasi-Hamiltonian $G\times H_i$ space for $i=1,2$ their 
fusion product $$M_1\fus M_2$$ 
is defined to be the
quasi-Hamiltonian $G\times H_1\times H_2$-space 
obtained from the
quasi-Hamiltonian $G\times G \times H_1 \times H_2$-space 
$M_1\times M_2$ by
fusing the two factors of $G$.
If $\cC\subset G$ is a conjugacy class then the quasi-Hamiltonian reduction at $\cC$ may be defined
$$M\spqa{\cC} G = \mu^{-1}(\cC)/G \cong (M\fus \cC') \spq G$$
by fusing $M$ with the inverse conjugacy class $\cC'$ and then reducing at the value $1$ of the resulting moment map. (Note in the algebraic context that care is needed, and will be taken, if $\cC$ is not semisimple and thus not affine.)

It is convenient to formalise the notion of {\em gluing} 
quasi-Hamiltonian spaces, as follows. 
Given a quasi-Hamiltonian $G\times G\times H$-space $M$, 
we may fuse the two $G$ factors to obtain a quasi-Hamiltonian $G\times H$ space. Then, if the quotient is well-defined, we may
reduce by the $G$ factor (at the identity of $G$) to obtain a quasi-Hamiltonian $H$-space, the {\em gluing} of the two $G$-factors.
Thus for example if $M_i$ is a  quasi-Hamiltonian $G\times H_i$-space for $i=1,2$ then $M_1$ and $M_2$ may be glued to obtain a quasi-Hamiltonian $H_1\times H_2$ space (if it is a manifold) by gluing %
 their product:

$$M_1\glue{G} M_2 := (M_1 \fusion{G} M_2)\spq G.$$

If $G=\GL(V)$ for some vector space $V$, the subscript $G$ will often be replaced by $V$, and if the factors to be glued are clear from the context the subscript will be omitted and we will write $M_1\glu M_2$.
In most of the cases considered here the $G$-action will be free with a global slice so there is no problem performing the gluing.
Note that whereas fusion is only commutative up to isomorphism, the gluing operation is actually commutative (when it is defined).

\subsection{General linear fission spaces}\label{ssn: gl fission spaces}

Some new examples of complex quasi-Hamiltonian spaces were defined in \cite{saqh02, saqh, fission, gbs}, the {\em (higher) fission spaces}.
Here we will use some of these spaces in the case of general linear groups.
Gluing on these spaces defines a new operation, {\em fission},
enabling the structure group to be changed, and leads to many new examples of symplectic manifolds.

Suppose $V$ is a finite dimensional complex vector space with an ordered grading, i.e. with a decomposition 
$$V = V_1\oplus V_2\oplus \cdots \oplus V_k$$
for some $k$.
Then we may consider the groups $G=\GL(V)$ and $H=\Prod \GL(V_i)\subset G$
together with the parabolic subgroup $P_+ \subset G$  stabilising the flag
$$F_1\subset F_2\subset \cdots \subset F_k=V$$
where  $F_i= V_1\oplus\cdots \oplus V_i$,
and the opposite parabolic $P_- \subset G$ 
 stabilising the flag
$$F'_k\subset \cdots \subset F'_2\subset F'_1=V$$
where  $F'_i= V_i\oplus\cdots \oplus V_k$,
 and
their unipotent radicals $\U_\pm\subset P_\pm$.
Thus in an adapted basis $U_{+/-}$ is the subgroup of block upper/lower triangular matrices with $1$'s on the diagonal (respectively) and $H$ is the block diagonal subgroup of $G$.
Then given an integer $r\ge 1$, define the 
{\em higher fission space}
$$\cA^r(V)=\cA^r(V_1,\ldots,V_k) = \gahr := G \times (U_+\times U_-)^{r}\times H.$$
If $r=1$ the superscript will be omitted so that $\cA(V) = \cA^1(V)$.
A point of $\cA^r(V)$ is given by specifying $C\in G, h\in H$ and 
$\bS\in (U_+\times U_-)^{r}$ with $\bS = (S_1,\ldots,S_{2r})$ where
$S_{\text{even}}\in U_-$ and $S_{\text{odd}}\in U_+.$
The unipotent elements $S_i$ will usually be referred to as ``Stokes multipliers''.
The group  $G\times H$ acts on $\cA^r(V)$ as follows:
$$(g,\eta)(C,\bS,h)  = (\eta Cg^{-1}, \eta\bS \eta^{-1}, \eta h\eta^{-1})$$
where $(g,\eta)\in G\times H$ and $\eta\bS \eta^{-1} = (\eta S_1\eta^{-1},\ldots,\eta S_{2r}\eta^{-1})$.

\begin{thm}(\cite{gbs})\label{thm: gl fission spaces}
Suppose $V$ is an ordered graded vector space and $r \ge 1$ is an integer.
Then $\cA^r(V)$ is a quasi-Hamiltonian 
$G\times H$-space, with moment map
$$\mu(C,\bs,h) = (C^{-1} h S_{2r}\cdots S_2 S_1 C,\,\, h^{-1}) \in G\times H.$$
\end{thm}

For example if $k=1$ so $V$ only has one graded piece (and $r$ is arbitrary) then $H=G$ and both groups $U_\pm$ are trivial, so 
$\cA^r(V)=G\times G=D(G)$ is the ``double'' of \cite{AMM}.
In the case $H$ is a torus, $\cA^r(V)$ is the space $\wt\cC/L$ of 
\cite{saqh02, saqh} Remark 4.
If $r=1$ the spaces $\cA^r(V)$ appear in \cite{fission}. The general formula for the quasi-Hamiltonian two-form is given in \cite{gbs} equation (9).

If $r\ge 2$ we may reduce by $G$ at the identity to define a
quasi-Hamiltonian $H$-space
\beq
\cB^r(V) := \cA^r(V)\spq G = \mu_G^{-1}(1)/G
\eeq
where $\mu_G:\cA^r(V) \to G$ is the $G$-component of $\mu$.
If $r=2$ this space will be denoted $\cB(V):=\cB^2(V)$.
If $H$ is a torus each space $\cB^r(V)$ is a symplectic manifold.

After the first fission spaces were constructed  in \cite{saqh02} 
(when $H$ is a maximal torus),
M. Van den Bergh (see \cite{vdb-doublepoisson, vdb-ncqh, yamakawa-mpa}) found some further new examples of complex quasi-Hamiltonian spaces. 
He showed that, for any finite dimensional complex vector spaces $V_1,V_2$, the space
\beq\label{eq: vdb space}
\{(a,b)\in \Hom(V_2,V_1)\oplus\Hom(V_1,V_2)\st \det(1+ab)\ne 0\}.
\eeq
is a quasi-Hamiltonian $\GL(V_1)\times \GL(V_2)$ space 
with moment map $((1+ab)^{-1}, 1+ba)$.
(The underlying space, without quasi-Hamiltonian structure, is well-known in the local classification of $\cD$-modules, e.g. \cite{malg-book} p.31.)
Subsequently it was shown in \cite{gbs} Theorem 4.2 that if $V$
has exactly two graded pieces ($V=V_1\oplus V_2$, i.e. $k=2$) then 
the reduced fission space
$$\cB(V)=\cB(V_1,V_2)$$ 
is isomorphic (as a quasi-Hamiltonian space) to Van den Bergh's space \eqref{eq: vdb space}.
Whereas the spaces \eqref{eq: vdb space} facilitate the construction of the classical multiplicative quiver varieties, 
the more general spaces $\cB(V)$
will be key to construct generalisations here.

Many properties of the fission spaces were established in \cite{gbs}, and it is possible to show that up to isomorphism $\cA^r(V)$
only depends on the graded vector space $V$ and not on the choice of ordering of the graded pieces.
This follows from \cite{gbs} Theorem 10.4, but since it is important here we will explain how it follows in more detail.

\begin{prop}\label{prop: imd to change parab gln} %
If
$$V=V_1\oplus \cdots \oplus V_n$$ 
is an ordered graded vector space
and 
$V'$ is the same graded vector space but with a different ordering,
then $\cA^r(V)\cong \cA^r(V')$.
\end{prop}
\pf
The idea is the same for any $r$ so we will assume $r=1$ (and omit to write it).
Further we may suppose the orderings of $V$ and $V'$ differ only by swapping two consecutive subspaces, say $V_k,V_{k+1}$.
Recall that a sequence of unipotent groups $\cU=(U_1,\ldots U_r)$
of $G$ ``directly spans'' a unipotent group $U\subset G$
if the product map $\cU\to U; (u_1,\ldots,u_r)\mapsto 
u_1\cdots u_r$ is an isomorphism of spaces.
In our situation $\U_\pm$
have direct spanning decompositions (in any order)
$$U_+  = \prod_{i<j} U_{ij},\qquad U_-  = \prod_{i>j} U_{ij}$$
where $U_{ij} = 1+\Hom(V_j,V_i)\subset G$.
Further one may check they also have direct spanning decompositions
\beq \label{eq: DiSpe1}
U_+  = U_{k, k+1} U_+^\circ, \qquad 
U_-  =   U_{k+1,k} U^\circ _-
\eeq
where $$U_+^\circ = \prod_{\overset{i<j}{(ij)\neq (k,k+1)}} U_{ij},\qquad
U_-^\circ = \prod_{\overset{i>j}{(ij)\neq (k+1,k)}} U_{ij}.$$
On the other hand $\cA(V') = G\times H\times U'_+\times U'_-$
where $\U'_\pm$ have direct spanning decompositions
\beq \label{eq: DiSpe2}
U'_+  =  U_+^\circ U_{k+1, k}, \qquad 
U'_-  =   U^\circ _-  U_{k,k+1}.
\eeq
Thus an isomorphism between $\cA(V)$ and $\cA(V')$ is obtained by first using the direct spanning equivalences \eqref{eq: DiSpe1}, \eqref{eq: DiSpe2}
(cf. \cite{gbs} \S6.1)
and then the isomonodromy isomorphism (of \cite{gbs} \S6.2) corresponding to 
$$
(U_{k, k+1},  U_+^\circ, U_{k+1,k}, U^\circ _-) 
\ \mapsto\ 
( U_+^\circ, U_{k+1,k} , U^\circ _-, U_{k, k+1}). $$
This gives an isomorphism of spaces, and the quasi-Hamiltonian structures match up due to 
\cite{gbs} Lemma 6.1 and Proposition 6.3.
\epf

Another operation studied in \cite{gbs} was the possibility to ``nest'' the fission spaces. 
In the case of general linear groups
the nesting results imply that the fission spaces attached to ordered graded vector spaces behave well under refining the gradings, as follows.
Suppose $V=\bigoplus_1^k V_i$ is an ordered graded vector space and for some fixed index $i=f$ that we have an ordered grading 
$V_f = \bigoplus_1^l W_j$
of the vector space $V_f$. 
Then we can consider the refined ordered grading $\bU$ of $V$ 
given by
$$\bU=V_1\oplus V_2\oplus\cdots\oplus V_{f-1}\oplus
W_1\oplus W_2\oplus\cdots\oplus W_l\oplus V_{f+1}\oplus\cdots\oplus V_k.$$ 
\begin{prop} \label{prop: gl nesting} %
The gluing of $\cA^r(V_1,\ldots,V_k)$ and $\cA^r(W_1,\ldots,W_l)$ via the action of $\GL(V_f)$ is isomorphic 
to $\cA^r(\bU)$:
\beq\label{eq: splaying}
\cA^r(V_1,\ldots,V_k)\glue{V_f}\cA^r(W_1,\ldots,W_l) \ \cong\ \cA^r(\bU).
\eeq
\end{prop}
\pf
If we write 
$K=\Prod \GL(V_i)$ and $H=\Prod\GL(U_i)$, where $U_i$ is the $i$th vector space in $\bU$, then the left-hand side of \eqref{eq: splaying}
is isomorphic $\cA^r(V)\glue{K}\papk{K}{H}{r}$ (because $\papk{K}{H}{r}$ is just the product of $\cA^r(W)$ and lots of doubles, and the doubles disappear in the gluing). Then the result follows from the nesting result \cite{gbs} Corollary 6.5.
\epf

\section{Quasi-Hamiltonian spaces of graph representations}

This section sets up some further definitions 
related to graphs and 
then establishes a relation between some fission spaces and some representations of graphs---this is the key step to define the multiplicative quiver varieties in the next section.

\subsection{Orderings of graphs}
Recall (from \S\ref{ssn: reps of graphs}) that we have defined the notion of representation of a graph $\Ga$ with nodes $I$
and have recalled that  a graph is complete $k$-partite if there is a partition $I=\bigsqcup_J I_j$ of its nodes into $k=\#J$ parts, and two nodes are connected by an edge if{f} they are in different parts.

\begin{defn} \label{defn: ordering of ckpg}
An {\em ordering} of a complete $k$-partite graph $\Ga$ with nodes $I$ is the choice of :

1) a total ordering of its parts (so we can label the parts $I_1,\ldots ,I_k$), and 

2) a total ordering of the nodes in each part $I_j$.
\end{defn}

By declaring $I_i<I_j$ if $i<j$ (and using the given ordering in each part) such an ordering determines a total order of $I$.
Note that given such an ordering one then gets an orientation of $\Ga$ (so that edges go $i\to j\in I$ if $i>j$) but in general the ordering contains more information than the associated orientation. 

Recall that a ``quiver'' is an oriented graph, i.e. a graph together with a choice of orientation for each edge. 
In other words a quiver is obtained by gluing together some oriented edges.
In turn an edge is a complete bipartite graph $\Ga(1,1)$, and an orientation of $\Ga(1,1)$ is equivalent to an ordering of it (in the sense of Definition \ref{defn: ordering of ckpg}).
Thus we may generalise the notion of  `oriented edge' as follows.

\begin{defn} \label{}
A {\em simple coloured quiver} $\Ga$ is an ordered complete $k$-partite graph.
\end{defn}

Below (in Definition \ref{defn: coloured quiver}) we will generalise the notion of a quiver by gluing together lots of simple coloured quivers---the case when each simple coloured quiver is an edge (a copy of $\Ga(1,1)$) will yield the usual notion of quiver.

\subsection{Quasi-Hamiltonian spaces of graph representations}

Suppose $\Ga$ is a simple coloured quiver with nodes $I$, and $V$ is a fixed $I$-graded vector space, so we may consider the space $\Rep(\Ga,V)$ of representations of the graph $\Ga$ on $V$.
(Recall $\Rep(\Ga,V)$ denotes the space of representations of the graph underlying $\Ga$.)
Then, considering the parts $I_j\subset I$, 
we can define ordered graded vector spaces
$$W_j = \bigoplus_{i\in I_j} V_i$$
and $W=\bigoplus W_j$ (which is $V$ with a different grading).
Thus there is a nested sequence of groups:
$$H:=\Prod \GL(V_i) \ \subset\ 
 K:=\Prod \GL(W_j) \ \subset\ 
 G:=\GL(V).$$
This enables us to consider the fission spaces $\cA(W_j)$ for each $j$, and the quasi-Hamiltonian $K$-space 
\beq
\cB(W)= \cA^2(W)\spq G.
\eeq
The basic result linking graphs to the fission spaces is the following.

\begin{prop}\label{prop: open subset of reps}
Suppose $\Ga$ is a simple coloured quiver with 
nodes $I$ and $V$ is an $I$-graded vector space.
Then there is a canonical $H$-invariant nonempty open subset
$$\Rep^*(\Ga,V) \subset \Rep(\Ga,V)$$ 
 of the space of representations of the graph $\Ga$ 
on $V$ which is a smooth affine variety and a 
quasi-Hamiltonian $H$-space, canonically isomorphic to
$$ \cB(W)\glue{K} \Prod_j\cA(W_j).$$
\end{prop}

\pf
Given a representation $(v_{ij})$ of $\Ga$ on $V$, use the ordering of 
$I$ to define the following  unipotent elements %
\beq\label{eq: vpm from rho} 
v_+ = 1 + \sum_{i<j} v_{ij},\qquad v_-=1 + \sum_{i>j} v_{ij}
\eeq
of $\GL(V)$, where we set $v_{ij}\in \Hom(V_j,V_i)$ to be zero if $i,j$ are in the same part of $I$.
Then consider the subset $R^*=\Rep^*(\Ga,V)$ of $\Rep(\Ga,V)$ such that $v_-v_+$ is in the opposite big cell of $\GL(V)$ determined by the ordered grading of $V$, i.e. so that we may write
\beq\label{eq: big cell factn}
v_-v_+ = w_+gw_-
\eeq
for some $g\in H$ and unipotent elements 
$w_+ = 1 + \sum_{i<j} w_{ij},
w_-=1 + \sum_{i>j} w_{ij}$ 
with $w_{ij}\in \Hom(V_j,V_i)$ (which are allowed to be nonzero for $i,j$ in the same part).
Note that if such a factorisation \eqref{eq: big cell factn} exists, then it is unique.
This subset $R^*$ is defined by the nonvanishing of the polynomial function
\beq\label{eq: minor prod}
f:= \Prod_{i\in I} \Delta_i : \Rep(\Ga,V)\to \IC
\eeq
where $\Delta_i$ is the minor of $v_-v_+$ corresponding to the summand
$\bigoplus_{j\ge i} V_j$ of $V$.
$R^*$ is nonempty as it contains the representation with each $v_{ij}=0$.
Thus on one hand $R^*$ is isomorphic to the affine variety 
$\{z.f=1\}$ in $\IC\times \Rep(\Ga,V)$. On the other hand 
it is a straightforward unwinding of the definitions to identify $R^*$ with $\cB(W)\glue{K} \Prod_j\cA(W_j)$.
Indeed specifying a point of this quasi-Hamiltonian space amounts to solving
$$\kappa S_4S_3S_2S_1 = 1  \in \GL(W), \qquad
\kappa_j = h_j s_{2j}s_{1j} \in \GL(W_j)$$
where the elements $S_\bullet, s_{\bullet j}$ are the Stokes multipliers in $\cB(W)$ and $\cA(W_j)$ respectively,
where   $\kappa\in K$ has components $\kappa_j$ and where 
$h_j\in \Prod_{i\in I_j}\GL(V_i)\subset \GL(W_j)$.
Thus, setting $v_+= S_1, v_- = S_2$  it is easy to translate these equations into the form \eqref{eq: big cell factn}.
Observe in particular that $g\in H$ 
(from \eqref{eq: big cell factn}) then has components $h^{-1}_j$
and so (from Theorem \ref{thm: gl fission spaces}) 
$g:\Rep^*(\Ga,V)\to H$ is the moment map for the $H$ action.
\epf

Thus for example if $W=\bigoplus_1^k W_j$ 
the space $\cB(W)$ itself is an open subset of the space of representations on $W$ of the complete graph with $k$ nodes.
The $k=2$ case of this example is
the Van den Bergh space \eqref{eq: vdb space}, 
corresponding to the interval, and then the nontrivial matrix entries of $v_+,v_-$ are the elements $(a,b)$ in \eqref{eq: vdb space}, as in \cite{gbs} \S4.

\begin{rmk}
Given an irregular type $Q=A_2w^2+A_1w\in \lt[w]$ as in \S\ref{sn: nqv}   
(with $\lt\subset \End(V)$) 
Proposition \ref{prop: open subset of reps}  says, 
in effect, that the space $\cB(Q)=\cA(Q)\spq\GL(V)$ 
is isomorphic to an open subset of 
$\Rep(\Ga(Q),V)$, where $\Ga(Q)$ is the fission graph of  $Q$ (from \cite{rsode} Appendix C and \S\ref{sn: snova} above)  and
 $\cA(Q)$ is from \cite{gbs} Theorem 7.6 (with $z=1/w$; 
this will be fleshed out in 
the proof of Proposition \ref{prop: wcv and omqv} below).
The same result holds with similar proof for any irregular type, not just the simply-laced case.
This yields even more general multiplicative quiver varieties, to be studied elsewhere.
\end{rmk}

\subsection{Coloured graphs and quivers}
Now we wish to consider graphs built out of coloured pieces, each of which is a complete $k$-partite graph.
Let $C$ be a nonempty finite set (which we will call the set of {\em colours}). 

\begin{defn}
A ``coloured graph'' is a graph $\Ga$  (with nodes $I$ and edges $\Ga$), together with
a colouring of each edge (i.e.  a map $\ga:\Ga\to C$ to the set $C$ of colours), such that each connected component of each monochromatic subgraph 
$$\Ga_c=\ga^{-1}(c)\subset \Ga$$ 
is a complete $k$-partite graph for some $k$ (dependent on $c\in C$).
\end{defn}

For simplicity here we assume $\Ga_c\subset \Ga$ is an embedded copy of a complete $k$-partite graph (so it is connected and no nodes get identified), but it is convenient when drawing pictures to allow $\Ga_c$ to be disconnected. 
We will say a coloured graph is {\em  classical} if each 
connected component of each monochromatic subgraph $\Ga_c$ consists of just one edge (for example if the map $\ga$ is injective). 
For example any graph has a classical ``tautological colouring'' by taking $C=\Ga$ and $\ga$ to be the identity map.

\begin{defn}\label{defn: coloured quiver}
A ``coloured quiver'' is a coloured graph $\Ga$  (with nodes $I$, edges $\Ga$ and colour map $\ga$), together with
the choice of ordering of each subgraph $\Ga_c\subset \Ga $ (as in Definition \ref{defn: ordering of ckpg}).
\end{defn}

\noindent
(In general the orderings chosen for the nodes of different intersecting monochromatic subgraphs may be completely independent.)
Note that a classical/tautological coloured quiver is just a quiver:
if $\Ga$ is classical then
choosing an ordering is the same as choosing an orientation of each edge.

\begin{cor}\label{cor: open subset of reps 2}
Suppose $\Ga$ is a coloured quiver with 
nodes $I$, $V=\bigoplus V_i$ is an $I$-graded vector space and
$H=\Prod\GL(V_i)$.
Then there is a canonical nonempty  $H$-invariant open subset
$$\Rep^*(\Ga,V) \subset \Rep(\Ga,V)$$ 
 of the space of representations of the graph $\Ga$ 
on $V$ which is a smooth affine variety. 
Further,
given a choice of ordering of the colours $C$ for each node $i\in I$,
then $\Rep^*(\Ga,V)$
is a
quasi-Hamiltonian $H$-space.
\end{cor}

Henceforth the open subset $\Rep^*(\Ga,V)$ will be referred to as the set of ``invertible'' representations of $\Ga$ on $V$.
The term ``open multiplicative quiver variety'' will also be used, alluding to the fact that we have not yet quotiented by $H$ (cf. \cite{rsode}). 

\pf
Given a colour $c\in C$, let $I_c\subset I$ be 
the nodes of $\Ga_c\subset \Ga$, and let 
$V_c=\bigoplus_{i\in I_c} V_i$, $H_c=\Prod_{i\in I_c}\GL(V_i)$.
Thus
$\Rep^*(\Ga_c,V_c)$ is a quasi-Hamiltonian $H_c$-space, by Proposition \ref{prop: open subset of reps}. 
Denote by $g_{ci}$ 
 the $\GL(V_i)$ component of the moment map 
$g_c:\Rep^*(\Ga_c,V_c)\to H_c$.
Then consider the space 
$$\Rep^*(\Ga,V) := \Prod_{c\in C} \Rep^*(\Ga_c,V_c).$$
By Proposition \ref{prop: open subset of reps}
this is
an affine variety, a nonempty  open subset of $\Rep(\Ga,V)$, and
a quasi-Hamiltonian $\Prod{H_c}$-space. 
If we fix a node $i\in I$ and use an ordering of the colours $C$ (at this node $i$) %
then we may (internally) 
fuse together all the copies of the group $\GL(V_i)\subset H_c$ 
for the colours $c$ meeting the node $i$.
Repeating for each node this gives $\Rep^*(\Ga,V)$
the structure of quasi-Hamiltonian $H$-space. 
The moment map $\mu:\Rep^*(\Ga,V)\to H$ has components 
\beq\label{eq: mmap cmpts on IM}
\mu_i = \prod_{\{c\,|\,i\in I_c\}} g_{ci}\in\GL(V_i)
\eeq
where the factors are ordered according to the chosen ordering of 
the colours.
\epf

As for any affine variety with an action of a reductive group, 
the stable points of $\Rep^*(\Ga,V)$ are defined as the 
points whose $H$-orbit is closed and of dimension $\dim(H)-\dim(\Ker)$,
where $\Ker$ is the kernel of the action, i.e. the kernel  of the map
$H\to \Aut(\Rep^*(\Ga,V))$, which has dimension $1$ in the present set-up.
The following simple remark will be very useful:

\begin{lem}\label{lem: stablem}
An invertible graph representation $\rho\in\Rep^*(\Ga,V)$
is stable if and only if $\rho$ is a stable point of $\Rep(\Ga,V)$, i.e. (via \cite{king-quivers}) if and only if $\rho$ is irreducible.
\end{lem}
\pf
$\Rep^*(\Ga,V)$ is defined as the nonvanishing locus of a function $F:\Rep(\Ga,V)\to \IC$ (obtained by taking the product of the functions \eqref{eq: minor prod} over each monochromatic component).
This function is constant on the $H$-orbit $\cO$ of $\rho$, and so $\cO$ is closed in $\Rep^*(\Ga,V)$ if and only if it is closed in $\Rep(\Ga,V)$.  
\epf

\begin{rmk}\label{rmk: subreps}
In fact if $\rho\in\Rep^*(\Ga,V)$ and $\rho'\in \Rep(\Ga,U)$
is a graph subrepresentation of $\rho$ for some graded subspace $U\subset V$
then $\rho'\in \Rep^*(\Ga,U)$.
One way to see this is to let $\rho''=V/U$ and note that 
$\rho'\oplus \rho''$ represents a point in the closure of the orbit of 
$\rho$ in $\Rep(\Ga,V)$ (via the relation between $1$-parameter subgroups and filtrations \cite{king-quivers} p.521), and so is in $\Rep^*(\Ga,V)$ as argued in the proof of the lemma.
Consequently $\rho'\in \Rep^*(\Ga,U)$ as all the minors 
\eqref{eq: minor prod} factor, corresponding to the $I$-graded vector space decomposition $V\cong U\oplus V/U$.
\end{rmk}

To avoid further discussion of the choice of ordering of the colours $C$ for each node we could fix once and for all 
an ordering of the colours $C$ and always use this at each node.
In fact,
up to isomorphism, the space of invertible representations does not depend on any such choices:

\begin{thm}\label{thm: reordering}
Fix a coloured graph $\Ga$ with nodes $I$ and an $I$-graded vector space $V$.
The space $\Rep^*(\Ga,V)$ 
is independent of all the further choices of ordering needed 
to define it, up to isomorphism of quasi-Hamiltonian $H$-spaces.
\end{thm}
\pf
That it is independent of the orderings of the nodes (and parts) of the monochromatic subgraphs $\Ga_c$ follows from the isomonodromy isomorphisms as in Proposition
\ref{prop: imd to change parab gln}.
That it is independent of the orderings 
of the colours at each node, follows from the braid isomorphisms of \cite{AMM} Theorem 6.2, i.e. that fusion is commutative up to isomorphism.
\epf

\subsection{How to colour a supernova graph.}
To end this section recall that a (simply-laced) supernova graph was defined (p.\pageref{dfn: sn graph}) by gluing some legs on to a complete $k$-partite graph, and that the central complete $k$-partite graph was called the core.
In this article, unless otherwise stated, a supernova graph will always be coloured so that its core is monochromatic i.e.:

{\em Two edges should have the same colour if and only if they are in 
the core.}

If a supernova graph is star-shaped (i.e. if its  core has the form $\Ga(1,n)$ for some $n$)---then 
we will see that in fact this colouring is {\em equivalent} to the tautological/classical colouring (in the sense that the corresponding spaces of invertible representations are isomorphic), 
but that this is not the case in general.

By definition a {\em supernova quiver} is then
a supernova graph together with a choice of ordering of each monochromatic subgraph.

\section{Multiplicative quiver varieties.}\label{sn: col quiver vars}

In this section we will define a class of varieties attached to coloured graphs.
These varieties generalise the classical multiplicative quiver varieties suggested in \cite{CB-Shaw} (see also \cite{vdb-doublepoisson, yamakawa-mpa}).
Looking at the case of a triangle shows that in general 
new spaces are obtained this way. 
Subsequently we will show that for any supernova quiver the corresponding multiplicative quiver varieties arise as moduli spaces of Betti data for meromorphic connections on curves, and so have hyperk\"ahler metrics.

Let $\Ga$ be a coloured quiver with nodes $I$ and colours $C$.
Choose a dimension vector
$$d\in \IZ^I$$
where $d=(d_i)_{i\in I}$ has $d_i\ge 0$,
and some parameters
$$ q \in (\IC^*)^I$$
such that $$ q^d := \Prod_{i\in I}q_i^{d_i} = 1.$$
Further suppose we have made 
an additional choice of an ordering of the set $C$ of colours for each node $i\in I$ (for example we could choose and fix an ordering of $C$ before hand).
Then define $V_i = \IC^{d_i}$ and let $V=\bigoplus V_i$ be the 
corresponding $I$-graded vector space. 
By Corollary \ref{cor: open subset of reps 2}
the space $\Rep^*(\Ga,V)$ is a quasi-Hamiltonian $H$-space where 
$H=\Prod\GL(V_i)$, and we
now identify $q$ with the point $(q_i \Id_{V_i})$ of $H$.

\begin{defn}
The ``\,{\em multiplicative quiver variety}'' of $\Ga,q,d$ is the quasi-Hamiltonian reduction of $\Rep^*(\Ga,V)$
at the value $q$ of the moment map:
$$\cM=\cM(\Ga,q,d)  =  \Rep^*(\Ga,V)\spqa{q}H = \mu^{-1}(q)/H.$$ 
\end{defn} 

\noindent
Here the quotient by $H$ on the right is the affine geometric invariant theory quotient, taking the affine variety associated to the ring of $H$-invariant functions on the affine variety $\mu^{-1}(q)\subset \Rep^*(\Ga,V)$.
(Sometimes these varieties will be referred to as ``coloured multiplicative quiver varieties'' although this could be confusing since the classical case is the most colourful.)
Note that since $q$ is in the centre of $H$, only a cyclic ordering of the colours at each node is needed (rather than a total ordering).

\begin{rmk}\label{rmk: qtothed}
Observe that the relation \eqref{eq: big cell factn} implies $\det(g_c)=1$ for any colour $c\in C$.
Together with the determinants of the moment map conditions
$\prod_{\{c\,|\,i\in I_c\}} g_{ci} = q_i\Id_{V_i}$
this implies $q^d=1$. Thus if we chose $q,d$ so that $q^d\neq 1$ then $\cM$ would be empty.
\end{rmk}

One can further consider the open subset $\cM^{st}\subset \cM$ of {\em stable} points; by definition the stable points correspond to the closed orbits in $\mu^{-1}(q)$ of dimension $\dim(H)-1$ (noting that the diagonal copy of $\IC^*$ embedded in $H$ 
always acts trivially).

\begin{thm} \label{thm: stable cmqv}
$\cM^{st}(\Ga,q,d)$ is a smooth algebraic symplectic manifold
which is either empty or of dimension $2-(d,d)$, where $(\ , \ )$
is the bilinear form \eqref{eq: KM bil form} on the root lattice of 
$\Ga$. 
In terms of representations of graphs, the points of $\cM^{st}$
correspond to the $H$-orbits in $\mu^{-1}(q)$ of irreducible representations of the graph $\Ga$.
\end{thm}
\pf
All the work has been done to set this up as a quasi-Hamiltonian quotient and the result now follows as for 
classical multiplicative quiver varieties, which is explained nicely in 
\cite{yamakawa-mpa}. 
The final statement follows from Lemma
\ref{lem: stablem} and \cite{king-quivers}.
\epf

\begin{rmk}\label{rmk: phi-mqv}
Similarly to the additive case (see Remark \ref{rmk: theta-nqv}), 
and to the classical multiplicative case (see \cite{yamakawa-mpa}), 
given a weight 
$\phi\in\IQ^I$ with $\sum \phi_id_i=0$ 
we can define more general multiplicative quiver varieties
$\cM^{st}(\Ga,\phi,q,d)\subset \cM(\Ga,\phi,q,d)$
by considering the linearization $\cL\to\Rep^*(\Ga,V)$
defined by the character $h\mapsto \Prod\det(h_i)^{-m\phi_i}$ of $H$.
In other words $\cM(\Ga,\phi,q,d)$ is the proj of the graded ring of $H$-invariant sections of powers of $\cL$ restricted to $\mu^{-1}(q)$.
As in Lemma \ref{lem: stablem}
$$\Rep^*(\Ga,V)^{\phi\text{-st}/\phi\text{-ss}} =  
\Rep(\Ga,V)^{\phi\text{-st}/\phi\text{-ss}} \cap \Rep^*(\Ga,V)$$
and so set theoretically 
$\cM(\Ga,\phi,q,d)$ is the set of $H$-orbits in 
$\mu^{-1}(q)\cap \Rep(\Ga,V)^{\phi\text{-ss}}$
which are closed (in this subset) and 
$\cM^{st}(\Ga,\phi,q,d)$ is just the set of closed $H$-orbits in 
$\mu^{-1}(q)^{\phi\text{-st}}=\mu^{-1}(q)\cap\Rep(\Ga,V)^{\phi\text{-st}}$.
To keep the paper within a reasonable length the affine case $\phi=0$ will be our focus here. Note that 
$\phi$ will match up with the weights for Betti data of connections (i.e. the $\ga_i$ of \cite{wnabh} Remark 8.2, or the $\phi$ in \cite{logahoric}),
and such varieties will parameterise filtered Stokes $G$-local systems (\cite{gbs} Remark A.5), as shown in \cite{yamakawa-mpa} in the tame/star-shaped case.
Note that the isomorphisms in \S\S\ref{sn: mqvs and wcvs},\ref{sn: reflns} generalise immediately to the case with $\phi\neq 0$, by taking nontrivial linearisations on both sides, since in essence we are saying they are the same multiplicative symplectic quotient of the same smooth affine variety.   
\end{rmk}

Given a fixed dimension vector $d$, we will say that the parameters $q$ are {\em generic} if they 
obey the condition $$q^\al\ne 1$$
for any $\al$ in the finite set 
$$R_\oplus(d):= 
\{\al \in \IZ^I\st (\al,\al)\le 2 \text{ and } 0\le \al_i\le d_i 
\text{ for all $i$}\}\setminus \{0,d\}.$$
(Beware that the set of generic $q$ is not always dense in $\{q\st q^d=1\}$, for example if $d$ is an integer multiple of another dimension vector.) 

\begin{prop}
If the parameters are generic then all points of the 
multiplicative quiver variety are stable, 
$\cM(\Ga,q,d) = \cM^{st}(\Ga,q,d)$,
and so it is smooth.
\end{prop}
\pf 
If a representation $\rho$ is not stable then it has a stable subrepresentation, of dimension vector $\al$ say.
One may check this subrepresentation again obeys the moment map conditions (for the same $q$), cf. Remark \ref{rmk: subreps}.
Thus $q^\al=1$
as in Remark \ref{rmk: qtothed}.
But it is stable so $\cM^{st}(\Ga,q,\al)$ is nonempty and so Theorem \ref{thm: stable cmqv} implies $(\al,\al)\le 2$, which is not possible if $q$ is generic \epf

Let us record here the main reordering and recolouring results:

\begin{cor}{\em(Reordering.)} \label{cor: reordering}
Up to isomorphism %
$\cM(\Ga,q,d)$ 
only depends on $q,d$ and the underlying coloured graph $\Ga$, and not on any of the chosen orderings.
\end{cor}
\pf
This follows from Theorem \ref{thm: reordering}.
\epf

\begin{thm} {\em(Recolouring.)} \label{thm: recolouring}
Suppose $\Ga$ is a coloured graph such that for some colour $c\in C$ the subgraph $\Ga_c\subset \Ga$ is star-shaped.
Then  let $\Ga'$ be the new coloured graph obtained from $\Ga$ by choosing distinct new colours for each edge of $\Ga_c$ and adding these new colours to $C$.
Then $\cM(\Ga,q,d)\cong \cM(\Ga',q,d)$ for any $q,d$.
\end{thm}
\pf
This will follow immediately from  Theorem \ref{thm: basic isom} below.
\epf

\subsection{Example multiplicative quiver varieties}
\label{ssn: example mqv}

Up to isomorphism we are thus able to associate a variety to a coloured graph $\Ga$ and data  $d\in \IZ^I, q\in (\IC^*)^I$, where $I$ is the set of nodes of $\Ga$.

If the underlying graph is star-shaped 
then
all the different colourings of the graph
give the same varieties (via Theorem \ref{thm: recolouring}); they are all classical multiplicative quiver varieties.

Thus the simplest case where new varieties may appear is the case when 
$\Ga$ is a triangle. 
By recolouring  (as in Theorem \ref{thm: recolouring}) 
we see there are at most two inequivalent ways to colour the triangle. Namely either 1) all the edges are different colours, or 2) all the edges are the same colour. 
(If just two edges were the same colour, they make up a star-shaped subgraph and we could recolour them.)
Thus in case 1) we are considering the classical multiplicative quiver variety of the triangle (in the sense of Crawley-Boevey--Shaw \cite{CB-Shaw}), and in 2) we are considering a non-classical (supernova)  quiver variety.

\begin{figure}[ht]
	\centering
	\input{trianglesbw2.pstex_t}
	\caption{Left: the classical colouring, and right: the monochromatic/supernova colouring of a triangular graph $\Ga$.}\label{fig: two triangles}
\end{figure}
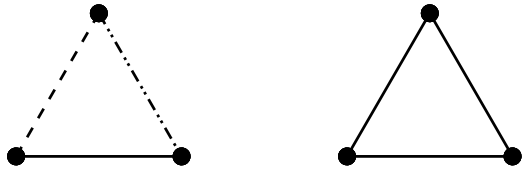

Our aim here is to show that these two varieties are not isomorphic
in general when the dimension vector is $d=(1,1,1)$, in which case the varieties
are complex surfaces (complex dimension two).
Later we will show that the supernova quiver varieties are isomorphic to wild character varieties and thus carry complete \hk metrics; on the other hand it is not clear if the spaces 1) arise in relation to nonabelian Hodge theory.

\begin{prop}\label{prop: triangle variety descriptions}
Let $\Ga$ be the triangle and fix dimension vector $d=(1,1,1)$. 
Suppose $q$ is generic. Then the corresponding varieties are smooth complex algebraic surfaces and 

1) the classical multiplicative quiver variety of $\Ga,q,d$ is obtained by blowing up three points on a triangle of lines in $\IP^2$ and then removing the strict transform of the triangle,

2) the supernova variety of $\Ga,q,d$ is obtained by blowing up three points on a conic in $\IP^2$ and then removing the strict transform of the union of the conic and a line.
\end{prop}

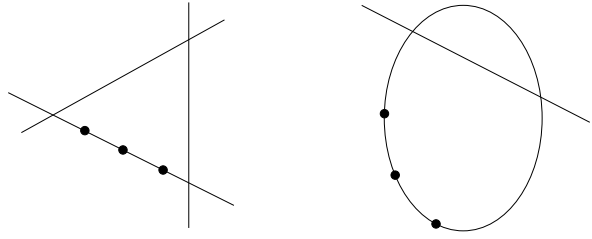
\begin{figure}[ht]
	\centering
	\input{3pointblowups.pstex_t}
	\caption{The two types of complex surfaces if $\Ga$ is a triangle.}\label{fig: blow up for triangle}
\end{figure}

Before proving this we deduce that it implies they are not algebraically isomorphic. This follows from the following fact (the author is grateful to J\'er\'emy Blanc  \cite{blanc-email} for explaining this):

\begin{prop}\label{prop: ncps}
For $i=1,2$ let $S_i$ be a smooth complex projective surface and let $D_i\subset S_i$ be a  curve consisting of exactly $a_i$ irreducible components. If $S_1\setminus D_1$ is isomorphic to $S_2\setminus D_2$ then $(K_{S_1})^2 + a_1 = (K_{S_2})^2 + a_2$.
\end{prop}

Thus taking $S_1$ to be $\IP^2$ blown up in three points on a line, and $S_2$ to be $\IP^2$ blown up in three points on a conic,
$(K_{S_1})^2 = (K_{S_2})^2 = 6$, but, taking $D_1$ to be the strict transform of the triangle of lines, and $D_2$ to be the strict transform of the union of a conic and a line, $a_1 = 3\ne a_2 = 2$
so they are not isomorphic.

\pfms (of Proposition \ref{prop: triangle variety descriptions}).
Since $d=(1,1,1)$ we have $V_i=\IC$ for each $i$ so we may identify each of the six linear map $v_{ij}$ with a complex number.
Thus in both cases we are considering subvarieties of the
quotient $\IC^6/H$ where $k\in H\cong (\IC^*)^3$ acts diagonally as 
$k(v_{ij}) = k_iv_{ij}/k_j$ (and so there is a $\IC^*$ subgroup of $H$ acting trivially).
The invariants for this action are generated by 
$$a= v_{12}v_{21},\quad b= v_{13}v_{31},\quad c= v_{23}v_{32},$$
$$ p = v_{12}v_{23}v_{31},\qquad r =  v_{21}v_{13}v_{32}$$
which satisfy the relation $abc=pr$. Thus $\IC^6/H$ is identified with the hypersurface in $\IC^5$ cut out by $abc=pr$.

Now, using an appropriate ordering, 
the relations defining the supernova variety fix the diagonal matrix $h$ in the factorisation $v_+v_-=w_-hw_+$. Thus $h_1$ is the top-left matrix entry of $v_+v_-$ and $h_1h_2$ is the determinant of the top-left $2\times 2$ submatrix of $v_+v_-$.
Expanding and rewriting in terms of the invariants  gives the relations
$$1+a+b = h_1,\quad 1+b+c+ac-p-r=h_1h_2$$
where $h_1,h_2\in \IC^*$ are constants (directly related to $q$).
Using the first equation to eliminate $b$, and the second to eliminate $r$, then setting $w=p+a$ (and translating the remaining variables 
$w,a,c$ by suitable constants, then relabelling them $x,y,z$) enables us to convert the remaining relation $abc=pr$ into a relation of the form $$xyz - x^2 + c_1 x + c_2 y+ c_3 z + c_4$$
defining a hypersurface in $\IC^3$, for constants $c_i$.
If we now project this to the $x$-$y$ plane, there is a unique $z$ over each point $x,y$ not on the conic $P(x,y)=xy+c_3$.
(In fact $c_3= h_1\neq 0$.) Indeed we may rewrite the relation as $P(x,y)z  + Q(x,y)$ with 
$Q =  c_1 x + c_2 y + c_4 -x^2$, 
and the hypersurface is then expressed 
as the blow up of $\IC^2$ at the three finite points of intersection of the conics $P$ and $Q$ (minus the strict transform of $P$).
(In particular $z$ can take any value over these three points, and the complement of the image of the projection is the conic $P$ minus these three points.)
Rephrasing this in terms of the projective plane, compactifying $\IC^2$ by adding the line at infinity,  yields 2).

On the other hand for 1) the relations  may be taken to be
$$1+a+b+ab=q_1,\quad 1+c= q_2(1+a),\quad abc=pr.$$
Using the second relation to eliminate $c$, and then combinations of the first and third to determine $b$ we get an equation of the form
$$W(a)+ (1+a)pr=0$$
for a cubic polynomial $W$.
Projecting this to the $a$-$p$ plane, the value of $r$ is uniquely determined by the point $(a,p)$ if it is off of the two lines $\{a=-1\} \union \{p=0\}$. 
Over these two lines lie just three affine lines, each above one of the three points  of the line $p=0$ where $a$ is a root of $W$ (and none of these roots is ever at $a=-1$).
Rephrasing this in terms of the triangle consisting of these two lines together with  the line at infinity in $\IP^2$ gives the desired description.
\epfms

\begin{rmk}\label{rmk: square example etc}
Similarly one can consider the case when $\Ga$ is a square with dimension vector $d=(1,1,1,1)$ and again the spaces corresponding to the two inequivalent colourings are not isomorphic (the explicit descriptions are as above but with a fourth blow up point added, to the same component, in each diagram in  Fig. \ref{fig: blow up for triangle}).
Note that the next type $A$ affine Dynkin diagram $\wh A_4$, the pentagon, is not a supernova graph, so we cannot further extend this discussion in this direction
(this is related to the fact that there is no second order Painlev\'e equation with this affine Weyl symmetry group, nor is there a two dimensional Hitchin system for $\wh A_4$).  
However, although beyond the scope of the present discussion,
one can consider the affine $A_1$ Dynkin graph (two nodes connected by two edges). 
The corresponding (non-simply-laced) supernova varieties are of the form 
$\cA^3(V_1,V_2) \spq_{\!q} \GL(V)\times\GL(V_1)\times\GL(V_2)$ 
(cf. \cite{rsode} appendix C)
and for $\dim(V_1)=\dim(V_2)=1$ this is a surface, and again not isomorphic to the corresponding classical multiplicative quiver variety
(the explicit descriptions are as above, but with only two blow up points, on the same component, in each diagram in  Fig. \ref{fig: blow up for triangle}).
\end{rmk}

Note that there is a conjectural list of all the
wild Hitchin spaces which are complex surfaces (\hk manifolds of real dimension four), in \cite{ihptalk} \S3.2, and the above examples are the Betti descriptions of some of them.
The next simplest class of examples seem to be the 
{\em higher Painlev\'e spaces} (\cite{rsode}, \cite{slims}\S11.4): 
for each of the complex 
surfaces $\cM$ listed in 
\cite{ihptalk} \S3.2 there is a sequence of spaces $\cM_{n}$ of complex dimension $2n$ for each $n\in\IZ_{\ge 1}$.
It is conjectured in \cite{slims} that $\cM_{n}$  is diffeomorphic to the Hilbert scheme of $n$-points on $\cM$ (this conjecture has been proved recently in the tame/star-shaped case in 
\cite{groechenig.hilbert.schemes}).
If we consider the cases where $\cM$ is of type $D^{(1)}_4,A^{(1)}_3,A^{(1)}_2$ respectively then the Betti description of $\cM_n$ is the multiplicative quiver variety associated to the coloured graphs in Figure \ref{fig: hpcg 6-4}, and this seems to be the most direct way to construct the underlying manifold.

\begin{figure}[ht] 
	\centering
	\input{p654.v2.pstex_t}
	\caption{} \label{fig: hpcg 6-4}
\end{figure}

\section{Isomorphism theorems}\label{sn: isos}

In this section we will establish some quasi-Hamiltonian isomorphisms that will be useful in later sections.
These isomorphisms relate incredibly complicated explicit expressions for the quasi-Hamiltonian two-forms, and at first sight look to be completely out of reach---but it turns out a quite geometric, inductive, proof is possible using the fission idea, although the simplest case is still too complicated to do by hand.
(The reader is strongly encouraged to draw some diagrams in order to understand the proofs of these isomorphisms; this is how they were found.) 
In the complete bipartite case
these isomorphisms are closely related to the work of Malgrange \cite{malg-book} and Balser--Jurkat--Lutz \cite{BJL81} 
computing the action of Fourier--Laplace on certain spaces of Betti data, i.e. Stokes and monodromy data (see also \cite{k2p} (36) and Remark 20 for more details on \cite{BJL81}). 
Even in these cases we are ``upgrading'' their isomorphisms to the quasi-Hamiltonian level.\footnote{More precisely the isomorphisms here are slightly different to those of \cite{malg-book} XII, since for example we do not reverse the orientation.
However it is straightforward to adjust what we do here to show the exact maps in \cite{malg-book} are also quasi-Hamiltonian.}

\subsection{Linear algebra} \label{ssn: lin alg}
First we will prepare some basic linear algebra.
Let $V,W$ be two finite dimensional complex vector spaces.
Choose elements $x\in \Hom(W,V)$, $y\in\Hom(V,W)$.
Note the following elementary fact:

\begin{lem}
$1+xy\in \End(V)$ is invertible if and only if $1+yx\in \End(W)$ is invertible.
\end{lem}

Occasionally  we will 
refer to $1+yx$ as the `dual' of $1+xy$ (if the choice of the maps $x,y$ is clear).
Now suppose $W$ has an ordered grading $W=\bigoplus_1^s W_i$ 
(i.e. $W$ has a direct sum decomposition into subspaces $W_i$ and we have chosen an ordering of the subspaces).
Let $\iota_i,\pi_i$ be the corresponding inclusions and projections
 and write
$x_i=x\circ \iota_i\in \Hom(W_i,V)$, 
$y^i=\pi_i\circ y\in\Hom(V,W_i)$
for the corresponding components.
Write
$$\varphi_i := 1+ x_1y^1+\cdots + x_iy^i\in \End(V)$$
for $i=0,1,\ldots,s$. 
Suppose furthermore that $x,y$ are such that
$\varphi_i$ is invertible for each $i$.
Then define 
$$
\wh x_i = \varphi_{i-1}^{-1}\circ x_i,\qquad
\wh y^i = y^i \circ \varphi_{i-1}^{-1}
$$
and let $$h_i = 1+\wh y^i x_i = 1 + y^i\wh x_i\in \End(W_i),$$ 
$$T_i= 1 + x_i\wh y^i,\quad M_i= 1+ \wh x_iy^i\in \End(V).$$

\begin{lem}
For all $i$
$$\varphi_i = T_i\cdots T_1 = M_1\cdots M_i\in\Aut(V)$$
and so in particular each $T_i,M_i$ (and thus by the previous lemma $h_i$) is invertible.
\end{lem}
\pf
For $i=1$ this is clear. In general $\varphi_i = \varphi_{i-1} + x_iy^i$
and so by induction 
$\varphi_i=T_{i-1}\cdots T_1 + x_i\wh y^i (T_{i-1}\cdots T_1) = T_i\cdots T_1.$
Similarly
$\varphi_i= M_1\cdots M_{i-1} + (M_1\cdots M_{i-1})\wh x_i y^i = M_1\cdots M_i.$
\epf

\begin{lem} \label{lem: gen gram}
The element $1+yx\in \End(W)$ admits a block Gauss decomposition, with block diagonal entries $h_i$.
Explicitly $1+yx = u_-^{-1}hu_+ = v_-hv_+^{-1}$
where
$$u_- = 1-[\wh y^ix_j]_{i>j},\quad hu_+ = 1 + [ \wh y^i x_j]_{i\le j},$$
$$v_-h = 1+[y^i\wh x_j]_{i\ge j},\quad v_+ = 1 - [  y^i \wh x_j]_{i< j},$$
where, for $a_{ij}\in\Hom(W_j,W_i)$, the expression $[a_{ij}]_{i<j}$ denotes the element of $\End(W)$ with $a_{ij}$ in its $(ij)$ block for $i<j$ and zero elsewhere.
In particular $hu_+-u_-$ and $v_-h-v_+$ are the ``generalized block Gram matrices":
$$hu_+-u_-=[\wh y^ix_j],       \qquad v_-h-v_+=[y^i\wh x_j].$$
\end{lem}
\pf Compute the products $u_-(1+yx)$ and $(1+yx)v_+$.\epf

Note that $\varphi_s=1+xy$, which is `dual' to $1+yx$, and similarly $\varphi_i$ is `dual' to  the submatrix of $1+yx$ in $\End(\bigoplus_1^iW_j)$, explaining why such decomposition exists.

Thus the set of $(x,y)\in T^*\Hom(V,W)$ for which all the $\varphi_i$ are invertible, has different coordinates given by $(\wh x_i, y^i)$ or $(x_i,\wh y^i)$ restricted such that all the $h_i$ are invertible.

\begin{eg}
Suppose $V$ has a nondegenerate symmetric bilinear form 
$(\cdot,\cdot)$ and basis $x_1,\ldots,x_s\in V$ 
such that $(x_i,x_i)=2$ for each $i$.
Set $W_i=\IC$ for each $i$ and let
$\wh y^i = -(x_i,\cdot)\in V^*=\Hom(V,\IC)$.
Then $T_i=1-x_i(x_i,\cdot)$ 
is the orthogonal reflection associated to $x_i$.
If we identify $W=\bigoplus W_i$ with $V$ via the basis $x_i$ and work in this basis (so that $x$ is the identity matrix)
then the above results imply
$$T_s\cdots T_1 = 1+y = -u_-^{-1}u_+$$
where $u_+ + u_-$ is the (symmetric) 
Gram matrix with entries $(x_i,x_j)$ (and the element $h=-1$).
This triangular decomposition of the product of reflections goes back  at least to Killing and Coxeter (cf. \cite{k2p} Remark 18).
\end{eg}

Now suppose further that there is also a direct sum decomposition $V=\bigoplus_1^r V_i$. 
Then we may repeat the above discussion using this decomposition: 
Define 
$$\gamma_j := 1+ y_1x^1+\cdots y_jx^j\in \End(W)$$
for $j=0,1,\ldots,r$,
where $x^j:W\to V_j, y_j:V_j\to W$ are the components of $x$ and $y$. 
Then if $\gamma_j$ is invertible for all $j$ we may define 
$$
\wh y_j = \gamma_{j-1}^{-1}\circ y_j,\qquad
\wh x^j = x^j \circ \gamma_{j-1}^{-1}
$$
and 
$$g_j = 1+\wh x^j y_j = 1 + x^j\wh y_j\in \End(V_j),$$ 
$$R_j= 1 + y_j\wh x^j,\quad N_j= 1+ \wh y_jx^j\in \End(W)$$
so that 
$$\ga_j = R_j\cdots R_1 = N_1\cdots N_j\in\Aut(W)$$ 
is invertible for all $j$ (so $1+yx=\ga_r=N_1\cdots N_r$), 
and  $1+xy\in \Aut(V)$ admits a block Gauss decomposition with block diagonal entries $g_j$:
$$1+xy = w_-gw_+^{-1}\in \Aut(V),\qquad
w_-g - w_+ = [x^i\wh y_j]\in \End(V).$$

\subsection{Isomorphism theorems}

We will establish three isomorphism theorems of increasing complexity.
Recall from \S\ref{ssn: gl fission spaces} the quasi-Hamiltonian spaces
$\cA^r(V), \cB^r(V)$.
The first case  is the following.

\begin{thm} \label{thm: basic isom}
Suppose $V, W_1,\ldots, W_s$ are finite dimensional complex vector spaces and $W$ is the ordered graded vector space $\bigoplus_1^s W_i$. 
Then there is an (explicit) isomorphism
\beq\label{eq: star equality}
\cB(W_s,V) \fusion{V}\cdots \fusion{V}\cB(W_1,V) \quad \cong\quad  
\cA(W)\glue{W} \cB(W,V) 
\eeq

\noindent
of quasi-Hamiltonian $\GL(V)\times \Prod\GL(W_i)$ spaces.
\end{thm}

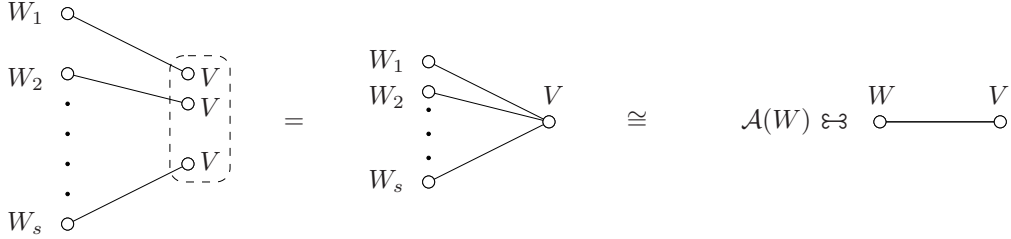
\begin{figure}[ht]
	\centering
	\input{splaying.isom3.pstex_t}
	\caption{First isomorphism.}\label{fig: splaying}
\end{figure}

In terms of graphs this isomorphism may be pictured as in Figure \ref{fig: splaying}. Namely on the left we are fusing $s$ edges together at one end, and on the right we are applying the fission operator $\cA(W)\glue{}(\,\cdot\,)$ to the edge representing $\cB(W,V)$, to break the vector space $W$ into pieces. This  operation is the multiplicative analogue of the splaying  operation of \cite{rsode} (also called $0$-fission in op. cit. Appendix C).

Note that this implies the recolouring result 
(Theorem \ref{thm: recolouring}), 
since  the spaces on the left would have each edge a distinct colour, whereas the spaces on the right correspond to a single monochromatic piece.

\pf
First, by a direct calculation postponed to \S\ref{ssn: s=2 case}, 
the $s=2$ case of this may be verified.
To deduce the general result define 
$W'=\bigoplus_1^{s-1} W_i$
and observe the following relations:
\begin{align*}
\cB(W_s,V) 
\fusion{V}
\Bigl(\cB(W_{s-1},&V)\fusion{V}\cdots\fusion{V}\cB(W_{1},V)
\Bigr)
\\
\ &\cong \
\cB(W_s,V) 
\fusion{V} 
\Bigl(\cA(W') 
\glue{W'}  
\cB(W',V)  \Bigr)
\qquad\text{(by induction)}\\
\ &= \ 
\cA(W')
\glue{W'}
\Bigl(
\cB(W_s,V)\fusion{V}\cB(W',V)
\Bigr)  
\qquad\text{}\\
\ &\cong \ 
\cA(W') 
\glue{W'}
\Bigl(\cA(W',W_s)\glue{W} \cB(W,V) \Bigr)\\  
\end{align*}
where the last isomorphism follows from the $s=2$ case.
Thus to establish the theorem it is  sufficient to find an isomorphism
$$\cA(W') 
\glue{W'}
\cA(W',W_s) \quad\cong\quad \cA(W). 
$$
But this is a special case of Proposition \ref{prop: gl nesting}.
\epf

\begin{rmk}
Using induction it follows that the map 
\eqref{eq: star equality}
used in the proof of 
Theorem \ref{thm: basic isom} 
is given explicitly by 
$$\{(\wh y^i,x_i)\}\mapsto (h,S_1,S_2,y,x)\in \cA(W)\glue{}\cB(W,V)$$ 
where $y,x$ have components $y^i,x_i$ resp. and 
$y^i = \wh y^i T_{i-1}\cdots T_1$ with $T_i = 1+x_i\wh y^i$ 
(as in \S\ref{ssn: lin alg})
so that $1+yx = u_-^{-1}hu_+=hS_2S_1$. 
\end{rmk}

Now suppose further that we also have an ordered grading $V=\bigoplus_1^r V_i$ of $V$.
Then a more symmetrical generalisation of the previous result is available, since we may apply fission operators to both sides of $\cB(W,V)$. %
 Let $H(V) = \Prod \GL(V_i)$, 
$H(W) = \Prod \GL(W_i)$ and $G=\GL(V\oplus W)$. 
From now on we will use the following notation for multiple fusions:
$$\bigoasterisk^{\mapleft{}}_{V}\cB(W_i,V) \ :=  \ 
\cB(W_s,V) \fusion{V}\cdots \fusion{V}\cB(W_1,V).$$

\begin{cor}
The following three quasi-Hamiltonian $H(W)\times H(V)$ spaces are isomorphic

\beq \label{eq: isom1.0}
\bigl(\cA(W)\glue{W}\cA^2(W,V)  \glue V \cA(V)\bigr)\spq G, 
\eeq

\beq \label{eq: isom1.1}
\left(\bigoasterisk^{s\mapleft{}1}_{V}\cB(W_i,V)\right)
\glue{V} \cA(V), %
\eeq

\beq \label{eq: isom1.2}
\left(\bigoasterisk^{r\mapleft{}1}_{W}\cB(V_i,W)\right)
\glue{W} \cA(W). %
\eeq

\end{cor}

Note that \eqref{eq: isom1.0} is isomorphic to 
$\cA(W)\glue{}\cB(W,V)  \glue{} \cA(V)$. 

\pf
An isomorphism between \eqref{eq: isom1.0} and \eqref{eq: isom1.1}
is immediate by upon 
applying $(\,\cdot\,)\glue{}\cA(V)$ 
to both sides of \eqref{eq: star equality}.  
An isomorphism between  \eqref{eq: isom1.0} and \eqref{eq: isom1.2}
arises similarly once we identify  $\cA^2(W,V)$ with $\cA^2(V,W)$ in \eqref{eq: isom1.0},  using Proposition 
\ref{prop: imd to change parab gln}.
\epf

In terms of graphs this gives three descriptions 
of the space 
$\Rep^*(\Ga(s,r),V\oplus W)$
for the complete bipartite graphs $\Ga(s,r)$.
Unlike the star-shaped case (with $r=1$), we now start to get spaces 
which are different to those obtained by fusing together the Van den Bergh spaces along the edges of the graph
(e.g. the reductions are different even if $r=s=2$ when $\Ga$ is the square, as mentioned in Remark \ref{rmk: square example etc}.)

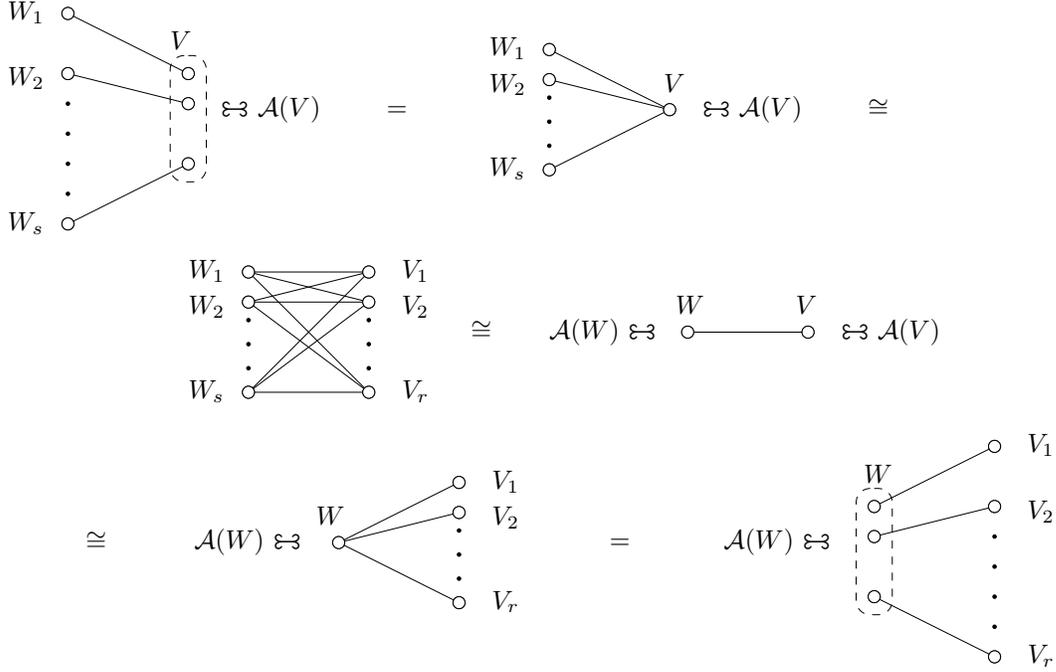
\begin{figure}[ht]
	\centering
	\input{second.isom.v3.pstex_t}
	\caption{Three ways to construct a complete bipartite graph.}\label{fig: second isom}
\end{figure}

The next generalisation is to replace the pair of ordered graded vector spaces $V,W$ appearing above by an arbitrary 
 $k$-tuple of ordered graded vector spaces.
First we will set up some new notation 
(beware of the unfortunate conflicts with that above).
Suppose $J=\{1,\ldots,k\}$ and we have a finite totally ordered 
set 
$I_j$ for each $j\in J$.
Let $I=\bigsqcup I_j$ be the disjoint union of these sets and suppose we have chosen a finite dimensional complex vector space 
$V_i$ for each $i\in I$.
Then we may define ordered graded vector spaces
$$W_j = \bigoplus_{i\in I_j} V_i, \qquad 
V  = \bigoplus_{j\in J} W_j,   \qquad 
U_j = V\ominus W_j = \bigoplus_{p\in J\setminus\{j\}} W_{p}.$$
Thus we view $V$ as graded by $J$ (and not by $I$) and similarly $U_j$ is graded by $J\setminus\{j\}$. 
Define 
$G= \GL(V)$ and
$H=\prod_{i\in I}\GL(V_i)$.

\begin{thm}\label{thm: main isom}
For any fixed integer $j\in\{1,\ldots,k\}$ the space

$$\left(\cA^2(V) 
\glue{W_1} \cA(W_1) 
\glue{W_2} \cA(W_2) %
\cdots
\glue{W_k} \cA(W_k) %
\right)\spq G$$

\noindent
is isomorphic to the space
$$
\left(
\bigoasterisk^{\mapleft{i\in I_j}}_{U_j}\cB(V_i,U_j)
\right)
\glue{U_j} \cA^2(U_j)  
\glue{W_1} \cA(W_1) 
\glue{W_2}
 \cA(W_2) %
\ \cdots\left(\begin{smallmatrix}\text{omitting} \\ \text{$j$th term}\end{smallmatrix}\right)\cdots
\glue{W_k} 
\cA(W_k)$$

\noindent
as a quasi-Hamiltonian $H$-space.
\end{thm}

Note that the first space here is 
$ \cB(W)\glue{K} \Prod_j\cA(W_j)$, in the notation of Prop. 
\ref{prop: open subset of reps}.

\pf
We will show this for $j=1$ since the other cases are similar 
(note that up to isomorphism the space $\cA^2(V)$ does not depend on the ordering of $J$---this follows from the isomonodromy isomorphisms of Proposition \ref{prop: imd to change parab gln}). 
First recall that
$\cB(W_1,\ldots,W_k) = \cA^2(V) \spq G,$
and let $\cB^{\!\!\fusion{}}_1$ denote the fusion 
$$ \bigoasterisk^{\mapleft{i\in I_1}}_{U_1}\cB(V_i,U_1).$$
Then the result will follow if we can show
\beq\label{eq: key step}
\cB(W_1,\ldots,W_k)\glue{W_1} \cA(W_1)\  \cong \ %
\cB^{\!\!\fusion{}}_1
\glue{U_1} \cA^2(U_1)   
\eeq
as quasi-Hamiltonian spaces,  
since, if so, we can just glue  $\cA(W_j)$ %
on to both sides for $j=2,\ldots,k$ to obtain the desired result.
Now observe Theorem \ref{thm: basic isom} implies
$$\cB(W_1,U_1)\glue{W_1} \cA(W_1)
\quad \cong \quad
\cB^{\!\!\fusion{}}_1. 
$$
If  we apply the higher fission operator $(\ \cdot\ )\glue{U_1} \cA^2(U_1)$ to both sides of this we see that \eqref{eq: key step} (and thus the theorem) will follow provided that 
$$\cB(W_1,U_1) \glue{U_1} \cA^2(U_1) \quad  \cong\quad 
\cB(W_1,\ldots,W_k).$$
In turn this will follow if there is an isomorphism
$$
\cA^2(W_1,U_1)
 \glue{U_1} \cA^2(U_1) \quad  \cong\quad \cA^2(V)$$
upon reducing by $G$. %
But this follows from the nesting result Proposition \ref{prop: gl nesting}.
\epf

In terms of graphs, this gives $k+1$ descriptions of the space 
$\Rep^*(\Ga,V)$ for any complete $k$-partite graph $\Ga$.

The higher fission operator $(\,\cdot\,)\glue{} \cA^2(V)$
is the multiplicative analogue of the $1$-fission operation of \cite{rsode} Appendix C:  In terms of graphs, if $V=\bigoplus_1^k V_i$, this operation breaks up a node with vector space $V$ into $k$ pieces (and repeats the original edges to each of the new nodes) 
and then each pair of the $k$ new 
nodes is connected by a single edge.
For example the key isomorphism 
$
\cB^{\!\!\fusion{}}_1
\glue{} \cA^2(U_1)\cong \cB(W_1,\ldots,W_k)\glue{} \cA(W_1) $  
of  \eqref{eq: key step} may be pictured
in terms of graphs as follows 
(in the case where  $W_1=V_1\oplus V_2$ and $k=3$, so $U_1=W_2\oplus W_3$):

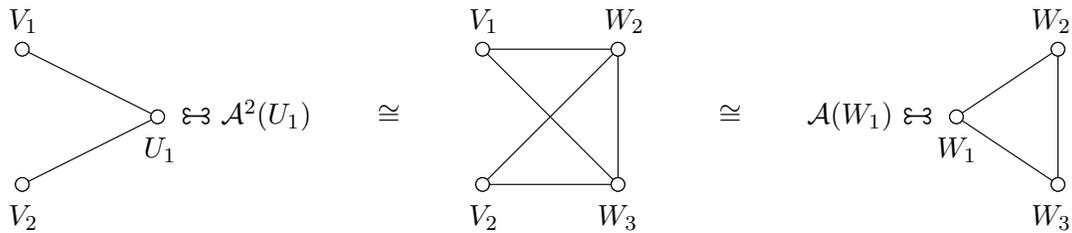
\begin{figure}[ht]
	\input{hfission4.pstex_t}
	\caption{Example isomorphism involving a higher fission operator $(\,\cdot\,)\glue{}\cA^2$}
\end{figure}

\subsection{Simplest isomorphism}\label{ssn: s=2 case}

This section gives more details of the simplest case of Theorem \ref{thm: basic isom}, i.e. the statement that 
\beq \label{eq: s=2 isom}
\cB(W_2,V)\fusion{V}\cB(W_1,V)\  \cong \ 
\cB(W,V)\glue{W} \cA(W_1,W_2)
\eeq
as quasi-Hamiltonian $\GL(V)\times\GL(W_1)\times\GL(W_2)$ spaces, if $W=W_1\oplus W_2$.

Suppose $x:W\to V$ and $y:V\to W$ represent a point $(y,x)$ of $\cB(W,V)$, so the moment maps are 
$((1+yx)^{-1}, 1+xy)\in \GL(W)\times\GL(V).$
Thus a point of 
$\cB(W,V)\glue{W}  \cA(W_1,W_2) $ 
is given by a tuple $y, x, S_1, S_2 ,h$ such that
\beq\label{eq: s=2 big cell}
1+yx = hS_2S_1\in \GL(W). %
\eeq
It will be convenient below to change notation slightly and set 
$$u_+=S_1,\qquad u_- = hS_2^{-1}h^{-1}$$ so that  $hS_2S_1 = u_-^{-1}hu_+$ 
(the quasi-Hamiltonian two-form is written in these coordinates in \cite{fission}).
Now let  $\pi_i:W\to W_i$ and $\iota_i:W_i\to W$ be the inclusion and projection between $W$ and its summand $W_i$, and define
$$ x_i = x\circ \iota_i : W_i\to V,\qquad
y^i  = \pi_i\circ y: V\to W_i$$
$$\varphi_0=1,\  \varphi_1 = 1 + x_1 y^1 \in \Aut(V),\qquad
\wh y^i  = y^i \circ \varphi^{-1}_{i-1} : V\to W_i,$$
$$h_i = 1+\wh y^ix_i\in \End(W_i), \qquad
T_i = 1 + x_i \wh y^i\in \End(V)$$
(Note that the existence of the decomposition \eqref{eq: s=2 big cell} is equivalent to $\varphi_1$ being  invertible.)
As in \S\ref{ssn: lin alg} we have that
\beq \label{eq: s=2 2nd gluing}
1+xy = T_2T_1 \in \GL(V).
\eeq
Thus we are led to take the corresponding point of $\cB(W_2,V)\fusion{V}\cB(W_1,V)$
to be $((\wh y^2,x_2),(\wh y^1,x_1))$.
This prescription defines an isomorphism between the two sides of 
\eqref{eq: s=2 isom}. 
It is clearly equivariant under the action of 
$\GL(W_1)\times \GL(W_2)\times \GL(V)$
and the moment maps match up due to \eqref{eq: s=2 2nd gluing} and since $h_i^{-1}$ is the $i$th component of $h^{-1}$.
Finally we need to check the quasi-Hamiltonian two-forms match up. This comes down to verifying the identity: 
\begin{align*}
\Tr(1+yx)^{-1}dy&\wedge dx - \Tr(1+xy)^{-1}dx\wedge dy
+\Tr (\bar\cU_- h\bar \cU_+ h^{-1})\\
&=\Tr  h_1^{-1}d\wh y^1\wedge dx_1 - \Tr T_1^{-1}dx_1\wedge d\wh y^1\\
&+\Tr  h_2^{-1}d\wh y^2\wedge dx_2 - \Tr T_2^{-1}dx_2\wedge d\wh y^2
-\Tr  T_1^{-1}T_2^{-1} dT_2\wedge dT_1
\end{align*}
where $\bar\cU_- = (du_-) u_-^{-1}$ and  
$\bar\cU_+ =  (du_+)u_+^{-1}$.
This will be left as an exercise (the author has checked it symbolically using a computer program).

\subsection{Alternative version.}

Suppose instead we set $M_i = 1+\wh x_i y^i\in \Aut(V)$ where 
$\wh x_i = \varphi^{-1}_{i-1}\circ x_i$. 
Then
\beq \label{eq: s=2 2nd gluing v2}
1+xy  = M_1M_2
\eeq
so we are led to consider the point 
$((y^1,\wh x_1),(y^2,\wh x_2))$
of
$\cB(W_1,V)\fusion{V}\cB(W_2,V)$.
This prescription defines an isomorphism
\beq \label{eq: s=2 isom v2}
\cB(W_1,V)\fusion{V}\cB(W_2,V)\  \cong \ 
\cB(W,V)\glue{W} \cA(W_1, W_2).
\eeq
It is clearly equivariant under the action of 
$\GL(W_1)\times \GL(W_2)\times \GL(V)$.
The moment maps match up 
due to \eqref{eq: s=2 2nd gluing v2} and since $h_i^{-1}$ is the $i$th component of $h^{-1}$ (and $h_i$ is dual to both $M_i$ and $T_i$).
The difference between the above two descriptions is just swapping the order of the two factors. Noting that $\varphi_0=1$ and $\varphi_1 = T_1 = M_1$ the corresponding isomorphism
\beq \label{eq: com isom}
\cB(W_1,V)\fusion{V}\cB(W_2,V)\  \cong \ 
\cB(W_2,V)\fusion{V}\cB(W_1,V) 
\eeq
is just 
$((y^1,x_1),(y^2,\wh x_2))\mapsto ((\wh y^2,x_2),(y^1,x_1))$ and $(\wh y^2,x_2) = T_1\cdot(y^2,\wh x_2)$
where the dot denotes the action of $\GL(V)$, i.e. it is the braid isomorphism of \cite{AMM} Theorem 6.2.
In particular the fact that \eqref{eq: s=2 isom v2} is a quasi-Hamiltonian isomorphism follows since
\eqref{eq: s=2 isom}  and \eqref{eq: com isom} are quasi-Hamiltonian isomorphisms.

\section{Wild character varieties}\label{sn: wcvs}

The wild character varieties  are a large class of algebraic symplectic/Poisson varieties, which parameterise %
 connections on bundles on curves.
They generalise the (tame) character varieties, i.e. spaces of fundamental group representations, which parameterise regular singular connections on bundles on curves. 
The first approach to their symplectic structures was analytic \cite{thesis, smid}, in the style of Atiyah--Bott. Subsequently a purely algebraic approach was developed extending the quasi-Hamiltonian theory \cite{saqh02, saqh, fission, gbs}. 

The general set-up is as follows (see \cite{gbs} for more details).
Fix a connected complex reductive group $G$, such as $G=\GL_n(\IC)$, and a maximal torus $T\subset G$.
A wild character variety $\MB(\Si)$ 
is then associated to an object called an ``irregular curve'' $\Si$. 
This generalises the notion of curve with marked points, in order to encode some of the boundary conditions for irregular connections:
An irregular curve $\Si$ (for fixed $G$) is a triple $(\Si,\ba,\bQ)$ consisting of a smooth compact complex algebraic curve $\Si$, plus some marked points 
$\ba=(a_1,\ldots,a_m)\subset \Si$, for some $m\ge 1$, and an ``irregular type'' $Q_i$ at each marked point $a_i$.
The extra data   is the irregular type: an irregular type at a point 
$a\in \Si$ is an arbitrary element 
$Q\in \lt(\wh \cK_a)/\lt(\wh \cO_a)$ where $\lt=\Lie(T)$.
If we choose a local coordinate $z$ vanishing at $a$ then 
$Q\in \lt\flp z \frp/\lt\flb z \frb$
and we may write
$$Q = \frac{A_r}{z^r}+\cdots + \frac{A_1}{z}$$ 
for  some elements $A_i\in \lt$, for some integer $r\ge 1$.
(The tame/regular singular case is the case when each $Q_i=0$, and in this case the wild character variety coincides with the usual space of fundamental group representations.)
There are several variations of this definition 
(bare irregular curves, twisted irregular curves) discussed in \cite{gbs} but they will not be needed here.

Given an irregular curve $\Si=(\Si,\ba,\bQ)$ we may consider algebraic connections on algebraic $G$-bundles on $\Si^\circ = \Si\setminus\{\ba\}$
such that near each puncture $a_i$ there is a local trivialisation such that the connection takes the form $dQ_i + \text{logarithmic terms}$, i.e. whose irregular part is determined up to isomorphism by $dQ_i$.
One way to state the irregular Riemann--Hilbert correspondence in this context is that: 
the category (groupoid) of such connections is equivalent to the category 
of {\em Stokes $G$-local systems} determined by $\Si$ (cf. \cite{gbs} Appendix A).
This generalises %
Deligne's equivalence \cite{Del70}  between regular singular connections on algebraic vector bundles on 
$\Si^\circ$ and local systems on $\Si^\circ$.

A Stokes $G$-local system is defined as follows. Given the irregular 
curve $\Si$ we may define two real surfaces with boundary 
$$\wt\Si\hookrightarrow \wh\Sigma \twoheadrightarrow  \Si$$ as follows:
$\wh\Si$ is the real oriented blow-up of $\Si$ at the points of $\ba$, so that each point $a_i$ is replaced by a circle $\partial_i$ parameterising the real oriented direction emanating from $a_i$, and the boundary of  $\wh \Si$ is $\partial\wh \Si=\partial_1\sqcup\cdots \sqcup \partial_m$.
In turn each irregular type $Q_i$ canonically determines three pieces of data (see \cite{gbs} for full details):

1) a connected reductive group $H_i\subset G$, the centraliser of $Q_i$ in $G$,

2) a finite set $\IA_i\subset \partial_i$ of singular/anti-Stokes directions for all $i=1,\ldots,m$, and

3) a unipotent group $\ISto_d\subset G$ normalised by $H_i$, the Stokes group of $d$, for each $d\in \IA_i$ for all $i=1,\ldots,m$.

\begin{figure}[ht]
	\centering
	\input{halo.on.curve2.pstex_t}
	\caption{The surface $\wt \Si$ with halo drawn}\label{fig: halo}
\end{figure}
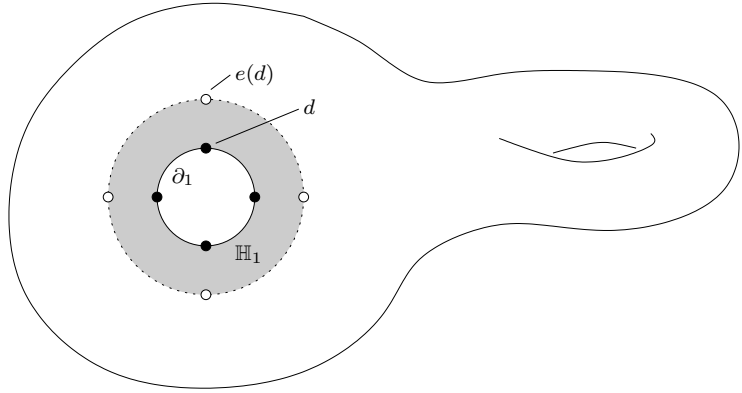

The surface $\wt \Si\subset \wh \Si$ is defined by puncturing $\wh \Si$ once at a point $e(d)$ in its interior near each singular direction $d\in \IA_i\subset \partial_i$, for all $i=1,\ldots,m$.
For example we could fix a small tubular neighbourhood of $\partial_i$ (a ``halo'')  $\IH_i\subset \wh \Si$  (so $\IH_i$ an annulus), and 
choose the extra puncture $e(d)\in\wh \Si$ to be on the {\em interior} boundary of $\IH_i$ near the singular direction $d$ (as pictured in Figure \ref{fig: halo}).

A Stokes $G$-local system for the irregular curve $\Si$ consists of 
 a $G$-local system on $\wt \Si$, with a flat reduction to $H_i$ in $\IH_i$ for each $i=1,\ldots,m$, such that the local monodromy around $e(d)$ is in 
$\ISto_d$ for any basepoint in $\IH_i$, for all $d\in \IA_i$ for all $i$.

Thus, via the irregular Riemann--Hilbert correspondence, the classification of connections with given irregular types is reduced to the classification of Stokes $G$-local systems, which goes as follows:
Choose a basepoint $b_i\in \IH_i$ for each $i$ and let 
$\Pi$ be the fundamental groupoid of $\wt\Si$ with basepoints 
$\{b_1,\ldots,b_m\}$.
Then we may consider the space $\Hom(\Pi,G)$ of $G$-representations of $\Pi$ 
and the subspace $\Hom_\IS(\Pi,G)$ of ``Stokes representations'' which are the representations such that the local monodromy around $\partial_i$ (based at $b_i$) is in $H_i$ and the local monodromy around $e(d)$ (based at $b_i$) is in $\ISto_d$,  for all $d\in\IA_i$, for all $i$.
The group $\bH:=H_1\times\cdots\times H_m$ acts naturally on $\Hom_\IS(\Pi,G)$ and the set of isomorphism classes of Stokes $G$-local systems for $\Si$ is in bijective correspondence with the set of $\bH$-orbits in $\Hom_\IS(\Pi,G)$.

The space of Stokes representations is naturally a smooth affine variety
and the  wild character variety $\MB(\Si)$ is defined to be the affine quotient
$$ \MB(\Si) = \Hom_\IS(\Pi,G)/\bH$$
taking the variety associated to the ring of $\bH$-invariant functions.
One of the main results of \cite{gbs} is that 
$\Hom_\IS(\Pi,G)$ is an algebraic quasi-Hamiltonian $\bH$-space, and this implies that  $\MB(\Si)$ is an algebraic Poisson variety, with symplectic leaves obtained by fixing 
the conjugacy classes of ``inner'' local monodromies, in $H_i$, of the Stokes local systems around $\partial_i$ for each $i$.
In full generality the space $\Hom_\IS(\Pi,G)$ of Stokes representations is identified (\cite{gbs} (36)) with the reduction
\beq\label{eq: big fusion prod}
\Hom_\IS(\Pi,G)\cong\left(
\ID^{\fus g}\fusion{G}\cA(Q_1)
\fusion{G}\cdots\fusion{G} \cA(Q_m)\right)\spq{} G\qquad\qquad\qquad\qquad\ \eeq
$$
\ \qquad\qquad\qquad\qquad\qquad\qquad\qquad\qquad\cong \left(\ID^{\fus g}\fusion{G}\cA(Q_1)
\fusion{G}\cdots\fusion{G} \cA(Q_{m-1})\right)\glue{G} \cA(Q_m)
$$ 

\noindent
involving the explicit quasi-Hamiltonian $G\times H_i$-spaces $\cA(Q_i)$ defined in \cite{gbs} Theorem 7.6. 
Here $\ID\cong G\times G$ is the internally fused double of \cite{AMM}, which is a quasi-Hamiltonian $G$ space with moment map given by the group commutator.
(In the tame  case $Q_i=0$, the space $\cA(Q_i)$ 
is just the double $D(G)\cong G\times G$, and if $Q_i$ has regular semisimple leading term then 
$\cA(Q_i)$ is the space $\wt \cC/L$ of \cite{saqh02, saqh} Remark 5.)
Given a conjugacy class $\bcC=\cC_1\times\cdots\times\cC_m\subset \bH$ there is a symplectic wild character variety
\beq\label{eq: swcv}
\MB(\Si,\bcC) = \mu_{\bH}^{-1}(\bcC)/\bH 
\subset  \MB(\Si)\eeq
obtained by fixing the inner local monodromy  conjugacy classes.
In turn we may consider the stable points (for the action of $\bH$ on $\Hom_\IS(\Pi,G)$) which may be characterised as the {\em irreducible} Stokes representations (see \cite{gbs} \S9) to define the subset of stable points
$$\MB^{st}(\Si,\bcC)\subset \MB(\Si,\bcC).$$
If the class $\bcC$ is sufficiently generic (see \cite{gbs} Corollaries 9.8, 9.9) then $\MB^{st}(\Si,\bcC)=\MB(\Si,\bcC)$ and it is a symplectic orbifold; in the case $G=\GL_n(\IC)$ it is smooth, and isomorphic to one of the complete \hk manifolds (wild Hitchin spaces) constructed in \cite{wnabh}---as mentioned in Rmk \ref{rmk: phi-mqv}, here we have set the Betti weights equal to zero.

\subsection{Type $\bf 3$ wild character varieties}
In this article we will focus on a special class of wild character varieties.
From the general viewpoint described above they appear to be very special, but from the point of view of quivers, we will see it is still a somewhat vast class.
Suppose $G=\GL(V)$ is a general linear group, with $V=\IC^n$.

\begin{defn}
An irregular curve $\Si$ is of {\em``type $3$''} if $\Si=(\IP^1,\infty,Q)$ where $Q$ has pole order $r \le 2$, i.e. the underlying curve is the Riemann sphere with just one marked point $\infty$ and $dQ$ has a pole of order at most  $3$.
\end{defn}

Thus if $\Si=(\IP^1,\infty,Q)$ is an irregular curve of type $3$
we may write
\beq\label{eq: type 3 irr type}
Q= \frac{A}{2}z^2+Tz = \frac{A}{2w^2}+\frac{T}{w}\qquad
\text{so that}\qquad dQ = (Az+T)dz
\eeq
for some diagonal matrices $A,T\in \End(V)$, where $w=1/z$ is a coordinate vanishing at $z=\infty$. Thus $dQ$ has a pole of order at most $3$ at $\infty$.

More generally we will consider irregular curves of {type $3+1^m$}, for integers $m\ge 1$: an irregular curve $\Si$ is of type $3+1^m$, if it is of type $3$, except it has $m$ other marked points 
$t_1,\ldots,t_m\in \IC=\IP^1\setminus \infty$ each with zero irregular type.

The next main result we will prove may be summarised as follows:

\begin{thm}\label{thm: reduce to three}
Suppose $G=\GL(V)$ is a general linear group and $\Si$ is an irregular curve of type $3+1^m$ for some $m\ge 1$.
Then for any choice of conjugacy classes $\bcC$ 
there is another vector space $\wh V$ and a type $3$ irregular curve 
$\wh\Si$
(for the group $\GL(\wh V)$) such that
the wild character variety
$$\MB^{st}(\Si,\bcC)$$
is isomorphic as an algebraic symplectic manifold to the type $3$ wild character variety
$$\MB^{st}(\wh\Si,\wh\bcC)$$
for some conjugacy classes $\wh \bcC$.
\end{thm}

In particular we could take $\Si$ to be tame (i.e. of type $1^m$) and still find it is isomorphic to a type $3$ wild character variety: any genus zero tame character variety is isomorphic to a type $3$ wild character variety.

The way we will prove this is via multiplicative quiver varieties. The main statements are the following.

\begin{thm}\label{thm: wcv2mqv}
Suppose $G=\GL(V)$ is a general linear group and $\Si$ is an irregular curve of type $3+1^m$ for some $m\ge 0$.
Then for any choice of conjugacy classes $\bcC$ 
there is a supernova graph $\Ga$ (with monochromatic core), and data $q,d$ such that
the wild character variety of $\Si$ with classes $\bcC$
is isomorphic  as an algebraic symplectic manifold to the multiplicative quiver variety of $\Ga$ with parameters
$q,d$:
$$\MB^{st}(\Si,\bcC)\cong\cM^{st}(\Ga,q,d).$$
\end{thm}

Conversely all such multiplicative quiver varieties arise as type $3$ 
wild character varieties:

\begin{thm}\label{thm: mqv2wcv}
Suppose $\Ga$ is a supernova graph (with monochromatic core) and $q,d$ are given.
Then there is a vector space $V$ and a type $3$ irregular curve $\Si$
(for the group $G=\GL(V)$) and conjugacy classes 
$\bcC$ such that the multiplicative quiver variety of $\Ga$
is isomorphic  as an algebraic symplectic manifold to the type $3$ wild character variety associated to $\Si,\bcC$:
$$\cM^{st}(\Ga,q,d) \cong \MB^{st}(\Si,\bcC).$$
\end{thm}

It is clear that these two results will imply Theorem \ref{thm: reduce to three}.
Further, as will be explained in \S\ref{sn: reflns}, the reflection isomorphisms of 
Theorem \ref{thm: refln isoms} follow almost directly 
from the passage between wild character varieties and multiplicative quiver varieties in Theorems \ref{thm: wcv2mqv}, \ref{thm: mqv2wcv}.
Theorems \ref{thm: wcv2mqv}, \ref{thm: mqv2wcv} will be proved in the next section.

\section{Multiplicative quiver varieties and wild character varieties}\label{sn: mqvs and wcvs}

\subsection{Type $\bf 3$ wild character varieties as multiplicative quiver varieties}\label{ssn: 3wcv as mqvs}

First we will establish Theorem \ref{thm: mqv2wcv} and the case $m=0$ of Theorem \ref{thm: wcv2mqv}.
The strategy is to consider an intermediate space 
$\Rep^*(\Ga,V){\spq{}}_{\bcC}\,\, H$ of reductions of an open multiplicative quiver variety at some conjugacy classes and then show the subspace of stable points of this is isomorphic both to a 
type $3$ wild character variety and to a multiplicative quiver variety.

Suppose $G=\GL(V)$ where $V=\IC^n$ and 
$\Si=(\IP^1,\infty,Q)$ is a type $3$ irregular curve with irregular type $Q=Az^2/2 + Tz$, where $A,T\in \End(V)$ are diagonal matrices.

Then $Q$ determines a complete $k$-partite graph $\Ga=\Ga(Q)$, the fission graph of $Q$, as recalled in \S\ref{sn: snova} (with $w=1/z$).
$\Ga(Q)$ may be defined as follows.
Let $J$ be the set of eigenspaces of $A$, and $k=\#J$.
Let
 $W_j\subset V$ be the corresponding eigenspace of $A$
for each $j\in J$,
and let $I_j$ be the set of eigenspaces of $T\bigl\vert_{W_j}$.
For any $i\in I_j$ let  $V_i\subset W_j\subset V$ be the corresponding eigenspace of $T\bigl\vert_{W_j}$.
Let $I=\bigsqcup I_j$ be the disjoint union of all the $I_j$.
Then $\Ga$ is the complete $k$-partite graph with nodes $I$ 
determined by this partition (i.e two nodes are joined by a single edge if{f} they are in different parts $I_j$ of $I$).
Further $Q$ determines a grading $V=\bigoplus V_i$ of $V$ by $I$.

Thus, viewing $\Ga$ as a simple (monochromatic) coloured quiver, for any choice of ordering there is a quasi-Hamiltonian $H$-space
$$\IM:=\Rep^*(\Ga,V)$$
as in Proposition \ref{prop: open subset of reps},
with moment map $\mu:\IM\to H$, where $H=\Prod\GL(V_i)$.
The stable points $\IM^{st}\subset \IM$ for the $H$-action consist of the 
irreducible (invertible) graph representations (as in 
Theorem \ref{thm: stable cmqv}), 
and given any conjugacy class $\bcC\subset H$ we may consider the 
stable points of the reduction at $\bcC$:
\beq\label{eq: redn of simple quiver}
\IM^{st}\spqa{\bcC} H = \{\rho\in \IM^{st} \st \mu(\rho)\in \bcC\}/H.
\eeq
This is an open subvariety of the affine variety associated to the ring of 
$H$-invariant functions on the affine variety $\mu^{-1}(\overline\bcC)$, 
where $\overline{\bcC}\subset H$ is the closure of the conjugacy class $\bcC$.

Note that the group $\bH$ attached to the irregular curve (the centralizer of $Q$) is the same as the group $H$ acting on $\IM$.

\begin{prop}\label{prop: wcv and omqv}
For any conjugacy class $\bcC\subset H$ the wild character variety 
$\MB^{st}(\Si,\bcC)$ is isomorphic to the reduction 
\eqref{eq: redn of simple quiver} of the open multiplicative quiver variety $\IM$.  %
\end{prop}
\pf
$\MB^{st}(\Si,\bcC)$ is the set of stable points of 
$\MB(\Si,\bcC)$. 
In turn the wild character variety $\MB(\Si,\bcC)$ 
is constructed in \cite{gbs} as the multiplicative symplectic quotient
of $\Hom_\IS(\Pi,G)$ by $H$ at the class $\bcC\subset H$, where $G=\GL(V)$.
Now recall from \eqref{eq: big fusion prod} (or \cite{gbs} (36)),
that $\Hom_\IS(\Pi,G)$
may be identified as an explicit quasi-Hamiltonian reduction.
In the current situation $g=0, m=1$ so the fusion product \eqref{eq: big fusion prod} just has one term $\cA(Q)$, where $Q=Az^2/2+Tz$ is the irregular type we fixed above.
The space $\cA(Q)$ attached to $Q$ was defined in 
\cite{gbs} Theorem 7.6, and was shown there to be a quasi-Hamiltonian $G\times H$-space\footnote{One can view $\cA(Q)$ as $\Hom_\IS(\Pi(\Si'),G)$
where $\Si'$ is either the irregular curve of type $3+1$ obtained 
by adding a single extra, tame, marked point to $\Si$ at $z=0$, 
or $\Si'$ is the analytic irregular curve obtained by taking a disk around 
$\infty$ in $\Si$.}. 
This was proved by showing that there was an isomorphism
$$\cA(Q)\  \cong\  \papt{G}{K} \glue{K} \pap{K}{H}$$
between $\cA(Q)$ and the gluing of two fission spaces
with respect to the intermediate group $K:=C_G(A)=\Prod\GL(W_j)$, so that 
$H\subset K\subset G$
(see \cite{gbs} Proposition 7.12, and proof of Lemma 7.11, where $K$ is denoted $H_2$).
(Beware that if $K\neq H$ then $\cA(Q)$ is not isomorphic to $\!\!\papt{G}{H}$, and will have different dimension.)
Thus we see
$$\Hom_\IS(\Pi,G) = \cB(Q) = \cA(Q)\spq{} G 
\cong (\papt{G}{K} \glue{K} \pap{K}{H})\spq{}G \ =\  
\cB(W)\glue{K} \pap{K}{H}$$
since, by definition $\cB(W) = \papt{G}{K}\spq{}G$, where $W=\bigoplus W_j$
(which is $V$ with a different grading).
Finally to identify $\Hom_\IS(\Pi,G)$ with $\Rep^*(\Ga,G)$ (as defined in 
Proposition \ref{prop: open subset of reps})
we just need to remark that 
$$\pap{K}{H} = \Prod_{j\in J} \cA(W_j)$$
since $K$ is a product, where $W_j$ is graded by $I_j$.

Thus the reduction \eqref{eq: redn of simple quiver} and the wild character variety $\MB^{st}(\Si,\bcC)$ are  the reduction of the {\em same} smooth affine variety
$$\Hom_\IS(\Pi,G) = \cB(Q) \ \cong\  \Rep^*(\Ga,V) = \IM$$
with the same structure of quasi-Hamiltonian $H$-space, by the same group $H$ at the same conjugacy class $\bcC\subset H$, and so are identified. (The stable points were defined intrinsically in terms of the $H$-action and so are identified on both sides.)
\epf

Now to identify the wild character variety with a multiplicative quiver variety we will explain how to 
encode the conjugacy class $\bcC$ by gluing on some legs to the fission graph 
$\Ga$, converting it into a supernova graph.

Since $\bcC\subset H=\Prod\GL(V_i)$ we may write 
$\bcC= \Prod\breve\cC_i$ where 
$\breve\cC_i$ is a conjugacy class of $\GL(V_i)$.
The extra data needed to determine a supernova graph is a ``marking'' of each
conjugacy class $\breve\cC_i$.

\begin{defn}\label{defn: marking}
Suppose  $\cC\subset \GL_n(\IC)$ is a conjugacy class.
A {\em `marking'} of $\cC$ is
a finite ordered set $(\xi_1,\xi_2,\ldots,\xi_w)$ of invertible complex numbers such that $\Prod_1^w(M-\xi_i) = 0$ for any $M\in \cC$.
\end{defn}

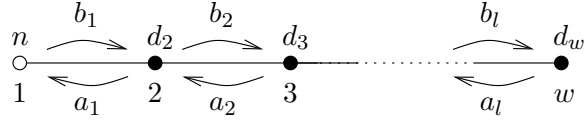
\begin{figure}[h]
	\centering
	\input{leg4.pstex_t}
	\caption{Representation of a type $A_w$ Dynkin graph.}\label{fig: leg}
\end{figure}

Equivalently a marking of $\cC$ is the choice of a monic annihilating polynomial $f\in \IC[x]$ with $f(0)\neq 0$, so that $f(M)=0$ for all $M\in \cC$, together with a choice of ordering of the multiset of roots of $f$.
A marking will be said to be {\em minimal} if $w=\deg(f)$ is minimal (so that $f$ is the minimal polynomial of $M$).
A marking is {\em special} if the first root is the identity ($\xi_1=1$).
Given a marking of $\cC\subset \GL_n(\IC)$, 
define invertible complex numbers 
\beq\label{eq: la defn}
q_i = \xi_{i}/\xi_{i-1}
\eeq
(including $q_1=\xi_1$) %
and integers 
$$d_i = \rank (M-\xi_1)\cdots (M-\xi_{i-1})$$ 
$i=1,2,\ldots$ (for any $M\in \cC$) so that $d_1=n$.
The marking and the dimensions $d_i$ are well-known to determine the class $\cC$. The following lemma gives a quiver-theoretic approach to this fact.
Consider the type $A_w$ Dynkin graph (a leg) with $w$ nodes and $l:=w-1$ edges, as in Figure \ref{fig: leg}.

\begin{lem} (cf. \cite{CB-Shaw}).\label{lem: orbits}  
If $\{(a_i,b_i)\}$ is an invertible representation of this leg 
(type $A_w$ Dynkin graph) on 
the vector space $V=\bigoplus_1^w \IC^{d_i}$ such that each 
$b_i$ is surjective and each $a_i$ is injective and the moment map conditions
$$ M = q_1(1+a_1b_1),\ 
1+b_1a_1 = q_2(1+a_2b_2),\ \ldots, \ 
1+b_{l-1}a_{l-1} = q_l(1+a_lb_l),\ 
1+b_la_l = q_w$$
hold, then $M\in \cC$.
\end{lem}
\pf
Writing $q_i=\xi_i/\xi_{i-1}$ and setting $\la_i=\xi_i-\xi_{i-1}, p_i=\xi_ib_i$, this reduces to the additive version, for adjoint orbits $\cO\subset \gl_n(\IC)$ (see \cite{gbs} Lemma 9.10).
Here is a proof, for completeness: 
It is clear that the conjugacy class of $M$ is uniquely determined by these conditions (cf. Proposition \ref{prop: class relations}).
Thus we just need to check that the class determined in this way is indeed the class $\cC$ we started with.
But if $M$ is any element of $\cC$ we may define an invertible 
representation satisfying these conditions by setting $V_1=\IC^n$ and then inductively 
$b_i=(M/\xi_i-1)\bigl\vert_{V_i}, V_{i+1}=\Im(b_i)$ and
$a_i$ to be the inclusion 
$V_{i+1}=
\Prod_1^i(M-\xi_j)(V)
\hookrightarrow \Prod_1^{i-1}(M-\xi_j)(V)= V_i$.
\epf

Note that if these moment map conditions hold and $M\in \cC$ then
\begin{align} \label{eq: detC}
\det(M) &= q_1^{d_1}\det(1+a_1b_1)\\
&=q_1^{d_1}\det(1+b_1a_1)=
q_1^{d_1}q_2^{d_2}\det(1+a_2b_2)=\cdots=\Prod_1^wq_i^{d_i}.\notag
\end{align}

Returning to the global picture, we will say that the type $3$ irregular curve $\Si$ is {\em marked} if we have chosen a conjugacy class 
$\breve\cC_i\subset \GL(V_i)$
for each $i$ (as above) together with a marking of each class $\breve\cC_i$.
Let $\IL_i$ denote the leg determined (as above) by the marking of 
$\breve\cC_i$.
Thus we may construct a larger graph $\wh \Ga$ by gluing the left-hand node of $\IL_i$ onto the node $i\in I$ of $\Ga$.
Let $\wh I$ denote the set of nodes of $\wh \Ga$ (equal to the disjoint union of the nodes of $\IL_i$ for all $i\in I$).
By definition $\wh \Ga$ is a simply-laced supernova graph.
Since each node of each leg comes with a vector space 
(i.e. $\IC^{d_i}$ in the above example) each node of  $\wh \Ga$ comes with a vector space, so we obtain an $\wh I$ graded vector space $\wh V$, and we will denote its graded pieces by $V_i$ for all $i\in \wh I$, and from now on we will let $d_i=\dim(V_i)$ for all $i\in \wh I$.
Similarly we obtain a scalar $q_i\in \IC^*$ for each $i\in \wh I$ from the scalars at each node of each leg.
We colour the graph $\wh \Ga$ so that the core $\Ga$ is monochromatic. 
Thus we have all the data required to define a (supernova) multiplicative quiver variety
$\cM(\wh \Ga,q,d)$.
Note that from \eqref{eq: detC}, it follows that
 $\Prod_I \det(M_i) = q^d$ (for any $M_i\in \breve\cC_i$), and both the quiver variety 
$\cM(\wh \Ga,q,d)$ and the 
wild character variety $\MB^{st}(\Si,\bcC)$ will be empty unless this common value is $1$.

\begin{thm}\label{thm: 3wcv as mqv}
The multiplicative quiver variety $\cM^{st}(\wh \Ga,q,d)$
is isomorphic to the type $3$ wild character variety $\MB^{st}(\Si,\bcC)$, as an algebraic symplectic manifold.
\end{thm}
\pf
By Proposition \ref{prop: wcv and omqv} it is sufficient to prove that $\cM^{st}(\wh \Ga,q,d)$ is isomorphic to $\IM^{st}\spq{}_{\bcC}\, H=
\Rep^*(\Ga)^{st}\spq{}_{\bcC}\, H$.
This result is the multiplicative analogue of 
\cite{slims} Theorem 9.11, and we have now set things up so that the same proof works verbatim (noting Lemma \ref{lem: stablem} above), with \cite{gbs} Lemma 9.10 replaced by Lemma \ref{lem: orbits} above.
In particular stability implies all the maps $b_i$ down the legs are surjective and all the maps $a_i$ up the legs are injective, and so the conjugacy class of the element $M_i=q_{i1}^{}(1+a_{i1}b^{}_{i1})$
(from the $i$th leg) is fixed. 
Given \cite{saqh, gbs}, identifying the symplectic structures comes down to the fact (\cite{yamakawa-mpa} Remark 4.2) that the quasi-Hamiltonian reduction of the representations of the leg $\IL_i$ gives the quasi-Hamiltonian form on the conjugacy class $\breve\cC_i$.
\epf

This establishes the $m=0$  case of Theorem \ref{thm: wcv2mqv}.
The above proof really shows the equivalence between the supernova multiplicative quiver varieties $\cM^{st}(\wh \Ga,q,d)$
and the reductions  $\IM^{st}\spq{}_{\bcC}\, H$, modulo the choice of a marking of each conjugacy class.
Thus (via Proposition \ref{prop: wcv and omqv}) this also establishes
Theorem \ref{thm: mqv2wcv}.

This shows how to ``read'' a wild character variety from a 
supernova multiplicative quiver variety. 
This is the ``generic reading'' (the other $k$ readings will be described in the following subsection).
In summary the basic data correspond as follow:

1) the rank $n$ of the structure group $G=\GL_n(\IC)$ of Stokes local systems 
appearing in the wild character variety is the sum of the dimensions on the nodes $I$ of the core $\Ga$ of the supernova graph $\wh \Ga$:
$n=\sum_{i\in I} d_i,$

2) If the core $\Ga\subset \wh\Ga$ is a complete $k$-partite graph, then $k$ is the number $\# J$ of eigenspaces $J$ of the leading term $A$ of the irregular type $Q$ at $\infty$ on the corresponding irregular curve (as in \eqref{eq: type 3 irr type}),

3) the number of nodes $I_j$ in the $j$th part of the core $\Ga$ is the number of eigenspaces of $T\bigl\vert_{W_j}\in\End(W_j)$ 
where $W_j\subset \IC^n$ is the eigenspace of $A$ corresponding to $j\in J$, and where $T$ is the next coefficient of the irregular type $Q$ (see \eqref{eq: type 3 irr type}),

4) The eigenspace $V_i$ of $T\bigl\vert_{W_j}\in\End(W_j)$  corresponding to 
$i\in I_j\subset I\subset \wh I$ has dimension $d_i$, and in turn $\dim(W_j) = \sum_{I_j} d_i$,

5) The (formal) monodromy conjugacy class $\bcC\subset H=\Prod_I \GL(V_i)$
of the Stokes local systems is a product of conjugacy classes 
$\breve\cC_i\subset \GL(V_i)$ for $i\in I$,
and $\breve\cC_i$ is the conjugacy class determined by the (data on the) $i$th leg of the supernova graph, as in Lemma \ref{lem: orbits}.

6) The coefficients $A,T$ of the irregular type $Q$
can be encoded as data on the graph as follows: The eigenvalues of $A$
correspond to assigning a distinct complex number for each part of the core
(i.e. to specifying an injective map $J\hookrightarrow \IC$). In turn the eigenvalues of $T$ correspond to choosing a scalar $t_i\in \IC$ for each node of $i\in I$ of the core, such that $t_{i_1}\neq t_{i_2}$ 
whenever $i_1\neq i_2$ are in the same part of $I$.

Finally recall that the passage from the wild character variety to the multiplicative quiver variety required the choice of a marking of each conjugacy class $\breve\cC_i$: changing the choice of marking immediately yields many isomorphisms between multiplicative quiver varieties:
This yields the reflection isomorphisms
for the nodes not in the core (see \S\ref{sn: reflns}).
To obtain the other reflection isomorphisms for the core nodes (for example for $i\in I_j\subset I$), in the 
following subsection we will use the main isomorphism
(Theorem \ref{thm: main isom}) to realize 
$\cM^{st}(\wh \Ga,q,d)$ as a wild character variety 
of type $3+1^m$, where $m=\#I_j$, parameterising certain 
Stokes local systems of smaller rank $\sum_{I\setminus I_j} d_i$. 
(This will also complete the proof of 
Theorem \ref{thm: wcv2mqv}.)

\subsection{Type $\bf 3+1^m$ wild character varieties as multiplicative quiver varieties}

Suppose $G=\GL(U)$ where $U=\IC^r$ and 
$\Si$ is a type $3+1^m$ irregular curve with underlying curve $\IP^1$, with irregular type $Q=Az^2/2 + Tz$ at $z=\infty$, where $A,T\in \End(U)$ are diagonal matrices, and there are $m\ge 1$ further marked points $t_1,\ldots,t_m\in \IC\subset \IP^1$ each with trivial irregular type.

Then $Q$ determines a complete $k'$-partite graph $\Ga'=\Ga(Q)$, the fission graph of $Q$, as in the previous subsection:
Let $J'$ be the set of eigenspaces of $A$, and $k'=\#J'$.
Since the eigenspaces are parameterised by the eigenvalues of $A$ we can (and will) identify $J'$ with a subset of $\IC$. 
Let
$W_j\subset U$ be the corresponding eigenspace of $A$
for each $j\in J'$,
and let $I_j$ be the set of eigenspaces of $T\bigl\vert_{W_j}$.
For any $i\in I_j$ let  $V_i\subset W_j\subset U$ be the corresponding eigenspace of $T\bigl\vert_{W_j}$.
Let $I'=\bigsqcup I_j$ be the disjoint union of all the $I_j$.
Then $\Ga'$ is the complete $k'$-partite graph with nodes $I'$ 
determined by this partition (i.e two nodes are joined by a single edge if{f} they are in different parts $I_j$ of $I'$).
Further $Q$ determines a grading $U=\bigoplus V_i$ of $U$ by $I'$.

Up to a slight change of notation all this is as before. 
Now we set $k=k'+1$ and construct a complete $k$-partite graph $\Ga$ 
by adding an extra part of size $m$ (the number of extra marked points):
Let $J=J'\sqcup\{\infty\}$, let $I_\infty=\{t_1,\ldots,t_m\}$,
let $I=\bigsqcup_J I_j$ and let
$\Ga$ be the complete $k$-partite graph determined by this 
partition of $I$ into $k$ parts.

As in \eqref{eq: swcv},
the extra data needed to determine a symplectic wild character variety
for the irregular curve $\Si$ consists of a conjugacy class $\bcC$ of the group $G^m\times\prod_{I'} \GL(V_i)$, i.e. a conjugacy class 
$\cC_i\subset G=\GL(U)$ for each $i\in I_\infty$ and a class 
$\breve\cC_i\subset \GL(V_i)$ for each $i\in I'$.

Further, we will say the type $3+1^m$ irregular curve $\Si$ is {\em marked} if we have chosen conjugacy classes $\bcC$ (as above) plus 
a marking of each class $\cC_i,\breve \cC_i$, such that the markings 
of each $\cC_i$ are special (for $i\in I_\infty$), 
cf. Definition \ref{defn: marking}.
The extra choice of marking is enough to 
 determine a supernova graph with core $\Ga$, as follows.
Let $\cC'_i\subset \GL(U)$ denote the inverse conjugacy class of $\cC_i$ (for  $i\in I_\infty$).
By inverting each $\xi$, the marking of $\cC_i$ determines a marking of $\cC'_i$ (for $i\in I_\infty$).
As in Lemma \ref{lem: orbits} each choice of marking determines a leg: 
Let $\wh \IL_i$ denote the leg determined by the marking of the inverse class $\cC'_i$ ($i\in I_\infty$) and let $\IL_i$ denote the leg of $\breve\cC_i$ 
($i\in I'$).
As usual note that each node of each leg comes with a parameter in $\IC^*$ 
and a dimension in $\IZ_{\ge 0}$.

Now, for each $i\in I_\infty$ let
$\IL_i$ denote the leg obtained by removing the left most edge of $\wh \IL_i$ (i.e. we remove the node with dimension $\dim(U)$).
Now we have a leg $\IL_i$ for any $i\in I=I'\sqcup I_\infty$, and we define a supernova graph $\wh \Ga$ 
by gluing the left-most node of $\IL_i$ 
onto the node $i$ of $\Ga$ for all $i\in I$.
Let $\wh I$ denote the set of nodes of $\wh \Ga$, and so from the 
parameters and the dimensions on the nodes of the legs we obtain 
$q\in (\IC^*)^{\wh I}$ and $d\in \IZ^{\wh I}$ (with components $d_i\ge 0$),
and thus may consider the multiplicative quiver variety 
$\cM(\wh\Ga,q,d)$, where $\wh \Ga$ is coloured so its core $\Ga$ is monochromatic.

\begin{thm}\label{thm: main reduced isom}
The multiplicative quiver variety $\cM^{st}(\wh \Ga,q,d)$
is isomorphic to the type $3+1^m$ wild character variety $\MB^{st}(\Si,\bcC)$, as an algebraic symplectic manifold.
\end{thm}
\pf
For $i\in I_\infty$ let $V_i=\IC^{d_i}$ and let 
$\breve \cC_i\subset \GL(V_i)$ denote the %
conjugacy class determined by the leg $\IL_i$
(here $d_i$ is the dimension of the left-most node of $\IL_i$, which equals $\rank(M-\id_U)$ for any $M\in \cC_i$).
Thus we now have a class $\breve \cC_i\subset \GL(V_i)$ for
all $i\in I$. 
Let $\breve\bcC\subset H=\Prod_I\GL(V_i)$ denote this collection of conjugacy classes.
Let $V=\bigoplus_I V_i$. Thus $V = U\oplus W_\infty$ where 
$W_\infty = \bigoplus_{I_\infty} V_i$, but we will always view $V$ as graded by $I$.
We will break the proof in several steps.

1) First, 
via  Proposition \ref{prop: wcv and omqv} and 
Theorem \ref{thm: 3wcv as mqv}
 the 
multiplicative quiver variety $\cM^{st}(\wh \Ga,q,d)$
may be identified with $\IM^{st}\spq_{\breve\bcC} H$
where $\IM=\Rep^*(\Ga,V)$ with $V=\bigoplus_I V_i$
(this is really what is shown in the proof of Theorem \ref{thm: 3wcv as mqv}).

2) The main isomorphism, Theorem \ref{thm: main isom}, then shows that
$\IM=\Rep^*(\Ga,V)$
is isomorphic as a (smooth, affine) quasi-Hamiltonian $H$-space to  
\beq\label{eq: bigfusion1}
\IM\cong \left(\,
\bigoasterisk^{\mapleft{i\in I_\infty}}_{U}\cB(V_i,U)
\right)
\glue{U} \cA(Q).
\eeq
In particular the stable points (for the action of $H$) are identified.
Here we have identified 
$$\cA(Q) \cong \cA^2(U)  
\glue{W_1} \cA(W_1) 
\glue{W_2}\cA(W_2)
\ \cdots
\glue{W_{k'}} \cA(W_{k'})
$$

\noindent
as quasi-Hamiltonian $G\times H'$ spaces, as in  Proposition \ref{prop: wcv and omqv} (but here $Q$ takes values in the Cartan subalgebra of $\gl(U)$),
where $J'$ is identified with $\{1,\ldots,k'\}$, 
$W_j=\bigoplus_{I_j} V_i$, $G=\GL(U)$ and $H'=\Prod_{I'} \GL(V_i)$,
so that $H=H'\times H_\infty$, where $H_\infty = \Prod_{I_\infty} \GL(V_i)$.

3) On the other hand, the wild character variety $\MB^{st}(\Si,\bcC)$ 
is constructed in \cite{gbs} as the quasi-Hamiltonian reduction of
\beq\label{eq: prewcv}
\left(\cC_m\fusion{G}\cdots\cC_2\fusion{G}\cC_1\glue{G} \cA(Q)\right)^{st}
\eeq
by $H'$ at the conjugacy class $\breve\cC'$.
More precisely the general construction of  \cite{gbs} (cf. \eqref{eq: big fusion prod} above)
identifies $\MB^{st}(\Si,\bcC)$ with the symplectic leaf of 
$$(D(G)^{\fus m}\glue{G} \cA(Q))^{st} / (G^m\times H')$$
corresponding to $\bcC\subset G^m\times H'$, since at a tame point  (with zero irregular type $Q_i$) the fission space $\cA(Q_i)$ equals the double $D(G)\cong G\times G$.
Since the action of $G^m$ is free (and $D(G)/G\cong G$)
we may identify  the stable points 
$(D(G)^{\fus m}\glue{} \cA(Q))^{st}$ (for the action of $G^m \times H'$)
with 
the stable points of $(G^m\glue{} \cA(Q))^{st}$ (for the $H'$ action).
In turn the space \eqref{eq: prewcv} is defined to be the subset of 
$(G^m\glue{} \cA(Q))^{st}$
 obtained by restricting to $\cC_m\fus\cdots\fus\cC_1\subset G^m$.
It is an open subset of the affine variety $\overline\cC_m\fus\cdots\fus\overline\cC_1\glue{}\cA(Q)$.
This now looks quite close to \eqref{eq: bigfusion1}.

4) Let $\mu_\infty:\IM\to H_\infty$ be the $H_\infty$ component of the moment map on \eqref{eq: bigfusion1}, and let 
$\breve \bcC_\infty=\Prod_{I_\infty}\breve \cC_i\subset H_\infty$.
Consider the space 
$$\mu_\infty^{-1}(\breve \bcC_\infty)^{st} := 
\mu_\infty^{-1}(\breve \bcC_\infty) \cap \IM^{st}$$
where $\IM^{st}$ denotes the stable points for the full $H$ action (i.e. the irreducible invertible representations of the graph $\Ga$).
To complete the proof we claim 

\begin{prop}\label{prop: stable isom}
There is a well-defined map
from $\mu_\infty^{-1}(\breve \bcC_\infty)^{st}/H_\infty$
to the space \eqref{eq: prewcv} which is an isomorphism
of quasi-Hamiltonian $H'$-spaces.
\end{prop}
\pf
Via the isomorphism  \eqref{eq: bigfusion1} a graph representation 
$\rho\in \Rep^*(\Ga,V)$ determines elements
$$(a_i,b_i)\in  \cB(V_i,U)$$
for any $i\in I_\infty$, so $a_i:U\to V_i$ and 
$b_i:V_i\to U$.
If $\rho$ is stable then $a_i$ is surjective and $b_i$ is injective:
via Lemma \ref{lem: stablem} $\rho$ is irreducible as a graph representation,
and $a_i$ encodes all the maps to $V_i$, and $b_i$ encodes all the maps out of $V_i$, and so if $a_i$ was not surjective then we could replace $V_i$ by the image of $a_i$ to obtain a subrepresentation, and similarly if $b_i$ was not injective, then we could replace $V$ by the kernel of $b_i$.
Now the $\GL(V_i)$ component of $\mu_\infty$ is $(1+a_ib_i)^{-1}$, and so 
fixing this to be in $\breve \cC_i$ is equivalent to fixing 
$1+b_ia_i$ to be in $\cC_i\subset \GL(U)$, by Lemma \ref{lem: orbits}
(this is why we attached $\wh \IL_i$ to the inverse of $\cC_i$ in the preparation for the theorem).
Thus there is a map from  
$$\mu_\infty^{-1}(\breve \bcC_\infty)^{st}\subset 
\left(\,
\bigoasterisk^{\mapleft{i\in I_\infty}}_{U}\cB(V_i,U)
\right)
\glue{U} \cA(Q)\qquad
\text{to}\qquad 
\cC_m\fusion{}\cdots\cC_2\fusion{}\cC_1\glue{U} \cA(Q)$$
 taking 
$(a_i,b_i)\in \cB(V_i,U)$ to $1+b_ia_i\in \cC_i\subset G=\GL(U)$.
Due to the injectivity/surjectivity conditions
this map has
 fibres exactly the $H_\infty$-orbits.
It is clearly $H'$-equivariant and the reader may readily verify, by comparing exactly what the stability conditions say, that it is 
an isomorphism onto the stable subset \eqref{eq: prewcv}:
In more detail the stable points of \eqref{eq: prewcv} may be described as follows.
A point of \eqref{eq: prewcv} is determined by matrices $T_i\in \cC_i$
and a representation $\rho'\in \Rep(\Ga',U)$ of $\Ga'$ on 
$U=\bigoplus_{I'} V_i$. The representation $\rho'$ determines elements
$v_\pm \in U_\pm\subset \GL(U)$
(cf. \eqref{eq: vpm from rho}), and  
these data should satisfy the condition
\beq\label{eq: explicit full space}
T_m\cdots T_1 v_-v_+ = w_+gw_-\in \GL(U)
\eeq
for some $w_\pm\in U_\pm,g\in H'$.
This point of $\cC_m\fus\cdots\fus\cC_1\glue{}\cA(Q)$ is stable 
if{f} there is no proper nontrivial subrepresentation $W\subset U$ of the graph representation $\rho'$ such that $T_i(W)\subset W$ for all $i$. 
\epf

This completes the proof of the theorem since both sides are now the reduction by $H'$ of the same space.
\epf

In the next section this will be used to deduce the desired Weyl group isomorphisms, since exactly the same class of multiplicative quiver varieties arises from both the type $3$ and the type $3+1^m$ wild character varieties.
First we will give a corollary, summarise the dictionary and give examples.

\begin{cor}
If $\Ga$ is a simply-laced supernova graph, with monochromatic core, then the multiplicative quiver variety 
$\cM^{st}(\Ga,q,d)$ is a \hk manifold.
\end{cor}
\pf
From Theorem \ref{thm: 3wcv as mqv}, $\cM^{st}(\Ga,q,d)$ is a type $3$ wild character variety $\MB^{st}(\Si,\bcC)$. 
Then the irregular Riemann--Hilbert correspondence identifies 
this wild character variety with one of the spaces of stable meromorphic connections shown to be \hk in \cite{wnabh}. 
As mentioned in Rmk \ref{rmk: phi-mqv}, here we have set the Betti weights equal to zero.
\epf

\begin{rmk}
We conjecture that the other $k$ (non-generic) readings of $\Ga$ yield isometric \hk manifolds (this should follow from an extension of the results of \cite{szabo-nahm} on the Nahm transform).
\end{rmk}

\subsection{Dictionary}\label{ssn: dictionary}

In summary the dictionary between type $3+1^m$ wild character varieties and supernova multiplicative quiver varieties is as follows (cf. 
\cite{rsode, slims} for the additive version).
Suppose $\wh \Ga$ is a simply-laced supernova graph with nodes $\wh I$, and core $\Ga$  with core nodes $I$, and we fix data 
$d\in \IZ^{\wh I}, q\in (\IC^*)^{\wh I}$. 
Then (if $\wh \Ga$ is coloured so as to have monochromatic core) 
we may consider the multiplicative quiver variety 
$\cM^{st}(\wh \Ga,q,d)$, and this may be ``read'' in terms of wild character varieties in $k+1$ different ways (if $\Ga$ is a complete $k$-partite graph, corresponding to a partition $I=I_1\sqcup\cdots\sqcup I_k$).
To describe these readings first recall that we can attach the following data to $\wh \Ga,q,d$:
Let $\IL_i$ be the leg of $\wh \Ga$ attached to the node $i\in I$.
This leg (and the data on it) determines a conjugacy class $\breve \cC_i\subset \GL(V_i)$ as in Lemma \ref{lem: orbits}
where $V_i=\IC^{d_i}$.
For any $i\in I$ let $U_i$ be the direct sum of the vector spaces at all the nodes of the other parts of $I$, not containing $i$:
If $i\in I_j\subset I$ then
$$U_i = \bigoplus_{l\in I\setminus I_j} V_l.$$

\begin{rmk}
Note that the injectivity/surjectivity conditions (as at the start of the proof of Prop. \ref{prop: stable isom}) show that $\cM^{st}(\wh \Ga,q,d)$
will be empty unless
$\dim(V_i) \le \dim(U_i)$ for all $i\in I$.
\end{rmk}

In turn the leg $\IL_i$ may be lengthened to construct a new leg $\wh \IL_i$ by adding a single node (connected by a single edge to the left end of $\IL_i$) with vector space $U_i$ and parameter $1\in \IC^*$.
The leg $\wh \IL_i$ determines a conjugacy class 
$\cC'_i\subset \GL(U_i)$ (as in Lemma \ref{lem: orbits})
for any $i\in I$, and we define $\cC_i\subset \GL(U_i)$ to be the inverse of the class $\cC'_i$ (if $g\in \cC_i$ then $g^{-1}\in \cC_i'$).
In general we will say that the class $\breve \cC_i\subset \GL(V_i)$ is the {\em child} of the class $\cC_i\subset \GL(U_i)$ (cf. Appendix 
\ref{apx: relating orbits} for the exact relation between the Jordan forms).

Then the basic data for the corresponding wild character varieties is as in the following table (there is a reading of the graph for any of the $k$ parts $I_j\subset I$, plus the generic reading, choosing no such part). The rank is the rank of the bundles on $\Si$, $m$ is the number of tame marked points in the finite plane, $H=H(Q)$ is the centralizer of the irregular type $Q$,
$\#A$ denotes the number of eigenvalues of the leading term $A$ of $Q$, and 
$\#T_l$  denotes the number of eigenvalues of the second term $T$ of $Q$
that are in the $l$th eigenspace of $A$. 
\renewcommand{\arraystretch}{2}
\begin{center}
  \begin{tabular}{ c | c | c | c | c | c | c | c}
    Reading & rank $r$           & $m$  & $H(Q)\subset \GL(\IC^r)$             & $\breve\bcC\subset H(Q)$ 
& $m$ classes 
& $\#A$ & $\#T_l$ 
\\ \hline\hline 
    generic & $\sum_I d_i$ & $0$  & $\Prod_I\GL(V_i)$  & $\breve\cC_i\bigl\vert_{i\in I}$  & $-$  & $k$ & $\#I_l$     \\ \hline 
 $I_j\subset I$  & $\sum_{I\setminus I_j} d_i$ & $\#I_j$  & 
$\Prod_{I\setminus I_j}\GL(V_i)$  & $\breve\cC_i\bigl\vert_{i\in I\setminus I_j}$  & $\cC_i\bigl\vert_{i\in I_j}$  & $k-1$ & $\#I_l$       \\ \hline
  \end{tabular}
\end{center}

\subsection{Examples}

For example suppose $U$ is a finite dimensional complex vector space and 
$\cC_1,\ldots,\cC_m,\breve \cC \subset \GL(U)$ are $m+1$ 
arbitrary conjugacy classes.
Consider the tame character variety parameterising irreducible representations
 of $\pi_1(\IC\setminus \text{$m$-points})$ in $\GL(U)$ with these fixed local conjugacy classes.

Let $d_i = \rank(M_i-1)$ for any $M_i\in \cC_i$, and let 
$V_i=\IC^{d_i}$, $V=\bigoplus_1^m V_i$.
Consider the groups 
$H=\Prod_1^m \GL(V_i)\subset G:=\GL(V)$
and define new conjugacy classes as follows:
Let  $\cC\subset G$ be the parent of the conjugacy class 
$\breve \cC \subset \GL(U)$
and let $\breve \cC_i\subset \GL(V_i)$ denote the child of $\cC_i$ (cf. Appendix \ref{apx: relating orbits}), and let 
$\breve \bcC = (\breve\cC_i)_1^m\subset H$. 
Recall the fission space $\gah=\cA(V)$ enables us to break the group $G$ to its subgroup $H$.

\begin{cor}
Then there is an isomorphism:
$$ (\cC_1\fusion{}\cdots \fusion{} \cC_m)^{st} \spqa{\breve \cC} \GL(U)
\ \cong\  (\cC \glue{G} \gah)^{st} \spqa{\breve \bcC} H$$

\noindent
from the tame character variety parameterising irreducible representations
 of $\pi_1(\IC\setminus \text{$m$-points})$ in $\GL(U)$, to a type $2+1$ 
wild character variety for the group $\GL(V)$.  
\end{cor}
\pf
These are the two non-generic readings of a star-shaped graph\footnote{The third reading of a star-shaped graph is as a type $3$ wild character
for the group $\GL(U\oplus V)$ and is of the form
$$\cB(Q)\spq_{\breve \cC\times\breve\bcC} \,\,\GL(U)\times H$$

\noindent
where $Q$ has coefficients $A,T$ such that $A$ has two eigenvalues (and has centralizer $\GL(U)\times\GL(V)$) and $T$ breaks the group further to $\GL(U)\times H$.}.
Since the core is bipartite ($k=2$), in each case $A$ only has one eigenvalue, which we can take to be zero.
Further, one of the parts only has one node, so in one of the readings
$T$ has only one eigenvalue, which again we may take to be zero, so in that reading $Q=0$, i.e. it is a tame case. 
\epf

Up to some relabellings this is Theorem \ref{thm: intro isom} from the introduction.

\begin{figure}[ht]
	\centering
	\input{star2.pstex_t}
	\caption{Star-shaped supernova graph, with core $\Ga(4,1)$}
\end{figure}
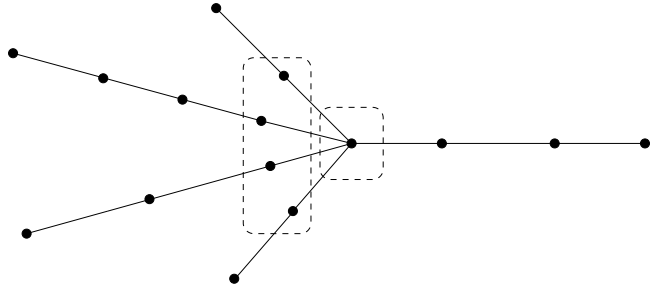

\begin{rmk}
The additive analogue of this example is easier, e.g. since the analogue of the fission operation $\cC\mapsto \cC\glue{}\gah$ is just restriction to $H$: 
A coadjoint orbit $\cO\subset \g^*$ is naturally a Hamiltonian $H$-space for any subgroup $H\subset G$, with moment map the dual of the derivative of the inclusion $\iota:H\hookrightarrow G$ 
(cf. the ``splaying'' operation in \cite{rsode}).
This motivates the operation $(\,\cdot\,)\glue{}\gah$ as the multiplicative analogue of this restriction operation (cf. \cite{fission}).
The additive analogue takes the form:
$$ (\cO_1\times\cdots \times \cO_m)^{st} \spqa{\breve \cO} \GL(U)
\ \cong\  (\cO)^{st} \spqa{\breve \bcO} H$$

\noindent
and essentially goes back to \cite{AHH-dual, Harn94}.
(See also \cite{iastalk07-nom, rsode, slims} for the quiver approach---this precise statement, with stability conditions and possibly non-semisimple orbits, is a special case of the results in \cite{slims} \S9.)
\end{rmk}

\section{Reflection isomorphisms and global Weyl group actions}
\label{sn: reflns}

In this section we will deduce the desired reflection isomorphisms
(to prove Theorem \ref{thm: refln isoms-intro} stated in the introduction).
This shows that many not-necessarily-affine Kac--Moody Weyl groups 
play a role in gauge theory (i.e. the theory of connections on bundles) on global Riemann surfaces.
This contrasts with the usual understanding of the special case of {\em affine} Kac--Moody Weyl groups (see e.g. \cite{pr-segal} \S5.1) 
in relation to loop groups (gauge theory on a circle).
A weaker (additive) version of these results appears in \cite{slims}.
Also, in the case of a star-shaped graph, similar (multiplicative) results appear in increasing generality in \cite{k2p} (for affine $D_4$), \cite{CB-Shaw, yamakawa-mpa}.\footnote{ 
The approach of \cite{k2p} (using the complex Fourier--Laplace transform of \cite{BJL81})  
is equivalent to the approach of \cite{CB-Shaw} (using the middle convolution of \cite{Katz-rls})---indeed 
it is easy to see that, generically, the scalar shift used in \cite{k2p}, which already appeared in \cite{BJL81} \S4.4, is middle convolution
(cf.  \cite{Katz-rls} \S2.10).
In fact it is straightforward to deduce explicit formulae for the complex version of Katz's middle convolution (equivalent to \cite{DettReit-katz}) from (prior) 
formulae in \cite{BJL81, malg-book} for the Fourier--Laplace transform on Betti data, a special case of which was used in \cite{k2p} (see also \cite{pecr}). %
}
For more general graphs the present approach will be less direct, although it is expected to be related to the Fourier--Laplace transform.

Recall (\S \ref{sn: mqvs and wcvs}) that given a marked 
irregular curve $\Si$ (of type $3+1^m$ for some $m\ge 0$), 
we have defined a supernova  graph $\wh \Ga$ and shown that the wild character variety of $\Si$ is isomorphic to a multiplicative quiver variety for the graph $\wh \Ga$.
Here we will recall how to attach a Kac--Moody   root system and Weyl group to a graph, and thus to a marked irregular curve.
This is the {\em global Weyl group} and root system of $\Si$.
In turn we will show how the global Weyl group acts to give isomorphisms  between multiplicative quiver varieties, corresponding to  isomorphisms 
between wild character varieties (often for different irregular curves).

\subsection{The Kac--Moody root system and Weyl group of a graph}\label{sn: km}

Let $\Gamma$ be a graph with no edge loops.
This determines a (symmetric) Kac--Moody algebra (see \cite{Kac-book}), and its
root system and Weyl group may be defined as follows.
Let $I$ be the set of nodes of 
$\Ga$ and let $n=\#I$ be the number of nodes.
Define the $n\times n$ (symmetric) Cartan matrix to be 
$$C=2\ \id - A$$
where $A$ is the adjacency matrix of $\Gamma$; 
the $i,j$ entry of $A$ is the number of edges connecting the nodes $i$ and $j$.
The {\em root lattice} $\IZ^I=\bigoplus_{i\in I} \IZ\eps_i$ inherits a bilinear form defined by 
\beq\label{eq: KM bil form}
(\eps_i,\eps_j) = C_{ij}.
\eeq

The simple reflections $s_i$, acting on the root lattice, are defined  by the 
formula
$$s_i(\be) := \be-(\be,\eps_i)\eps_i$$
for any $i\in I$.
They satisfy the relations
$$s_i^2=1,\qquad 
s_is_j=s_js_i \text{\ if $A_{ij}=0$,}\qquad
s_is_js_i=s_js_is_j \text{\ if $A_{ij}=1$}.
$$
By definition the Weyl group is the group generated by these simple reflections.
There are also dual reflections $r_i$ acting on the vector space 
$\IC^I$ by the formula 
$$r_i(\la) = \la - \la_i \al_i$$
where $\la = \sum_{i\in I} \la_i\eps_i\in \IC^I$ with $\la_i\in\IC$ and 
$\al_i := \sum_{j} (\eps_i,\eps_j)\eps_j\in \IC^I.$ 
By construction one has that $s_i(\be)\cdot r_i(\la) = \be \cdot \la$,
where the dot denotes the pairing given by 
$\eps_i\cdot\eps_j = \delta_{ij}$.
Exponentiating component-wise yields the multiplicative dual reflections (still denoted $r_i$) acting on $q\in (\IC^*)^I$, given in components by
$$r_i(q)_j =  q_i^{-(\eps_i,\eps_j)} q_j.$$
These satisfy $q^\be = r_i(q)^{s_i(\be)}$ where $q^\be := \Prod_I q_i^{\be_i}$.

The corresponding Kac--Moody root system is a subset of the {root lattice} 
$\IZ^I$.
It may be defined as the union of the set of real roots  and the set of imaginary roots, where

1) The simple roots are $\eps_i$ for $i\in I$,

2) The set of real roots is the Weyl group orbit of the set of simple roots,

3) Define the {\em fundamental region} to be the set of nonzero 
$\be\in\IN^I$ whose support is a connected subgraph of $\Gamma$ and such that $(\eps_i,\be)\le 0$ for all $i\in I$.
The set of imaginary roots is the union of the Weyl group orbit 
of the fundamental 
region and the orbit of minus the fundamental region.

This defines the root system (see \cite{Kac-book} Chapter 5 for the fact that this description does indeed give the roots of the corresponding Kac--Moody algebra).
By definition a root is positive if all its coefficients are $\ge 0$.
For example if $\Ga$ is an $ADE$ Dynkin diagram this gives the root system of the corresponding finite dimensional simple Lie algebra, or if $\Ga$ is an extended/affine $ADE$  Dynkin diagram then this is the root system of the corresponding affine Kac--Moody Lie algebra (closely related to the corresponding loop algebra), but of course there are many examples beyond these cases.

\subsection{Reflection isomorphisms}

Our aim is to establish the following.

\begin{thm}\label{thm: refln isoms}
Suppose $\wh \Ga$ is a simply-laced supernova graph with nodes $\wh I$, coloured so that its core $\Ga$ is monochromatic, $q,d$ are arbitrary, and 
$s_i, r_i$ are the corresponding simple reflections, generating the Kac--Moody Weyl group of $\wh\Ga$.
Suppose the support of $d$ intersects two distinct parts of the core $\Ga$ (so we ignore the trivial cases with just one part).
Then if $q_i\neq 1$ the multiplicative quiver varieties
\beq \label{eq: refln spaces}
\cM^{st}(\wh\Ga,q,d)\qquad\text{and}\qquad \cM^{st}(\wh\Ga,r_i(q),s_i(d))
\eeq
are isomorphic smooth symplectic algebraic varieties, for any node 
$i\in \wh I$.
\end{thm}
\pf
First suppose $i$ is not in the core (i.e. $i$ is in the interior of a leg). 
Then the desired reflection isomorphism arises 
simply by changing the choice of the marking, 
in the passage from wild character varieties 
in Theorem \ref{thm: 3wcv as mqv}: both of the multiplicative 
quiver varieties \eqref{eq: refln spaces}
are isomorphic to the same type $3$ wild character variety.
The details are now very similar to the additive case (\cite{slims} Corollary 9.12) so will be omitted here: in brief we swap the order of the two $\xi$'s appearing in the determination \eqref{eq: la defn} of $q_i$ from the marking of the corresponding leg.

Now suppose $i$ a node of the core, say $i\in I_j\subset I$.
Then, via Theorem \ref{thm: main reduced isom}, we may read  $\cM^{st}(\Ga,q,d)$ using the part $I_j\subset I$,
i.e. there is an isomorphic wild character variety $\MB^{st}(\Si,\bcC)$ for
a marked irregular curve $\Si$ of type $3+1^m$ with $m=\#I_j$, where the conjugacy classes 
$$\bcC\ =\ 
(\ \cC_i\bigl\vert_{i\in I_j}\ ,\ 
\breve\cC_i\bigl\vert_{i\in I\setminus I_j}\ )$$
for $\Si$ are determined by the legs of $\wh \Ga$,
as in \S\ref{ssn: dictionary}.
Then we can essentially proceed as above, changing the marking of the class 
$\cC_i$ corresponding to the chosen node $i$.
This is slightly complicated since the marking of $\cC_i$ must be special, but we can easily circumvent this, as follows. 
Note that once we have passed to the wild character variety we have the possibility to perform the following scalar shifts.

\begin{lem}
Suppose we choose $\ga_i\in \IC^*$ for all $i\in I_j$, and set 
$\ga = \Prod_{I_j} \ga_i$.
Define new conjugacy classes
$$\bcC'\ =\ 
(\ (\ga_i\cC_i)\bigl\vert_{i\in I_j}\ ,\ 
(\ga\breve\cC_i)\bigl\vert_{i\in I\setminus I_j}\ )$$
obtained by scaling the previous classes.
Then the wild character variety $\MB^{st}(\Si,\bcC')$
is isomorphic to $\MB^{st}(\Si,\bcC)$.  
\end{lem}
\pf
Straightforward (cf. \eqref{eq: explicit full space} and recall that $g$ is the moment map for $H'=\Prod_{I\setminus I_j}\GL(V_i)$, and we wish to  
impose $g_i\in\breve \cC_i$ for $i\in I\setminus I_j$, 
and $T_i\in\cC_i$).  
\epf

Now given a fixed node $i\in I_j\subset I$ we may change the marking of $\cC_i$, by swapping the order of the two $\xi$'s appearing in the determination \eqref{eq: la defn} of $q_i$, as before.
This new marking of  $\cC_i$ will in general no longer be special. But there is a unique scalar shift to make it special 
(with $\ga_l=1$ for all $l\in I_j\setminus\{i\}$).
Using this new marking, and this scalar shift, we can realise  
$\MB^{st}(\Si,\bcC')$ as a multiplicative quiver variety
$\cM^{st}(\Ga,q',d')$ for the same graph $\wh \Ga$ 
but with new data $q',d'$.
We claim that $q'=r_i(q), d'=s_i(d)$ and so the theorem follows.
The claim may be verified as follows (cf. \cite{slims} Corollary 9.18).
We will use subscripts ``$lk$'' to label the $k$th node down the  $l$th leg (for $l\in I, k=1,2,\ldots$).
Suppose the marking of $\breve \cC_l$ is 
$(\xi_{l1},\xi_{l2},\ldots)$ for each $l\in I$. 
The special marking of $\cC_i$ corresponding to the marking of 
$\breve \cC_i$ is thus
$$(1,\xi_{i1}^{-1},\xi_{i2}^{-1},\ldots).$$
Changing this marking by swapping the first two entries yields
$$(\xi_{i1}^{-1},1,\xi_{i2}^{-1},\ldots).$$
This is not special so we set $\ga=\ga_i=\xi_{i1}$ (and $\ga_l=1$ for all other $l\in I_j$) and perform the corresponding 
scalar shift to yield the special marking 
$$(1,\ga,\ga\xi_{i2}^{-1},\ldots)$$
of $\cC_i'=\ga\cC_i$.
The child  $\breve\cC_i'$ of $\cC_i'$ thus has marking
$$(\ga^{-1},\xi_{i2}/\ga,\ldots).$$
Also for $l\in I\setminus I_j$ the class $\breve\cC_l'=\ga\breve\cC_l$ has marking
$(\ga\xi_{l1},\ga\xi_{l2},\ldots)$. 
Now we can read off how the parameters $q$ have changed:
\begin{align*}
q_i=q_{i1}=\xi_{i1}\ &\mapsto\  \ga^{-1}=1/q_i,\\
q_{i2}= \xi_{i2}/\xi_{i1}\  &\mapsto\  (\xi_{i2}/\ga)/\ga^{-1}=\xi_{i2}=q_iq_{i2},\\
q_l=q_{l1}=\xi_{l1}\ &\mapsto\  \ga\xi_{l1} = q_iq_l,
\end{align*}
for all $l\in I\setminus I_j$, and all the 
other components of $q$ are unchanged. These are the components of $r_i(q)$.
For the dimension vector $d$, $d_{i1}$ is the only component which is changed, and the new value $d'_{i1}$ may be computed as follows.
Suppose $T_i\in \cC_i$. Then
$$d_{i1} = \rank(T_i-1) = \dim(U) -\dim\ker(T_i-1)$$
$$d_{i2} =  \dim(U) - \dim\ker(T_i-1) - \dim\ker(T_i-\xi_{i1}^{-1})$$
$$d'_{i1} = \rank(\ga T_i-1) =\rank( T_i-\ga^{-1}) = 
\dim(U) -\dim\ker(T_i-\xi_{i1}^{-1})$$

\noindent
so that $d'_{i1} = \dim(U) + d_{i2} - d_{i1}$, which is the
corresponding component of $s_i(d)$, given that 
$\dim(U) = \sum_{l\in I\setminus I_j} d_{l1}$.
\epf

\begin{rmk} (Affinized global Weyl groups, cf. \cite{rsode} Remark p.27.)
If we work with meromorphic connections (with parabolic/parahoric structures), rather than monodromy/Stokes data as we are here, then
in general a larger {\em affinized} Weyl group $W\sdp L$ appears, where
$L\cong\{\la\in \IZ^{\wh I}\st \la\cdot d=0\}$, 
since in essence one then has a choice of logarithm of each 
element $q_i$.
\end{rmk}

\section{Graphical Deligne--Simpson problems}\label{sn: gdsp}

We can now describe some basic linear algebra problems in terms of graphs, and present a conjectural solution.
Choose an ordered complete $k$-partite graph $\Ga$ with nodes $I$, for some integer $k$.
(As in \S\ref{sn: snova}
the graph $\Ga$ determines and is determined by a partition 
$\phi:I\onto J$ of a finite set $I$ 
into parts $I_j=\phi^{-1}(j)$ where $j\in J$.)
Let $V = \bigoplus_{i\in I} V_i$ be a finite dimensional $I$-graded vector space, and choose a conjugacy class $\breve \cC_i\subset \GL(V_i)$ for each $i\in I$.
Using the chosen ordering define the following subgroups of $\GL(V)$:
$$U_+ = \left\{ 1 + \sum_{i<j} u_{ij} \st u_{ij}:V_j\to V_i\right\},\qquad
U_- = \left\{ 1 + \sum_{i>j} u_{ij} \st u_{ij}:V_j\to V_i\right\}$$
and $H=\Prod_{i\in I}\GL(V_i)\subset \GL(V)$, the block-diagonal subgroup.
Thus we may consider the dense open subsets 
$$
U_-HU_+\subset \GL(V),\qquad
U_+HU_-\subset \GL(V)$$
of $\GL(V)$, which we will refer to as the big cell and the opposite big cell respectively (even though they are not cells).
Note that the choice of the classes $\breve\cC_i$ is equivalent to the choice of a conjugacy class  $\breve\cC = (\breve\cC_i) \subset H$.

Recall (from \S\ref{ssn: reps of graphs}) 
that a representation of $\Ga$ on $V$ is a choice of linear maps
$v_{ij}:V_j\to V_i$ for all $i,j\in I$ such that $\phi(i)\ne \phi(j)$.
Given such a representation we can use the ordering of $\Ga$ to define the following  unipotent elements of $\GL(V):$ 
$$v_+ = 1 + \sum_{i<j} v_{ij},\qquad v_-=1 + \sum_{i>j} v_{ij}.$$

The basic question then is: does there exist an irreducible representation 
$\{v_{ij}\}$ 
of $\Ga$  on $V$ such that $v_-v_+$ is in the opposite big cell of $\GL(V)$ with $H$ component  in the conjugacy class $\breve \cC$?
More explicitly: Is there an irreducible representation $\{v_{ij}\}$ of $\Ga$ on $V$,  
and elements $w_+\in U_+, w_-\in U_-$ 
and $g\in \breve \cC$ such that
$$v_-v_+ = w_+g w_-?$$
We will call this the graphical Deligne--Simpson problem for $\Ga$. 
(It is an explicit form of the irregular Deligne--Simpson problem $i$DS of 
\cite{gbs} \S9.4 for the case of type $3$ irregular curves---as explained there it encodes the question of when some of the \hk manifolds of \cite{wnabh} are nonempty.
The complete bipartite case was discussed in \cite{iastalk07-nom}.
Also the additive analogue was established in \cite{rsode, slims}; it is concerned with 
the spaces $\cM^*$ defined in \cite{smid}, and some of their generalisations.
Beware that Kostov \cite{kostov-nfDSp} has studied a {\em different, inequivalent} additive irregular Deligne--Simpson problem.)
Given what we have already done, the following propositions are now straightforward.

\begin{prop}
The solubility of the irregular Deligne--Simpson problem  for $\Ga$ is independent of the ordering of the graph $\Ga$.
\end{prop}

\begin{prop}
Suppose $\cC_1\ldots,\cC_m\subset \GL(V)$ are conjugacy classes. 
Consider the question:
is there a representation $\{v_{ij}\}$ of the graph $\Ga$ on $V$, and elements
$T_i\in\cC_i, w_\pm\in U_\pm, g\in\breve \cC$ such that
1) $T_1\cdots T_mv_-v_+ = w_+gw_-$
and 2) there are no nontrivial proper subrepresentations 
$U\subset V$ of $\Ga$ with 
$T_i(U)\subset U$ for $i=1,\ldots,m$.
Then this question is equivalent to a graphical Deligne--Simpson problem for a complete $(k+1)$-partite graph, obtained by adding a part of size $m$ to the nodes of $\Ga$.
\end{prop}

Taking $\Ga$ to be the graph with one node and no edges then yields:

\begin{prop}
The graphical Deligne--Simpson problem for star-shaped graphs is equivalent to the usual (tame) Deligne--Simpson problem.
\end{prop}
Thus we have a natural extension, involving a larger class of graphs.
In 2004 Crawley-Boevey \cite{CB-ihes}
made a conjecture, involving the Kac--Moody root system of a star shaped graph,
under which the tame Deligne--Simpson problem admits a solution.
(A proof of this conjecture was announced in 2006 \cite{CB.ICM}.)
See also for example
\cite{simpson-matrices, kostov-DSp99, CB-Shaw, Simpson-katz}.
This naturally leads us to conjecture the analogous result in the irregular case, now involving supernova graphs, as follows.  

Choose a marking of each conjugacy class $\breve \cC_i\subset \GL(V_i)$ for all $i\in I$ (see Definition \ref{defn: marking}). 
Then we get a supernova graph $\wh \Ga$ (with nodes $\wh I$) 
by gluing on legs, and parameters $q,d$, as in \S\ref{ssn: 3wcv as mqvs}.
The graph $\wh \Ga$ determines a Kac--Moody root system $\subset \IZ^{\wh I}$, as in \S\ref{sn: km}, and in particular the notion of positive roots.

\begin{conj}
The graphical Deligne--Simpson problem for the graph $\Ga$ and conjugacy classes 
$\breve \cC$ 
admits a solution if and only if
$d$ is a positive root, $q^d=1$, and, whenever it is possible to write 
$d=d_1+d_2+\cdots$  as a nontrivial sum of positive roots such that 
$q^{d_1}=q^{d_2}=\cdots =1$, then $\De(d)>\De(d_1)+\De(d_2)+\cdots$, where $\De(d) = 2-(d,d)$.
\end{conj}

\begin{rmk}
Similarly the analogous 
result should hold for all fission graphs $\Ga(Q)$ (see \S\ref{sn: snova}), 
not just the simply-laced ones appearing here, i.e. 
the nonemptiness of  $\Rep^*(\Ga(Q), V)^{st} \spq_{\breve\cC}\,\, H \cong 
\MB^{st}(\Si,\breve\cC)$ is controlled by the root system of the supernova graph $\wh \Ga(Q)$ determined by a marking of each class $\breve\cC_i$ in $\breve\cC$ (where 
$\Si=(\IP^1,\infty,Q)$). 
\end{rmk}

\section{Fission algebras}\label{sn: fissionalgs}

We will define some noncommutative 
algebras $\cF^q$ whose simple modules correspond to stable points 
of the multiplicative quiver varieties, generalising the multiplicative preprojective algebras of \cite{CB-Shaw}.
These algebras will be useful in the further study of 
the multiplicative quiver varieties, and are related to some generalisations of the double affine Hecke algebras.

Let $\Ga$ be a coloured quiver with nodes $I$  and colours $C$.
Choose some parameters $q\in (\IC^*)^I$ and a cyclic ordering of the colours for each node.
Construct a quiver $\wt \Ga$ as follows.
Take the double $\bar \Ga$ and for each colour $c\in C$ do the following:
1) for all $(i,j)\in I_c\times I_c$ with $i\neq j$ add an edge $w_{cij}$ 
from node $j$ to $i$,
2) for all $i\in I_c$ add a loop $\ga_{ci}$ from $i$ to $i$.
(Here we suppose each coloured graph $\Ga_c$ is connected.)
Let $\wt \Ga$ be the resulting quiver, and consider its path algebra 
$\cP(\wt \Ga)=\IC\wt\Ga$.\footnote{
It is the algebra having basis the paths in $\wt\Ga$---a {\em path} is a concatenation 
$a_n \cdots a_2a_1$ of oriented edges such that $h(a_{i})=t(a_{i+1})$.
(This includes for each $i\in I$ the path $e_i$ of length $0$ based at $i$.)
It is an associative algebra under path concatenation (setting to zero the product $a_2a_1$ if $h(a_{1})\neq t(a_2)$).
The identity is the element $1=\sum_{I} e_i$, and any subset $J\subset I$ determines an idempotent $\Id_J=\sum_{J} e_i$.}
For any $c\in C$ consider the elements
$$ 
w_{c+} = \Id_{I_c} + \sum_{i<j\in I_c} w_{cij},\qquad  
w_{c-} = \Id_{I_c} + \sum_{i>j\in I_c} w_{cij},\qquad
 \ga_c = \sum_{i\in I_c} \ga_{ci},$$
$$ 
v_{c+} = \Id_{I_c} + \sum_{i<j\in I_c} v_{cij},\qquad  
v_{c-} = \Id_{I_c} + \sum_{i>j\in I_c} v_{cij},$$
where $v_{cij}$ is the arrow in $\bar \Ga_c\subset \bar \Ga$ from 
$j$ to $i$ (or zero if $i,j$ are in the same part of $I_c$).
Then define $\cF^q=\cF^q(\Ga)$ by imposing the relations
$$v_{c-}v_{c+} = w_{c+}\ga_cw_{c-},\quad
\text{for each $c\in C$, and}\quad
\Prod_{\{c \,|\,i\in I_c\}} \ga_{ci} = q_i e_i
$$
using the ordering of $C$ at $i$, for each $i\in I$,  i.e. $\cF^q$ is obtained by quotienting $\cP(\wt \Ga)$ by the two-sided ideal generated by these relations.
This %
generalizes easily to the case when each subgraph 
$\Ga_c$ is an arbitrary fission graph, not just the simply-laced case.

Now suppose $V$ is a module for  $\cF^q$ such that $V$  is finite dimensional as a vector space over $\IC$.
Since it is a module for the path algebra of $\bar\Ga$, this determines a representation $\rho\in \Rep(\Ga,V)$ of the graph 
$\Ga$ on $V$.
By construction the representations which arise in this way 
are precisely those such that:

1) for any colour $c\in C$ the function (as in \eqref{eq: minor prod}) 
$$f_c=\Prod_{i\in I_c}\Delta_i:\Rep(\Ga_c,V)\to\IC$$
does not vanish on $\rho\bigl\vert_{\Ga_c}$ (the restriction of $\rho$ to the subgraph $\Ga_c$), so that $\rho$ is in the open subset $\Rep^*(\Ga,V)$ of $\Rep(\Ga,V)$, and

2)  $\mu(\rho) = q$ where $\mu:\Rep^*(\Ga,V)\to H$ is the moment map, as in \eqref{eq: mmap cmpts on IM}.

Thus, taking the S-equivalence class of $\rho\in \mu^{-1}(q)$, the $\cF^q$-module $V$ determines a point of the
 multiplicative quiver variety $\cM(\Ga,q,d)$,
where $d\in \IZ^I$ is the dimension vector of $V$.

\begin{prop}
Suppose $\Ga$ is a coloured quiver with nodes $I$, and  $q\in (\IC^*)^I, d\in \IZ_{\ge 0}^I$.
The multiplicative quiver variety $\cM(\Ga,q,d)$
has a stable point if and only if $\cF^q$ has a simple module with dimension vector $d$.
\end{prop}
\pf
By definition a module is simple if it has no proper nontrivial submodules. 
On the other hand an invertible graph representation is stable if and only if it has no proper nontrivial subrepresentations. As noted in 
Remark \ref{rmk: subreps}
any such subrepresentation is again invertible, and, considering the moment map condition, corresponds to a submodule for $\cF^q$. 
\epf

This is analogous to the deformed preprojective algebras \cite{CB-H}
(controlling the Nakajima quiver varieties \cite{nakaj-duke94}) and 
the
multiplicative preprojective algebras \cite{CB-Shaw}
(controlling the classical multiplicative quiver varieties).
Indeed if each monochromatic subgraph $\Ga_c$ just has one edge then the fission algebras are the same as the multiplicative preprojective algebras.

\begin{rmk}
Note that in Appendix A.1 of \cite{eor-gda}, Crawley-Boevey and Shaw explain 
a close relationship between the multiplicative preprojective
algebras and the ``generalized double affine Hecke algebras'' introduced by Etingof--Oblomkov--Rains. 
Since the  multiplicative preprojective algebras are the fission algebras of 
classically coloured quivers, this suggests there are yet more general double affine Hecke algebras, related to some of the other fission algebras.
\end{rmk}

\appendix

\section{Relating orbits} \label{apx: relating orbits}

Suppose $A:V\to W$ and $B:W\to V$ are linear maps between  two finite dimensional complex vector spaces $V,W$
such that $A$ is injective, $B$ is surjective, and the elements
$$T_V:= 1+ BA\in \End(V),\qquad T_W = 1+AB\in \End(W)$$
are invertible.
In this appendix we will recall the exact relation between the Jordan normal forms of $T_V$ and $T_W$.
Let $\cC\subset \GL(W)$ be the conjugacy class containing $T_W$ and let  
$\cC'\subset \GL(V)$ be the conjugacy class of $T_V$.

Recall that giving a Jordan form of an invertible matrix is equivalent to giving a partition (i.e. a Young diagram) $\varpi_s$
for each nonzero complex number $s$,
so that  $\varpi_s$ specifies the sizes of the Jordan blocks corresponding to the eigenvalue $s\in \IC^*$. 
For example the partition $\varpi_1=(2,2,1)$ specifies 
the $5\times 5$ unipotent matrix with three Jordan blocks of size $2,2$ and $1$ respectively
(and it corresponds to the Young diagram with three rows of lengths $2,2,1$).

\begin{prop}\label{prop: class relations}
Let $\{\varpi_s\st s\in \IC^*\}$ be the  partitions giving the Jordan form of $\cC$ and let $\{\varpi'_s\st s\in \IC^*\}$ be the  partitions giving the Jordan form of $\cC'$.
Then 
$\varpi'_s = \varpi_s \text{ if $s\neq$ 1, and}$
$\varpi'_1$ is obtained from 
$\varpi_1$ by deleting the first (i.e. longest) column of $\varpi_1$.
(In other words each part of $\varpi_1$ is decreased by one to obtain $\varpi'_1$.) 
\end{prop}

\pf
This follows immediately from the well-known relation between the adjoint orbit $\cO=\cC-1\subset \End(W)$ of $AB$ and the orbit 
$\cO'=\cC'-1\subset \End(V)$ of $BA$, which is reviewed for example in \cite{slims} Appendix D.
\epf

Let $\breve \cC\subset \GL(V)$ be the inverse of the class $\cC'$. 
In the body of the paper we will refer to $\breve \cC$ as the {\em child} of $\cC$ (and to $\cC$ as the {\em parent} of $\breve \cC$), if the vector spaces $V,W$ are clear from the context.

\renewcommand{\baselinestretch}{1}              %
\normalsize
\bibliographystyle{amsplain}    \label{biby}
\bibliography{../thesis/syr} 

\vspace{0.5cm}   

\noindent
D\'epartement de Math\'ematiques, \\
B\^{a}timent 425, \\
Facult\'e des Sciences d'Orsay \\
Universit\'e Paris-Sud \\
F-91405 Orsay Cedex \\

\noindent
boalch@gmail.com

\end{document}

%% file: macros.tex
\usepackage{verbatim}  
\usepackage{bbding}   

\def\reE@DeclareMathSymbol#1#2#3#4{%
    \let#1=\undefined
    \DeclareMathSymbol{#1}{#2}{#3}{#4}}
\DeclareSymbolFont{symbolsC}{U}{txsyc}{m}{n}
\SetSymbolFont{symbolsC}{bold}{U}{txsyc}{bx}{n}
\DeclareFontSubstitution{U}{txsyc}{m}{n}
\reE@DeclareMathSymbol{\strictiff}{\mathrel}{symbolsC}{76}

\newcommand\beq{\begin{equation}}
\newcommand\eeq{\end{equation}}
\newcommand\bal{\begin{align*}}
\newcommand\eal{\end{align*}}   
\newcommand\bmx{\left(\begin{matrix}}
\newcommand\emx{\end{matrix}\right)}
\newcommand\bsmx{\left(\begin{smallmatrix}}
\newcommand\esmx{\end{smallmatrix}\right)}

\newcommand{\onto}{\twoheadrightarrow}
\newcommand{\spq}{/\!\!/}

\providecommand{\spqa}[1]{\underset{#1}{/\!\!/}}
\providecommand{\sqpa}[1]{\underset{#1}{\setminus\!\!\setminus}}
\newcommand{\st}{\ \bigl\vert\ }

\def\part#1{\frac{\partial\phantom{q}}{\partial#1}}

\newcommand {\flb}{\lbrack\!\lbrack}
\newcommand {\frb}{\rbrack\!\rbrack}
\newcommand {\flp}{(\!(}
\newcommand {\frp}{)\!)}

\newcommand{\union}{\cup} 
 
\newcommand{\sdp}{{\ltimes}}
\newcommand{\glu}{\strictiff}
\newcommand{\glue}[1]{\underset{#1}{\strictiff}}
\newcommand{\fus}{\circledast}
\newcommand{\fusion}[1]{\underset{#1}{\circledast}}
\DeclareMathOperator*{\bigoasterisk}{\text{\rlap{$\bigoplus$}$\bigotimes$}}

\newcommand{\MB}{\mathcal{M}_{\text{\rm B}}}

\newcommand{\id}{\text{\rm Id}} 
\newcommand{\Id}{\text{\rm Id}}

\newcommand{\Lie}{{\mathop{\rm Lie}}}

\DeclareMathOperator{\ISto}{{\IS}to} 

\newcommand{\pap}[2]{{\ _{#1}\cA_{#2}}}
\newcommand{\papt}[2]{{\ _{#1}^{\phantom{2}}\cA_{#2}^2}}
\newcommand{\papk}[3]{  \ _{#1}^{\phantom{#3}}\cA_{#2}^{#3}     }
\newcommand{\gah}{\pap{G}{H}}  
\newcommand{\gahr}{\papk{G}{H}{r}}

\newcommand{\rank}{\mathop{\rm rank}}

\newcommand{\Prod}{\prod}

\newcommand{\Tr}{{\mathop{\rm Tr}}}
\DeclareMathOperator{\Hom}{Hom}         
\DeclareMathOperator{\Aut}{\mathop{\rm Aut}}

\newcommand{\GL}{{\mathop{\rm GL}}}

\newcommand{\U}{{\rm {U}}}	

\DeclareMathOperator{\Rep}{\rm Rep}
\renewcommand{\Im}{\mathop{\rm Im}}

\newcommand{\Ker}{\mathop{\rm Ker}}
\renewcommand{\ker}{\mathop{\rm Ker}}
\DeclareMathOperator{\End}{End}

\newcommand{\hk}{{hyperk\"ahler }}   

\newcommand{\ba}{{\bf a}}

\newcommand{\bH}{{\bf H}}

\newcommand{\bQ}{{\bf Q}}
\newcommand{\bs}{{\bf S}}
\newcommand{\bS}{{\bf S}}

\newcommand{\bU}{{\bf U}}

\newcommand{\IA}{\mathbb{A}}

\newcommand{\IC}{\mathbb{C}}
\newcommand{\ID}{\mathbb{D}}

\newcommand{\IH}{\mathbb{H}}

\newcommand{\IL}{\mathbb{L}}
\newcommand{\IM}{\mathbb{M}}
\newcommand{\IN}{\mathbb{N}}

\newcommand{\IP}{\mathbb{P}}                                     
\newcommand{\IQ}{\mathbb{Q}}                           
                           
\newcommand{\IS}{\mathbb{S}}

\newcommand{\IZ}{\mathbb{Z}}

\newcommand{\cA}{\mathcal{A}}
\newcommand{\cB}{\mathcal{B}}
\newcommand{\cC}{\mathcal{C}}

\newcommand{\cD}{\mathcal{D}}

\newcommand{\cF}{\mathcal{F}}

\newcommand{\cK}{\mathcal{K}}
\newcommand{\cL}{\mathcal{L}}

\newcommand{\cM}{\mathcal{M}}

\newcommand{\cN}{\mathcal{N}}
\newcommand{\cO}{\mathcal{O}}

\newcommand{\cP}{\mathcal{P}}

\newcommand{\cU}{\mathcal{U}}

\newcommand{\bcC}{\boldsymbol{\mathcal{C}}}
\newcommand{\bcO}{\boldsymbol{\mathcal{O}}}

\newcommand{\g}{       \mathfrak{g}     }

\newcommand{\lt}{\mathfrak{t}}
\newcommand{\lh}{\mathfrak{h}}

\newcommand{\gl}{       \mathfrak{gl}     } 

\newcommand{\wt}{\widetilde}

\newcommand{\wh}{\widehat}

\newcommand{\al}{\alpha}

\newcommand{\be}{\beta}
\newcommand{\ga}{\gamma}

\newcommand{\De}{\Delta}

\newcommand {\eps}{\varepsilon}

\newcommand{\Ga}{\Gamma}

\newcommand{\la}{\lambda}

\newcommand{\Si}{\Sigma}
\renewcommand{\th}{\theta}

\renewcommand{\bar}{\overline}

\makeatletter
 \newlength{\typesize}
\makeatother

\newlength{\vvoff}
\newlength{\hhoff}

\def\mapleft#1{\smash{
        \mathop{\longleftarrow}\limits^{#1}}}

\def\underset#1#2{\ \smash{\mathop{ #2 }\limits_{#1}}\ }

\newcommand{\pf}{\begin{bpf}}

\newcommand{\pfms}{\begin{bpfms}}
\newcommand{\epf}{\end{bpf}\hfill$\square$\\}           
\newcommand{\epfms}{\end{bpfms}\hfill$\square$\\}       

\newcommand{\idea}{\begin{bidea}}

\newcommand{\eidea}{\end{bidea}\hfill$\square$\\}           

\newcommand{\sk}{\begin{bsk}}    

\newcommand{\esk}{\end{bsk}\hfill$\square$\\}           
\newcommand{\sketch}{\begin{bsketch}}

\newcommand{\esketch}{\end{bsketch}\hfill$\square$\\}



%% file: macros-thm1.1.tex
\theoremstyle{plain}  \newtheorem{hypo}{Hypothesis}[section]
\newtheorem{thm}[hypo]{Theorem}
\newtheorem{prop}[hypo]{Proposition}

\newtheorem{conj}[hypo]{Conjecture}
\newtheorem{cor}[hypo]{Corollary}
\newtheorem{lem}[hypo]{Lemma}

\newtheorem {defn}[hypo]{Definition}

\theoremstyle{definition} \newtheorem{rmk}[hypo]{Remark}

\newtheorem{eg}[hypo]{Example}

%% file: quivers2.pstex_t
\begin{picture}(0,0)%
\includegraphics{quivers2.pstex}%
\end{picture}%
\setlength{\unitlength}{2901sp}%
\begingroup\makeatletter\ifx\SetFigFont\undefined%
\gdef\SetFigFont#1#2#3#4#5{%
  \reset@font\fontsize{#1}{#2pt}%
  \fontfamily{#3}\fontseries{#4}\fontshape{#5}%
  \selectfont}%
\fi\endgroup%
\begin{picture}(3706,3706)(2873,-4389)
\end{picture}%

%% file: trianglesbw2.pstex_t
\begin{picture}(0,0)%
\includegraphics{trianglesbw2.pstex}%
\end{picture}%
\setlength{\unitlength}{4558sp}%
\begingroup\makeatletter\ifx\SetFigFont\undefined%
\gdef\SetFigFont#1#2#3#4#5{%
  \reset@font\fontsize{#1}{#2pt}%
  \fontfamily{#3}\fontseries{#4}\fontshape{#5}%
  \selectfont}%
\fi\endgroup%
\begin{picture}(2806,885)(173,-10914)
\end{picture}%

%% file: 3pointblowups.pstex_t
\begin{picture}(0,0)%
\includegraphics{3pointblowups.pstex}%
\end{picture}%
\setlength{\unitlength}{2486sp}%
\begingroup\makeatletter\ifx\SetFigFont\undefined%
\gdef\SetFigFont#1#2#3#4#5{%
  \reset@font\fontsize{#1}{#2pt}%
  \fontfamily{#3}\fontseries{#4}\fontshape{#5}%
  \selectfont}%
\fi\endgroup%
\begin{picture}(5822,2299)(2689,-1898)
\end{picture}%

%% file: p654.v2.pstex_t
\begin{picture}(0,0)%
\includegraphics{p654.v2.pstex}%
\end{picture}%
\setlength{\unitlength}{4558sp}%
\begingroup\makeatletter\ifx\SetFigFont\undefined%
\gdef\SetFigFont#1#2#3#4#5{%
  \reset@font\fontsize{#1}{#2pt}%
  \fontfamily{#3}\fontseries{#4}\fontshape{#5}%
  \selectfont}%
\fi\endgroup%
\begin{picture}(5063,1218)(-1634,-11020)
\put(-1619,-10321){\makebox(0,0)[lb]{\smash{{\SetFigFont{12}{14.4}{\rmdefault}{\mddefault}{\updefault}{\color[rgb]{0,0,0}$n$}%
}}}}
\put(-719,-10321){\makebox(0,0)[lb]{\smash{{\SetFigFont{12}{14.4}{\rmdefault}{\mddefault}{\updefault}{\color[rgb]{0,0,0}$n$}%
}}}}
\put(-1034,-10906){\makebox(0,0)[lb]{\smash{{\SetFigFont{12}{14.4}{\rmdefault}{\mddefault}{\updefault}{\color[rgb]{0,0,0}$n$}%
}}}}
\put(2836,-10321){\makebox(0,0)[lb]{\smash{{\SetFigFont{12}{14.4}{\rmdefault}{\mddefault}{\updefault}{\color[rgb]{0,0,0}$n$}%
}}}}
\put(-1079,-10321){\makebox(0,0)[lb]{\smash{{\SetFigFont{12}{14.4}{\rmdefault}{\mddefault}{\updefault}{\color[rgb]{0,0,0}$2n$}%
}}}}
\put(-269,-10321){\makebox(0,0)[lb]{\smash{{\SetFigFont{12}{14.4}{\rmdefault}{\mddefault}{\updefault}{\color[rgb]{0,0,0}$1$}%
}}}}
\put(1531,-10321){\makebox(0,0)[lb]{\smash{{\SetFigFont{12}{14.4}{\rmdefault}{\mddefault}{\updefault}{\color[rgb]{0,0,0}$1$}%
}}}}
\put(766,-9961){\makebox(0,0)[lb]{\smash{{\SetFigFont{12}{14.4}{\rmdefault}{\mddefault}{\updefault}{\color[rgb]{0,0,0}$n$}%
}}}}
\put(-1034,-10006){\makebox(0,0)[lb]{\smash{{\SetFigFont{12}{14.4}{\rmdefault}{\mddefault}{\updefault}{\color[rgb]{0,0,0}$n$}%
}}}}
\put(1126,-10321){\makebox(0,0)[lb]{\smash{{\SetFigFont{12}{14.4}{\rmdefault}{\mddefault}{\updefault}{\color[rgb]{0,0,0}$n$}%
}}}}
\put(766,-10951){\makebox(0,0)[lb]{\smash{{\SetFigFont{12}{14.4}{\rmdefault}{\mddefault}{\updefault}{\color[rgb]{0,0,0}$n$}%
}}}}
\put(3331,-10321){\makebox(0,0)[lb]{\smash{{\SetFigFont{12}{14.4}{\rmdefault}{\mddefault}{\updefault}{\color[rgb]{0,0,0}$1$}%
}}}}
\put( 91,-10321){\makebox(0,0)[lb]{\smash{{\SetFigFont{12}{14.4}{\rmdefault}{\mddefault}{\updefault}{\color[rgb]{0,0,0}$n$}%
}}}}
\put(1891,-10006){\makebox(0,0)[lb]{\smash{{\SetFigFont{12}{14.4}{\rmdefault}{\mddefault}{\updefault}{\color[rgb]{0,0,0}$n$}%
}}}}
\put(1891,-10906){\makebox(0,0)[lb]{\smash{{\SetFigFont{12}{14.4}{\rmdefault}{\mddefault}{\updefault}{\color[rgb]{0,0,0}$n$}%
}}}}
\end{picture}%

%% file: splaying.isom3.pstex_t
\begin{picture}(0,0)%
\includegraphics{splaying.isom3.pstex}%
\end{picture}%
\setlength{\unitlength}{3315sp}%
\begingroup\makeatletter\ifx\SetFigFontNFSS\undefined%
\gdef\SetFigFontNFSS#1#2#3#4#5{%
  \reset@font\fontsize{#1}{#2pt}%
  \fontfamily{#3}\fontseries{#4}\fontshape{#5}%
  \selectfont}%
\fi\endgroup%
\begin{picture}(7493,1803)(-7214,-10165)
\put(-4499,-9781){\makebox(0,0)[lb]{\smash{{\SetFigFontNFSS{10}{12.0}{\rmdefault}{\mddefault}{\updefault}{$W_s$}%
}}}}
\put(-4499,-9151){\makebox(0,0)[lb]{\smash{{\SetFigFontNFSS{10}{12.0}{\rmdefault}{\mddefault}{\updefault}{$W_2$}%
}}}}
\put(-4499,-8881){\makebox(0,0)[lb]{\smash{{\SetFigFontNFSS{10}{12.0}{\rmdefault}{\mddefault}{\updefault}{$W_1$}%
}}}}
\put(-3194,-9151){\makebox(0,0)[lb]{\smash{{\SetFigFontNFSS{10}{12.0}{\rmdefault}{\mddefault}{\updefault}{$V$}%
}}}}
\put(-7199,-10096){\makebox(0,0)[lb]{\smash{{\SetFigFontNFSS{10}{12.0}{\rmdefault}{\mddefault}{\updefault}{$W_s$}%
}}}}
\put(-7199,-9016){\makebox(0,0)[lb]{\smash{{\SetFigFontNFSS{10}{12.0}{\rmdefault}{\mddefault}{\updefault}{$W_2$}%
}}}}
\put(-7199,-8521){\makebox(0,0)[lb]{\smash{{\SetFigFontNFSS{10}{12.0}{\rmdefault}{\mddefault}{\updefault}{$W_1$}%
}}}}
\put(-5759,-9016){\makebox(0,0)[lb]{\smash{{\SetFigFontNFSS{10}{12.0}{\rmdefault}{\mddefault}{\updefault}{$V$}%
}}}}
\put(-5759,-9241){\makebox(0,0)[lb]{\smash{{\SetFigFontNFSS{10}{12.0}{\rmdefault}{\mddefault}{\updefault}{$V$}%
}}}}
\put(-5759,-9646){\makebox(0,0)[lb]{\smash{{\SetFigFontNFSS{10}{12.0}{\rmdefault}{\mddefault}{\updefault}{$V$}%
}}}}
\put(-1709,-9331){\makebox(0,0)[lb]{\smash{{\SetFigFontNFSS{10}{12.0}{\rmdefault}{\mddefault}{\updefault}{$\cA(W)\glue{}$}%
}}}}
\put(-2564,-9331){\makebox(0,0)[lb]{\smash{{\SetFigFontNFSS{10}{12.0}{\rmdefault}{\mddefault}{\updefault}{$\cong$}%
}}}}
\put(-764,-9151){\makebox(0,0)[lb]{\smash{{\SetFigFontNFSS{10}{12.0}{\rmdefault}{\mddefault}{\updefault}{$W$}%
}}}}
\put(136,-9151){\makebox(0,0)[lb]{\smash{{\SetFigFontNFSS{10}{12.0}{\rmdefault}{\mddefault}{\updefault}{$V$}%
}}}}
\put(-5129,-9331){\makebox(0,0)[lb]{\smash{{\SetFigFontNFSS{10}{12.0}{\rmdefault}{\mddefault}{\updefault}{$=$}%
}}}}
\end{picture}%

%% file: second.isom.v3.pstex_t
\begin{picture}(0,0)%
\includegraphics{second.isom.v3.pstex}%
\end{picture}%
\setlength{\unitlength}{3315sp}%
\begingroup\makeatletter\ifx\SetFigFont\undefined%
\gdef\SetFigFont#1#2#3#4#5{%
  \reset@font\fontsize{#1}{#2pt}%
  \fontfamily{#3}\fontseries{#4}\fontshape{#5}%
  \selectfont}%
\fi\endgroup%
\begin{picture}(7860,5075)(-464,-7428)
\put(3151,-3706){\makebox(0,0)[lb]{\smash{{\SetFigFont{10}{12.0}{\rmdefault}{\mddefault}{\updefault}{\color[rgb]{0,0,0}$W_s$}%
}}}}
\put(3151,-3076){\makebox(0,0)[lb]{\smash{{\SetFigFont{10}{12.0}{\rmdefault}{\mddefault}{\updefault}{\color[rgb]{0,0,0}$W_2$}%
}}}}
\put(3151,-2806){\makebox(0,0)[lb]{\smash{{\SetFigFont{10}{12.0}{\rmdefault}{\mddefault}{\updefault}{\color[rgb]{0,0,0}$W_1$}%
}}}}
\put(4456,-3076){\makebox(0,0)[lb]{\smash{{\SetFigFont{10}{12.0}{\rmdefault}{\mddefault}{\updefault}{\color[rgb]{0,0,0}$V$}%
}}}}
\put(4681,-3256){\makebox(0,0)[lb]{\smash{{\SetFigFont{10}{12.0}{\rmdefault}{\mddefault}{\updefault}{\color[rgb]{0,0,0}$\glue{}\cA(V)$}%
}}}}
\put(5986,-3256){\makebox(0,0)[lb]{\smash{{\SetFigFont{10}{12.0}{\rmdefault}{\mddefault}{\updefault}{\color[rgb]{0,0,0}$\cong$}%
}}}}
\put(-449,-4111){\makebox(0,0)[lb]{\smash{{\SetFigFont{10}{12.0}{\rmdefault}{\mddefault}{\updefault}{\color[rgb]{0,0,0}$W_s$}%
}}}}
\put(-449,-3031){\makebox(0,0)[lb]{\smash{{\SetFigFont{10}{12.0}{\rmdefault}{\mddefault}{\updefault}{\color[rgb]{0,0,0}$W_2$}%
}}}}
\put(-449,-2536){\makebox(0,0)[lb]{\smash{{\SetFigFont{10}{12.0}{\rmdefault}{\mddefault}{\updefault}{\color[rgb]{0,0,0}$W_1$}%
}}}}
\put(766,-2761){\makebox(0,0)[lb]{\smash{{\SetFigFont{10}{12.0}{\rmdefault}{\mddefault}{\updefault}{\color[rgb]{0,0,0}$V$}%
}}}}
\put(2701,-4471){\makebox(0,0)[rb]{\smash{{\SetFigFont{10}{12.0}{\rmdefault}{\mddefault}{\updefault}{\color[rgb]{0,0,0}$V_1$}%
}}}}
\put(2701,-4741){\makebox(0,0)[rb]{\smash{{\SetFigFont{10}{12.0}{\rmdefault}{\mddefault}{\updefault}{\color[rgb]{0,0,0}$V_2$}%
}}}}
\put(2701,-5371){\makebox(0,0)[rb]{\smash{{\SetFigFont{10}{12.0}{\rmdefault}{\mddefault}{\updefault}{\color[rgb]{0,0,0}$V_r$}%
}}}}
\put(901,-5371){\makebox(0,0)[lb]{\smash{{\SetFigFont{10}{12.0}{\rmdefault}{\mddefault}{\updefault}{\color[rgb]{0,0,0}$W_s$}%
}}}}
\put(901,-4741){\makebox(0,0)[lb]{\smash{{\SetFigFont{10}{12.0}{\rmdefault}{\mddefault}{\updefault}{\color[rgb]{0,0,0}$W_2$}%
}}}}
\put(901,-4471){\makebox(0,0)[lb]{\smash{{\SetFigFont{10}{12.0}{\rmdefault}{\mddefault}{\updefault}{\color[rgb]{0,0,0}$W_1$}%
}}}}
\put(3376,-6046){\makebox(0,0)[rb]{\smash{{\SetFigFont{10}{12.0}{\rmdefault}{\mddefault}{\updefault}{\color[rgb]{0,0,0}$V_1$}%
}}}}
\put(3376,-6316){\makebox(0,0)[rb]{\smash{{\SetFigFont{10}{12.0}{\rmdefault}{\mddefault}{\updefault}{\color[rgb]{0,0,0}$V_2$}%
}}}}
\put(3376,-6946){\makebox(0,0)[rb]{\smash{{\SetFigFont{10}{12.0}{\rmdefault}{\mddefault}{\updefault}{\color[rgb]{0,0,0}$V_r$}%
}}}}
\put(136,-6496){\makebox(0,0)[lb]{\smash{{\SetFigFont{10}{12.0}{\rmdefault}{\mddefault}{\updefault}{\color[rgb]{0,0,0}$\cong$}%
}}}}
\put(2071,-6316){\makebox(0,0)[rb]{\smash{{\SetFigFont{10}{12.0}{\rmdefault}{\mddefault}{\updefault}{\color[rgb]{0,0,0}$W$}%
}}}}
\put(946,-6496){\makebox(0,0)[lb]{\smash{{\SetFigFont{10}{12.0}{\rmdefault}{\mddefault}{\updefault}{\color[rgb]{0,0,0}$\cA(W)\glue{}$}%
}}}}
\put(6166,-6001){\rotatebox{360.0}{\makebox(0,0)[rb]{\smash{{\SetFigFont{10}{12.0}{\rmdefault}{\mddefault}{\updefault}{\color[rgb]{0,0,0}$W$}%
}}}}}
\put(7381,-5776){\rotatebox{360.0}{\makebox(0,0)[rb]{\smash{{\SetFigFont{10}{12.0}{\rmdefault}{\mddefault}{\updefault}{\color[rgb]{0,0,0}$V_1$}%
}}}}}
\put(7381,-6271){\rotatebox{360.0}{\makebox(0,0)[rb]{\smash{{\SetFigFont{10}{12.0}{\rmdefault}{\mddefault}{\updefault}{\color[rgb]{0,0,0}$V_2$}%
}}}}}
\put(7381,-7351){\rotatebox{360.0}{\makebox(0,0)[rb]{\smash{{\SetFigFont{10}{12.0}{\rmdefault}{\mddefault}{\updefault}{\color[rgb]{0,0,0}$V_r$}%
}}}}}
\put(1081,-3256){\makebox(0,0)[lb]{\smash{{\SetFigFont{10}{12.0}{\rmdefault}{\mddefault}{\updefault}{\color[rgb]{0,0,0}$\glue{}\cA(V)$}%
}}}}
\put(2386,-3256){\makebox(0,0)[lb]{\smash{{\SetFigFont{10}{12.0}{\rmdefault}{\mddefault}{\updefault}{\color[rgb]{0,0,0}$=$}%
}}}}
\put(5716,-4921){\makebox(0,0)[lb]{\smash{{\SetFigFont{10}{12.0}{\rmdefault}{\mddefault}{\updefault}{\color[rgb]{0,0,0}$\glue{}\cA(V)$}%
}}}}
\put(3601,-4921){\makebox(0,0)[lb]{\smash{{\SetFigFont{10}{12.0}{\rmdefault}{\mddefault}{\updefault}{\color[rgb]{0,0,0}$\cA(W)\glue{}$}%
}}}}
\put(4546,-4741){\makebox(0,0)[lb]{\smash{{\SetFigFont{10}{12.0}{\rmdefault}{\mddefault}{\updefault}{\color[rgb]{0,0,0}$W$}%
}}}}
\put(5446,-4741){\makebox(0,0)[lb]{\smash{{\SetFigFont{10}{12.0}{\rmdefault}{\mddefault}{\updefault}{\color[rgb]{0,0,0}$V$}%
}}}}
\put(3016,-4921){\makebox(0,0)[lb]{\smash{{\SetFigFont{10}{12.0}{\rmdefault}{\mddefault}{\updefault}{\color[rgb]{0,0,0}$\cong$}%
}}}}
\put(4906,-6496){\makebox(0,0)[lb]{\smash{{\SetFigFont{10}{12.0}{\rmdefault}{\mddefault}{\updefault}{\color[rgb]{0,0,0}$\cA(W)\glue{}$}%
}}}}
\put(4051,-6496){\makebox(0,0)[lb]{\smash{{\SetFigFont{10}{12.0}{\rmdefault}{\mddefault}{\updefault}{\color[rgb]{0,0,0}$=$}%
}}}}
\end{picture}%

%% file: hfission4.pstex_t
\begin{picture}(0,0)%
\includegraphics{hfission4.pstex}%
\end{picture}%
\setlength{\unitlength}{3729sp}%
\begingroup\makeatletter\ifx\SetFigFont\undefined%
\gdef\SetFigFont#1#2#3#4#5{%
  \reset@font\fontsize{#1}{#2pt}%
  \fontfamily{#3}\fontseries{#4}\fontshape{#5}%
  \selectfont}%
\fi\endgroup%
\begin{picture}(7095,1566)(-6269,-10084)
\put(181,-9556){\rotatebox{360.0}{\makebox(0,0)[rb]{\smash{{\SetFigFont{11}{13.2}{\rmdefault}{\mddefault}{\updefault}{\color[rgb]{0,0,0}$W_1$}%
}}}}}
\put(811,-10006){\rotatebox{360.0}{\makebox(0,0)[rb]{\smash{{\SetFigFont{11}{13.2}{\rmdefault}{\mddefault}{\updefault}{\color[rgb]{0,0,0}$W_3$}%
}}}}}
\put(811,-8701){\rotatebox{360.0}{\makebox(0,0)[rb]{\smash{{\SetFigFont{11}{13.2}{\rmdefault}{\mddefault}{\updefault}{\color[rgb]{0,0,0}$W_2$}%
}}}}}
\put(-44,-9331){\rotatebox{360.0}{\makebox(0,0)[rb]{\smash{{\SetFigFont{11}{13.2}{\rmdefault}{\mddefault}{\updefault}{\color[rgb]{0,0,0}$\cong\qquad\cA(W_1)\glue{}$}%
}}}}}
\put(-3194,-8701){\makebox(0,0)[lb]{\smash{{\SetFigFont{11}{13.2}{\rmdefault}{\mddefault}{\updefault}{\color[rgb]{0,0,0}$V_1$}%
}}}}
\put(-3194,-10006){\makebox(0,0)[lb]{\smash{{\SetFigFont{11}{13.2}{\rmdefault}{\mddefault}{\updefault}{\color[rgb]{0,0,0}$V_2$}%
}}}}
\put(-2294,-8701){\makebox(0,0)[lb]{\smash{{\SetFigFont{11}{13.2}{\rmdefault}{\mddefault}{\updefault}{\color[rgb]{0,0,0}$W_2$}%
}}}}
\put(-2339,-10006){\makebox(0,0)[lb]{\smash{{\SetFigFont{11}{13.2}{\rmdefault}{\mddefault}{\updefault}{\color[rgb]{0,0,0}$W_3$}%
}}}}
\put(-5354,-9556){\makebox(0,0)[lb]{\smash{{\SetFigFont{11}{13.2}{\rmdefault}{\mddefault}{\updefault}{\color[rgb]{0,0,0}$U_1$}%
}}}}
\put(-6254,-10006){\makebox(0,0)[lb]{\smash{{\SetFigFont{11}{13.2}{\rmdefault}{\mddefault}{\updefault}{\color[rgb]{0,0,0}$V_2$}%
}}}}
\put(-6254,-8701){\makebox(0,0)[lb]{\smash{{\SetFigFont{11}{13.2}{\rmdefault}{\mddefault}{\updefault}{\color[rgb]{0,0,0}$V_1$}%
}}}}
\put(-5174,-9331){\makebox(0,0)[lb]{\smash{{\SetFigFont{11}{13.2}{\rmdefault}{\mddefault}{\updefault}{\color[rgb]{0,0,0}$\glue{}\cA^2(U_1)\qquad\cong$}%
}}}}
\end{picture}%

%% file: halo.on.curve2.pstex_t
\begin{picture}(0,0)%
\includegraphics{halo.on.curve2.pstex}%
\end{picture}%
\setlength{\unitlength}{2693sp}%
\begingroup\makeatletter\ifx\SetFigFont\undefined%
\gdef\SetFigFont#1#2#3#4#5{%
  \reset@font\fontsize{#1}{#2pt}%
  \fontfamily{#3}\fontseries{#4}\fontshape{#5}%
  \selectfont}%
\fi\endgroup%
\begin{picture}(6744,3570)(7176,-7684)
\put(8686,-5776){\makebox(0,0)[lb]{\smash{{\SetFigFont{8}{9.6}{\rmdefault}{\mddefault}{\updefault}{\color[rgb]{0,0,0}$\partial_1$}%
}}}}
\put(9901,-5146){\makebox(0,0)[lb]{\smash{{\SetFigFont{8}{9.6}{\rmdefault}{\mddefault}{\updefault}{\color[rgb]{0,0,0}$d$}%
}}}}
\put(9271,-4831){\makebox(0,0)[lb]{\smash{{\SetFigFont{8}{9.6}{\rmdefault}{\mddefault}{\updefault}{\color[rgb]{0,0,0}$e(d)$}%
}}}}
\put(9271,-6496){\makebox(0,0)[lb]{\smash{{\SetFigFont{8}{9.6}{\rmdefault}{\mddefault}{\updefault}{\color[rgb]{0,0,0}$\IH_1$}%
}}}}
\end{picture}%

%% file: leg4.pstex_t
\begin{picture}(0,0)%
\includegraphics{leg4.pstex}%
\end{picture}%
\setlength{\unitlength}{3729sp}%
\begingroup\makeatletter\ifx\SetFigFont\undefined%
\gdef\SetFigFont#1#2#3#4#5{%
  \reset@font\fontsize{#1}{#2pt}%
  \fontfamily{#3}\fontseries{#4}\fontshape{#5}%
  \selectfont}%
\fi\endgroup%
\begin{picture}(3722,846)(1282,-4)
\put(1297,512){\makebox(0,0)[lb]{\smash{{\SetFigFont{11}{13.2}{\rmdefault}{\mddefault}{\updefault}{\color[rgb]{0,0,0}$n$}%
}}}}
\put(1303,132){\makebox(0,0)[lb]{\smash{{\SetFigFont{11}{13.2}{\rmdefault}{\mddefault}{\updefault}{\color[rgb]{0,0,0}$1$}%
}}}}
\put(2203,132){\makebox(0,0)[lb]{\smash{{\SetFigFont{11}{13.2}{\rmdefault}{\mddefault}{\updefault}{\color[rgb]{0,0,0}$2$}%
}}}}
\put(3103,132){\makebox(0,0)[lb]{\smash{{\SetFigFont{11}{13.2}{\rmdefault}{\mddefault}{\updefault}{\color[rgb]{0,0,0}$3$}%
}}}}
\put(2197,512){\makebox(0,0)[lb]{\smash{{\SetFigFont{11}{13.2}{\rmdefault}{\mddefault}{\updefault}{\color[rgb]{0,0,0}$d_2$}%
}}}}
\put(3097,512){\makebox(0,0)[lb]{\smash{{\SetFigFont{11}{13.2}{\rmdefault}{\mddefault}{\updefault}{\color[rgb]{0,0,0}$d_3$}%
}}}}
\put(4897,512){\makebox(0,0)[lb]{\smash{{\SetFigFont{11}{13.2}{\rmdefault}{\mddefault}{\updefault}{\color[rgb]{0,0,0}$d_{w}$}%
}}}}
\put(4903,132){\makebox(0,0)[lb]{\smash{{\SetFigFont{11}{13.2}{\rmdefault}{\mddefault}{\updefault}{\color[rgb]{0,0,0}$w$}%
}}}}
\put(1711,659){\makebox(0,0)[lb]{\smash{{\SetFigFont{11}{13.2}{\rmdefault}{\mddefault}{\updefault}{\color[rgb]{0,0,0}$b_1$}%
}}}}
\put(2611,659){\makebox(0,0)[lb]{\smash{{\SetFigFont{11}{13.2}{\rmdefault}{\mddefault}{\updefault}{\color[rgb]{0,0,0}$b_2$}%
}}}}
\put(4411,659){\makebox(0,0)[lb]{\smash{{\SetFigFont{11}{13.2}{\rmdefault}{\mddefault}{\updefault}{\color[rgb]{0,0,0}$b_l$}%
}}}}
\put(1711, 74){\makebox(0,0)[lb]{\smash{{\SetFigFont{11}{13.2}{\rmdefault}{\mddefault}{\updefault}{\color[rgb]{0,0,0}$a_1$}%
}}}}
\put(2611, 74){\makebox(0,0)[lb]{\smash{{\SetFigFont{11}{13.2}{\rmdefault}{\mddefault}{\updefault}{\color[rgb]{0,0,0}$a_2$}%
}}}}
\put(4411, 74){\makebox(0,0)[lb]{\smash{{\SetFigFont{11}{13.2}{\rmdefault}{\mddefault}{\updefault}{\color[rgb]{0,0,0}$a_l$}%
}}}}
\end{picture}%

%% file: star2.pstex_t
\begin{picture}(0,0)%
\includegraphics{star2.pstex}%
\end{picture}%
\setlength{\unitlength}{2486sp}%
\begingroup\makeatletter\ifx\SetFigFont\undefined%
\gdef\SetFigFont#1#2#3#4#5{%
  \reset@font\fontsize{#1}{#2pt}%
  \fontfamily{#3}\fontseries{#4}\fontshape{#5}%
  \selectfont}%
\fi\endgroup%
\begin{picture}(6406,2806)(848,-3714)
\end{picture}%

%% file: cmqv.bbl
\def\cprime{$'$} \def\cprime{$'$} \def\cprime{$'$} \def\cprime{$'$}
\providecommand{\bysame}{\leavevmode\hbox to3em{\hrulefill}\thinspace}
\providecommand{\MR}{\relax\ifhmode\unskip\space\fi MR }
\providecommand{\MRhref}[2]{%
  \href{http://www.ams.org/mathscinet-getitem?mr=#1}{#2}
}
\providecommand{\href}[2]{#2}
\begin{thebibliography}{10}

\bibitem{AHH-dual}
M.~R. Adams, J.~Harnad, and J.~Hurtubise, \emph{Dual moment maps into loop
  algebras}, Lett. Math. Phys. \textbf{20} (1990), no.~4, 299--308.

\bibitem{AMM}
A.~Alekseev, A.~Malkin, and E.~Meinrenken, \emph{Lie group valued moment maps},
  J. Differential Geom. \textbf{48} (1998), no.~3, 445--495, math.DG/9707021.

\bibitem{AL-p6ims}
D.~Arinkin and S.~Lysenko, \emph{Isomorphisms between moduli spaces of {${\rm
  SL}(2)$}-bundles with connections on {${\bf P}^1\setminus \{x\sb 1,\cdots,
  x\sb 4\}$}}, Math. Res. Lett. \textbf{4} (1997), no.~2-3, 181--190.

\bibitem{adhm}
M.~F. Atiyah, V.~Drinfel{\cprime}d, N.~J. Hitchin, and Yu.~I. Manin,
  \emph{Construction of instantons}, Phys. Lett. A \textbf{65} (1978), no.~3,
  185--187.

\bibitem{AB83}
M.F. Atiyah and R.~Bott, \emph{The {Y}ang-{M}ills equations over {R}iemann
  surfaces}, Phil. Trans. R. Soc. London \textbf{308} (1983), 523--615.

\bibitem{BJL81}
W.~Balser, W.B. Jurkat, and D.A. Lutz, \emph{On the reduction of connection
  problems for differential equations with an irregular singularity to ones
  with only regular singularities, {I}.}, SIAM J. Math. Anal. \textbf{12}
  (1981), no.~5, 691--721.

\bibitem{wnabh}
O.~Biquard and P.~P. Boalch, \emph{Wild non-abelian {H}odge theory on curves},
  Compositio Math. \textbf{140} (2004), no.~1, 179--204.

\bibitem{blanc-email}
J.~Blanc, Email correspondence 19/3/2008.

\bibitem{thesis}
P.~P. Boalch, \emph{{{S}ymplectic geometry and isomonodromic deformations}},
  D.Phil. Thesis, Oxford University, 1999.

\bibitem{smid}
\bysame, \emph{{S}ymplectic manifolds and isomonodromic deformations}, Adv. in
  Math. \textbf{163} (2001), 137--205.

\bibitem{saqh02}
\bysame, \emph{Quasi-{H}amiltonian geometry of meromorphic connections},
  (2002), arXiv:0203161.

\bibitem{pecr}
\bysame, \emph{Painlev\'e equations and complex reflections}, Ann. Inst.
  Fourier \textbf{53} (2003), no.~4, 1009--1022, Proceedings of conference in
  honour of Fr\'ed\'eric Pham, Nice, July 2002.

\bibitem{k2p}
\bysame, \emph{From {K}lein to {P}ainlev\'e via {F}ourier, {L}aplace and
  {J}imbo}, Proc. London Math. Soc. \textbf{90} (2005), no.~3, 167--208,
  math.AG/0308221.

\bibitem{saqh}
\bysame, \emph{Quasi-{H}amiltonian geometry of meromorphic connections}, Duke
  Math. J. \textbf{139} (2007), no.~2, 369--405, (Beware section 6 of the
  published version is not in the 2002 arXiv version).

\bibitem{iastalk07-nom}
\bysame, \emph{Some geometry of irregular connections on curves}, 2007,
  November 27, Talk at Workshop on Gauge Theory and Representation Theory, IAS
  Princeton, transparencies available on author's webpage.

\bibitem{rsode}
\bysame, \emph{Irregular connections and {K}ac--{M}oody root systems}, 2008,
  arXiv:0806.1050.

\bibitem{nicetalk09}
\bysame, \emph{Irregular connections, {H}itchin systems and {K}ac--{M}oody root
  systems}, 2009, Talk at conference on fundamental groups in algebraic
  geometry, University of Nice, 26 May 2009, transparencies available here:
  http://www.math.polytechnique.fr/SEDIGA/nice0905.htm.

\bibitem{quad}
\bysame, \emph{Quivers and difference {P}ainlev\'e equations}, Groups and
  symmetries: From the Neolithic Scots to John McKay, CRM Proc. Lecture Notes,
  vol.~47, AMS, 2009, 25--51, arXiv:0706.2634.

\bibitem{fission}
\bysame, \emph{Through the analytic halo: {F}ission via irregular
  singularities}, Ann. Inst. Fourier (Grenoble) \textbf{59} (2009), no.~7,
  2669--2684, Volume in honour of B. Malgrange.

\bibitem{logahoric}
\bysame, \emph{Riemann--{H}ilbert for tame complex parahoric connections},
  Transform. Groups \textbf{16} (2011), no.~1, 27--50, arXiv:1003.3177.

\bibitem{hdr}
\bysame, \emph{Habilitation memoir}, Universit\'e Paris-Sud 12/12/12,
  (arXiv:1305.6593), 2012.

\bibitem{ihptalk}
\bysame, \emph{Hyperk\"ahler manifolds and nonabelian {H}odge theory of
  (irregular) curves}, 2012, text of talk at Institut Henri Poincar\'e,
  arXiv:1203.6607.

\bibitem{slims}
\bysame, \emph{Simply-laced isomonodromy systems}, Publ. Math. I.H.E.S.
  \textbf{116} (2012), no.~1, 1--68, arXiv:1107.0874.

\bibitem{gbs}
\bysame, \emph{Geometry and braiding of {S}tokes data; {F}ission and wild
  character varieties}, Annals of Math. \textbf{179} (2014), 301--365,
  arXiv:1111.6228.

\bibitem{CB.ICM}
W.~Crawley-Boevey, \emph{{Quiver algebras, weighted projective lines, and the
  Deligne-Simpson problem}}, Madrid ICM proceedings 2006, math.RA/0604273.

\bibitem{CB-additiveDS}
\bysame, \emph{On matrices in prescribed conjugacy classes with no common
  invariant subspace and sum zero}, Duke Math. J. \textbf{118} (2003), no.~2,
  339--352.

\bibitem{CB-ihes}
\bysame, \emph{Indecomposable parabolic bundles and the existence of matrices
  in prescribed conjugacy class closures with product equal to the identity},
  Publ. Math. Inst. Hautes \'Etudes Sci. (2004), no.~100, 171--207.

\bibitem{CB-H}
W.~Crawley-Boevey and M.~P. Holland, \emph{Noncommutative deformations of
  {K}leinian singularities}, Duke Math. J. \textbf{92} (1998), no.~3, 605--635.

\bibitem{CB-Shaw}
W.~Crawley-Boevey and P.~Shaw, \emph{Multiplicative preprojective algebras,
  middle convolution and the {D}eligne-{S}impson problem}, Adv. Math.
  \textbf{201} (2006), no.~1, 180--208.

\bibitem{Del70}
P.~Deligne, \emph{\'{E}quations diff\'erentielles \`a points singuliers
  r\'eguliers}, Springer-Verlag, Berlin, 1970, Lecture Notes in Mathematics,
  Vol. 163.

\bibitem{DettReit-katz}
M.~Dettweiler and S.~Reiter, \emph{An algorithm of {K}atz and its application
  to the inverse {G}alois problem}, J. Symbolic Comput. \textbf{30} (2000),
  no.~6, 761--798.

\bibitem{eor-gda}
P.~Etingof, A.~Oblomkov, and E.~Rains, \emph{{Generalized double affine {H}ecke
  algebras of rank 1 and quantized del {P}ezzo surfaces}}, Adv. in Math.
  \textbf{212} (2007), no.~2, 749--796.

\bibitem{groechenig.hilbert.schemes}
M.~Groechenig, \emph{Hilbert schemes as moduli of {H}iggs bundles and local
  systems}, arXiv:1206.5516.

\bibitem{Harn94}
J.~Harnad, \emph{Dual isomonodromic deformations and moment maps to loop
  algebras}, Comm. Math. Phys. \textbf{166} (1994), 337--365.

\bibitem{hi-ya-nslcase}
K.~Hiroe and D.~Yamakawa, \emph{Moduli spaces of meromorphic connections and
  quiver varieties}, arXiv:1305.4092, 2013.

\bibitem{Hit-sde}
N.~J. Hitchin, \emph{The self-duality equations on a {R}iemann surface}, Proc.
  London Math. Soc. \textbf{55} (1987), no.~3, 59--126.

\bibitem{Hit95long}
\bysame, \emph{Frobenius manifolds}, Gauge Theory and Symplectic Geometry
  (J.~Hurtubise and F.~Lalonde, eds.), NATO ASI Series C: Maths \& Phys., vol.
  488, Kluwer, 1995.

\bibitem{Kac-book}
V.~G. Kac, \emph{Infinite-dimensional {L}ie algebras}, third ed., C.U.P.,
  Cambridge, 1990.

\bibitem{Katz-rls}
N.~M. Katz, \emph{Rigid local systems}, Annals of Mathematics Studies, vol.
  139, Princeton University Press, Princeton, NJ, 1996.

\bibitem{king-quivers}
A.~D. King, \emph{Moduli of representations of finite-dimensional algebras},
  Quart. J. Math. Oxford Ser. (2) \textbf{45} (1994), no.~180, 515--530.

\bibitem{kostov-DSp99}
V.~P. Kostov, \emph{On the {D}eligne-{S}impson problem}, C. R. Acad. Sci. Paris
  S\'er. I Math. \textbf{329} (1999), no.~8, 657--662.

\bibitem{kostov-nfDSp}
\bysame, \emph{Additive {D}eligne--{S}impson problem for non-{F}uchsian
  systems}, Funk. Ekvac. \textbf{53} (2010), no.~3, 395--410.

\bibitem{Kraft-Procesi-InvMath79}
H.~Kraft and C.~Procesi, \emph{Closures of conjugacy classes of matrices are
  normal}, Invent. Math. \textbf{53} (1979), no.~3, 227--247.

\bibitem{Kron.ale}
P.~B. Kronheimer, \emph{The construction of {ALE} spaces as hyper-{K}\"ahler
  quotients}, J. Differential Geom. \textbf{29} (1989), no.~3, 665--683.

\bibitem{kron-nakaj90}
P.~B. Kronheimer and H.~Nakajima, \emph{Yang-{M}ills instantons on {ALE}
  gravitational instantons}, Math. Ann. \textbf{288} (1990), no.~2, 263--307.

\bibitem{malg-book}
B.~Malgrange, \emph{\'{E}quations diff\'erentielles \`a coefficients
  polynomiaux}, Progress in Mathematics, vol.~96, Birkh\"auser Boston Inc.,
  Boston, MA, 1991.

\bibitem{nakaj-duke94}
H.~Nakajima, \emph{Instantons on {ALE} spaces, quiver varieties, and
  {K}ac-{M}oody algebras}, Duke Math. J. \textbf{76} (1994), no.~2, 365--416.

\bibitem{nakaj-duke98}
\bysame, \emph{Quiver varieties and {K}ac-{M}oody algebras}, Duke Math. J.
  \textbf{91} (1998), no.~3, 515--560.

\bibitem{nakaj-sugaku}
\bysame, \emph{Quiver varieties and quantum affine algebras [translation of
  {S}\=ugaku {\bf 52} (2000), no. 4, 337--359; {MR1802956}]}, S\=ugaku
  Expositions \textbf{19} (2006), no.~1, 53--78.

\bibitem{pr-segal}
A.~Pressley and G.~Segal, \emph{Loop groups}, Oxford Mathematical Monographs,
  The Clarendon Press Oxford University Press, New York, 1986, Oxford Science
  Publications.

\bibitem{Sab99}
C.~Sabbah, \emph{Harmonic metrics and connections with irregular
  singularities}, Ann. Inst. Fourier \textbf{49} (1999), no.~4, 1265--1291.

\bibitem{Sim-hboncc}
C.~T. Simpson, \emph{Harmonic bundles on noncompact curves}, J. Am. Math. Soc.
  \textbf{3} (1990), 713--770.

\bibitem{simpson-matrices}
\bysame, \emph{Products of matrices}, Differential geometry, global analysis,
  and topology ({H}alifax, {NS}, 1990), CMS Conf. Proc., vol.~12, Amer. Math.
  Soc., Providence, RI, 1991, pp.~157--185.

\bibitem{Simpson-katz}
\bysame, \emph{Katz's middle convolution algorithm}, Pure Appl. Math. Q.
  \textbf{5} (2009), 781--852.

\bibitem{szabo-nahm}
S.~Szab{\'o}, \emph{Nahm transform for integrable connections on the {R}iemann
  sphere}, M\'em. Soc. Math. Fr. (N.S.) (2007), no.~110, ii+114 pp.

\bibitem{vdb-doublepoisson}
M.~Van~den Bergh, \emph{Double {P}oisson algebras}, Trans. Amer. Math. Soc.
  \textbf{360} (2008), no.~11, 5711--5769, arXiv:math/0410528.

\bibitem{vdb-ncqh}
\bysame, \emph{Non-commutative quasi-{H}amiltonian spaces}, Poisson geometry in
  mathematics and physics, Contemp. Math., vol. 450, Amer. Math. Soc.,
  Providence, RI, 2008, pp.~273--299, arXiv:math/0703293.

\bibitem{yamakawa-mpa}
D.~Yamakawa, \emph{Geometry of multiplicative preprojective algebra}, Int.
  Math. Res. Pap. IMRP (2008), 77pp.

\end{thebibliography}
